\documentclass[
    12pt,
    reqno,
    final,
    ]{amsart}

\usepackage{amsmath, 
    amsthm, 
    amssymb, 
    amsthm, 
    wasysym, 
    enumitem,
    mathrsfs,
    tikz,
    bbm,
    mathtools,
    thmtools,
    xcolor,
    hyperref
}
\usepackage[nameinlink]{cleveref}
\usetikzlibrary{patterns}
\usetikzlibrary{cd}
\usepackage{stackengine}
\usepackage{float}

\newcommand{\term}[1]{\textbf{#1}}
\newcommand{\mf}[1]{\mathfrak{#1}}
\newcommand{\mc}[1]{\mathcal{#1}}
\newcommand{\ms}[1]{\mathscr{#1}}

\renewcommand{\tilde}{\widetilde}
\renewcommand{\hat}{\widehat}
\renewcommand{\bar}{\overline}

\let\set\undefined
\newcommand{\set}[1]{\left\{#1\right\}}%
\newcommand{\inner}[1]{\left\langle#1\right\rangle}
\newcommand{\gen}[1]{\left\langle#1\right\rangle}
\newcommand{\sothat}{\;\middle|\;}

\newcommand{\sse}{\subseteq}

\newcommand{\defeq}{\coloneqq}
\newcommand{\onetoone}{\leftrightarrow}
\newcommand{\into}{\hookrightarrow}
\newcommand{\onto}{\twoheadrightarrow}
\newcommand{\lto}{\xleftarrow}
\newcommand{\rto}{\xrightarrow}
\newcommand{\res}[1]{|_{#1}}
\newcommand{\map}[3]{#1 \colon #2 \to #3}
\newcommand{\Id}{\mathbbm{1}}

\DeclareMathOperator{\im}{Im}

\newcommand{\DD}{\mathbb{D}}

\newcommand{\GG}{\mathbb{G}}

\newcommand{\NN}{\mathbb{N}}

\newcommand{\QQ}{\mathbb{Q}}
\newcommand{\RR}{\mathbb{R}}
\newcommand{\TT}{\mathbb{T}}

\newcommand{\ZZ}{\mathbb{Z}}

\DeclareMathOperator{\Hom}{Hom}

\newcommand{\actson}{\curvearrowright}

\newcommand{\dirsum}{\oplus}
\newcommand{\Dirsum}{\bigoplus}
\newcommand{\tensor}{\otimes}

\newcommand{\isom}{\cong}

\newcommand{\nsphere}[1][{}]{\mathbf{S}^{#1}}

\DeclareMathOperator{\diam}{diam}

\DeclareMathOperator{\dvol}{dvol}
\DeclareMathOperator{\vol}{vol}

\DeclareMathOperator{\Lip}{Lip}

\newcommand{\skel}[2]{{#1}^{(#2)}}

\newcommand{\Forms}[2][*]{\Omega^{#1}(#2)}

\newcommand{\boundary}{\partial}
\newcommand{\coboundary}{\delta}
\newcommand{\homotopic}{\simeq}
\newcommand{\rel}[1]{\operatorname{rel} #1}
\newcommand{\wedgesum}{\vee}

\newcommand{\trtdt}{\tensor \RR\gen{t,dt}}
\newcommand{\realize}[1]{\langle #1 \rangle}
\newcommand{\minimal}[1]{\mc{M}^*_{#1}}

\newenvironment{diagram}{\begin{center}\begin{tikzcd}[ampersand replacement=\&]}{\end{tikzcd}\end{center}}

\newcommand{\aisom}{\cong}
\newcommand{\simplex}[1][{}]{\Delta^{#1}}

\newcommand{\Dil}[2][{}]{\operatorname{Dil}_{#1}\l(#2\r)}
\newcommand{\rtdt}{\RR\gen{t,dt}}

\newcommand{\cochains}[2]{C^{#2}(#1)}

\newcommand{\isolipverbose}[3]{\operatorname{V\delta}^{#1}_{#2}(#3)}

\newcommand{\isolip}[1]{\operatorname{V\delta}(#1)}
\newcommand{\filvol}[1]{\operatorname{FV}(#1)}

\newcommand{\lcsaverbose}[3]{#1^{#2}_{#3}}
\newcommand{\lcsa}[2]{#1_{#2}}

\newcommand{\pinstep}[2]{\bar{\pi_{#1}}^{(#2)}}

\newcommand{\pidual}{V}
\newcommand{\picomplement}{\Lambda}
\newcommand{\vsdeg}[2]{{#1}_{#2}}
\newcommand{\vsstep}[2]{{#1}(#2)}

\newcommand{\pindual}[1]{\vsdeg{\pidual}{#1}}
\newcommand{\pindualverbose}[2]{\pindual{#1}(#2)}
\newcommand{\pindualupto}[2]{\vsdeg{\pidual}{#1}^{\le #2}}
\newcommand{\pindualstepverbose}[3]{\vsdeg{\picomplement}{#1}^{(#2)}(#3)}
\newcommand{\pindualstep}[2]{\vsdeg{\picomplement}{#1}^{(#2)}}

\newcommand{\mapdeg}[2]{#1_{#2}}
\newcommand{\mapstep}[2]{#1^{(#2)}}
\newcommand{\mapdegstep}[3]{\mapstep{\mapdeg{#1}{#2}}{#3}}

\newcommand{\postspace}[2]{{#1}_{#2}}
\newcommand{\postspacestep}[3]{\postspace{#1}{#2}^{#3}}
\newcommand{\postmap}[1]{\mapdeg{p}{#1}}
\newcommand{\postmapstep}[2]{\mapdegstep{p}{#1}{#2}}

\newcommand{\adeg}[2]{{#1}^{#2}}
\newcommand{\aplus}[1]{\adeg{#1}{+}}

\newcommand{\indecomp}[1]{Q{#1}}

\newcommand{\aprod}[1]{\wedge^{#1}}

\newcommand{\sulext}[2]{#1\gen{#2}}

\newcommand{\wt}{\mathrm{wt}}

\newcommand{\cautious}[2]{W_{#1}^{#2}}
\newcommand{\cautiousstep}[2]{M_{#1}^{#2}}
\newcommand{\naive}[2]{\tilde{W}_{#1}^{#2}}
\newcommand{\naivestep}[2]{\tilde{M}_{#1}^{#2}}

\newcommand{\len}[1]{\le_{#1}}
\newcommand{\lessn}[1]{<_{#1}}

\newcommand{\CS}{\cautiousstep}
\newcommand{\NS}{\naivestep}
\newcommand{\scomp}{d_{sim}}
\newcommand{\abcomp}{\delta}
\newcommand{\ncomp}{d_{nil}}
\newcommand{\ncompstep}[1]{\ncomp^{#1}}

\newcommand{\torsor}[1]{[{#1}]}

\newcommand{\diag}[1]{\mc{D}_{#1}}

\tikzset{commutative diagrams/.cd, mysymbol/.style={start anchor=center, end anchor=center,draw=none}}
\newcommand{\commutes}{\arrow[mysymbol]{}[description]{\text{\rotatebox{180}{\Large$\circlearrowright$}}}}
\newcommand{\homotopicCommuteSymbol}{\stackrel{\text{\rotatebox{180}{\Large$\circlearrowright$}}}{\text{\Large$\homotopic$}}}

\renewcommand{\SS}{\mathbf{S}}
\renewcommand{\DD}{\mathbf{D}}

\newcommand{\sspan}[1]{\operatorname{span}\set{#1}}

\newcommand{\quasiisom}{\simeq}

\newcommand{\openstar}[1]{\bigstar(#1)}

\newcommand{\flip}[1]{\rotatebox[origin=c]{180}{#1}}
\renewcommand{\phi}{\varphi}

\newcommand{\stackmap}[2]{\stackanchor{\flip{$#1$}}{$#2$}}

\newcommand{\basis}[1]{\overline{\mathbf{#1}}}
\newcommand{\inde}[1]{\mathbf{#1}}

\newcommand{\pp}{\oplus}

\newcommand{\aint}{\textstyle\int}
\renewcommand{\l}{\left}
\renewcommand{\r}{\right}

\newcommand{\wedgeofs}[1][{c}]{X_{#1}}
\newcommand{\torusofsplit}[1]{\l(\TT_{\splitt}\r)_{#1}}

\newcommand{\splitt}{\mc{T}}
\newcommand{\gluedt}{\TT_{\splitt_{A} \wedgesum \splitt_{B}}}

\newcommand{\highdeg}[2]{\zeta_{#1,#2}}
\newcommand{\highdegrep}[3]{\highdeg{#1}{#2}^{(#3)}}
\numberwithin{equation}{section}

\declaretheorem[style=plain,sibling=equation]{theorem,lemma,proposition,corollary,conjecture,question}

\declaretheorem[style=definition,sibling=equation]{definition}

\declaretheorem[style=plain,
within=equation]{claim}

\declaretheoremstyle[
    spaceabove=\parskip, 
    spacebelow=\parskip, 
    headfont=\itshape, 
    bodyfont =\normalfont,
    qed=$\lozenge$, 
    headpunct={:}]{proofofclaimstyle}
    
\declaretheorem[
    name={Proof of Claim}, 
    style=proofofclaimstyle, 
    unnumbered]{proofofclaim}

\declaretheorem[style=definition,numbered=no]{note,notation}

\crefname{maintheorem}{theorem}{theorems}
\crefname{conjecture}{conjecture}{conjectures}
\crefname{question}{question}{questions}
\crefname{claim}{claim}{claims}

\usepackage[top             =2cm, 
            bottom          =4.5cm, 
            left            =2.5cm, 
            right           =2.5cm, 
            marginparwidth  =2.0cm,
            ]{geometry}

\hypersetup{
    colorlinks=true, 
    linktoc=all,     
    urlcolor=blue,
    linkcolor=blue,  
    citecolor=red,
}

\usepackage[
    colorinlistoftodos,
    prependcaption,
    color=red!50,
    ]{todonotes}
\setuptodonotes{inline}

\title{Volumes of Nullhomotopies in Nilpotent Spaces}
\author{Kyle Hansen}

\begin{document}
\markboth{\theauthor}{\thetitle}

%-----------------------
% Front Matter
%-----------------------
\begin{abstract}  
The Shadowing Principle of Manin has proved a valuable tool for addressing questions of quantitative topology raised by Gromov in the late 1900s.
The principle informally provides a way for bounded algebraic maps between differential graded algebras to be translated into nearby genuine maps between their geometric realizations.
We extend this principle to finite towers of principal $K(G,n)$ fibrations, and in particular apply this construction to nilpotent spaces.
As a specific application of the extended principle, we provide upper bounds on the asymptotic behavior of volumes of nullhomotopies of Lipschitz maps into nilpotent spaces.
We further refine these bounds in the case when $c = 1$ to nearly meet those of the simply connected setting.
We similarly refine these bounds in the event the target space is coformal, and demonstrate that the bounds in this setting are nearly sharp.
\end{abstract}

\maketitle

%-----------------------
% Main Document
%-----------------------
\section{Introduction}\label{sec:introduction}
In \cite{gromov1978homotopical}, Gromov demonstrates for reasonable spaces that there is a polynomial in $L$ bounding the number elements of a homotopy class that admit $L$-Lipschitz representatives.
That is, he suggested that topological invariants be studied not only algebraically, but quantitatively as well.
Whether two manifolds are cobordant, or two maps are homotopic, are questions which in theory carry some geometric data along with them. 
Historically such equivalences have been dealt with in mostly algebraic ways, ignoring issues of computability or the potential sizes of their geometric realizations.
However, as Gromov has said in \cite{gromov1999quantitative},
\begin{quotation}
    ...we can organize an innocuous enough metric on $\SS^3$ so that it takes more than $10^{30}$ years to contract certain loops in the sphere and in the course of contraction we need to stretch the loop to something like $10^{30}$ light years in size. So if $10^{30}$ years is all the time you have, you conclude that the loop is not contractible and whether or not $\pi_1(\SS^3) = 0$ becomes a matter of opinion. 
\end{quotation}

Together with Gromov, the work of Nabutovsky and Weinberger (such as \cite{nabutovsky2000variational,nabutovsky2003area}) paved the way for further quantitative study of topological spaces. 
The past decade of work (e.g., \cite{manin2015asymptotic,manin2016volume,guth2017volumes,Chambers_DMW_2018_nullcobordism,Chambers_MW_2018_quantitative_nullhomotopy,manin2020integral}) by Chambers, Guth, Manin, and Weinberger, among others have developed Gromov's initial inquiries into a fruitful field of study, turning quantitative topology from an idea into a practical theory.

Sullivan's localization methods in homotopy theory have repeatedly proved to be an invaluable tool for understanding algebraic invariants of topological spaces, and proves useful in Gromov's suggested program as well.
As long as we only care about the rational homotopy type of a space (i.e., if we start by rationalizing our space, or localizing our space at $\set{0}$), standard algebraic invariants for a large class of spaces can be completely encoded by Sullivan's model.
This includes simply connected spaces, but extends more generally to \term{nilpotent spaces}, a class of spaces which show up naturally in the study of almost nonnegatively curved manifolds (see \cite{kapovitch2010nilpotency}).
Since Sullivan models completely encode the rational homotopy invariants of nilpotent spaces, such spaces make ideal objects of study in the story of quantitative topology.

One may naturally view quantitative topology as a cousin of the more familiar field of geometric group theory.
For example, recall that the \term{Dehn function} $\delta_{G}$ of a finitely presented group $G$ measures the size of a disk needed to fill a loop of length $\ell$ (i.e., the geometric realization of a word $w$ of wordlength $\ell$) in the Cayley graph of $G$.
If $G$ is additionally nilpotent then Malcev's correspondence places $G$ inside a Lie group $\GG$ as a cocompact lattice.
In particular, $G$ is quasi-isometric to $\GG$, and one may study the asymptotics of $\delta_{G}$ in $\GG$ instead.
Alternatively, viewing $\GG$ as the universal cover of $K(G,1)$, we view $\delta_{G}$ as an indicator of topological complexity.
The additional geometric data in this setting allows one to define various high-dimensional analogues of Dehn functions, all of which are isoperimetric in spirit.
That is, we ask how large an $n$-dimensional filling must be among all $(n-1)$-dimensional objects of a given size which admit a filling, where a ``filling'' may be defined (for example, and not exclusively) by the existence of a nullhomotopy.
See \cite{young2011homological,young2013filling} as well as \cite{brady2008snowflake,gromov1983filling,white1984mappings} for various notions of (not necessarily equivalent) higher order Dehn functions.

In this paper, we are primarily interested in a question of quantitative topology which can also be seen as an investigation of a certain higher order Dehn function.
\begin{question}\label{question:quantitative nullhomotopy}
    Let $Y$ be a closed Riemannian manifold.
    If $f : \SS^n \to Y$ is a nullhomotopic, $L$-Lipschitz map, what is the volume of an efficient nullhomotopy in terms of $L$? How does this depend on the topology of $Y$?
\end{question}
In the case $Y = K(G,1)$ with $G$ a finitely generated nilpotent group, this question can be seen as asking for the volume of an optimal (homotopic) filling of an $L$-Lipschitz map $F : \SS^n \to \tilde{Y} = \GG$, which in turn can be seen as a higher order Dehn function of $G$.
The current paper seeks to address this question for $Y$ a nilpotent space more generally.

In \cite{Chambers_DMW_2018_nullcobordism} the answer to the above question is shown to depend in general on the rational homotopy type of the target space $Y$.
In deeper conversation with a conjecture of Gromov's, Manin has made a study of this question when $Y$ is simply connected, and shows the following:
\begin{theorem}[{\cite[Theorem A(ii)]{Manin_2019}}]\label{thm:simply connected upper bounds}
    If $Y$ is simply connected and $\map{f}{\SS^{n}}{Y}$ is $L$-Lipschitz, then one can always find a nullhomotopy whose volume is at most $O(L^{2n})$.
\end{theorem}
In fact, the full result is a strengthening of the relevant conjecture of Gromov's; see \cite[Conjecture C; Theorem A]{Manin_2019}.
A family of examples shows that this bound is, in some sense, sharp:
\begin{theorem}[{\cite[\textsection 7]{Chambers_MW_2018_quantitative_nullhomotopy}}]\label{thm:simply connected lower bounds}
    Let $n \ge 2$. There is a simply connected space $Y_{n}$ and an $L$-Lipschitz map $\map{f}{\SS^{n}}{Y_{n}}$ such that any nullhomotopy of $f$ has volume $\Omega(L^{2(n-1)})$.
\end{theorem}

The proof of \Cref{thm:simply connected upper bounds} utilizes a powerful tool, developed by Manin for simply connected spaces in \cite{Manin_2019}, called the \term{Shadowing Principle}.
One consequence of the principle is that it provides a method for homotoping a map which is algebraically bounded (i.e., whose pullback is bounded on the level of DGAs) into a bounded geometric map.
Somewhat more concretely, suppose a genuine map $f : X \to Y$ into a simply connected space has an algebraic model $f^*m_Y : \minimal{Y} \to \Forms{X}$ which lies in the same homotopy class as a map whose dilatation is bounded by $L$. 
(One may think of the dilatation as an algebraic analogue of the Lipschitz constant.)
The principle then demonstrates that $f$ itself lies in the homotopy class of a genuine map whose Lipschitz constant is bounded (asymptotically) by $L$.

The guiding heuristic of the principle is that it is often easier to create bounded maps on the level of algebraic information, and that this is nearly enough to produce bounded geometric maps;
if you have an (unbounded) geometric map which has a bounded algebraic counterpart, then there is a bounded geometric map homotopic to the original one. 
Moreover, the algebraic counterpart to this new map is not too far away from the bounded algebraic counterpart of the original.
We encourage the reader to refer to the original paper \cite{Manin_2019} for further exposition and intuition for this principle, as well as its applications and the idea behind the principle's name.

One limitation of this principle is that one must manage to find bounded algebraic maps in the first place, to which one may then apply the principle.
This task may not be any easier than finding their geometric counterparts.
For instance, in \cite[\textsection 5.3]{Manin_2019}, a naive application of the principle produces nullhomotopies of volume $O(L^{n+2})$ for maps $\SS^{m} \to \SS^{2n}$ as long as $m \ge 2n - 1$.
However, \cite{berdnikov2022lipschitz} shows that it is possible to find linear nullhomotopies for maps $\SS^3 \to \SS^2$, and \cite[Theorem A(iii)]{berdnikov2022scalable} implies that this linear bound is always possible for maps into spheres.

\subsection{Summary of Results}
The Shadowing Principle in \cite{Manin_2019} requires the target spaces to be simply connected.
The main contribution of this paper is the extension of the Shadowing Principle to certain (towers of) principal fibrations (\Cref{thm:lifting shadows}).
Concretely, this results in a generalization of the principle to nilpotent targets (\Cref{thm:nilpotent shadows}).
As an application, we produce upper bounds on the volumes of nullhomotopies in general nilpotent spaces (\Cref{thm:nilpotent upper bounds}).
The upper bounds we attain by this method are not yet known to be sharp, and we expect they can be refined.
A more careful analysis on simple spaces, for instance, reduces the expected upper bound somewhat significantly (\Cref{cor:simple isolip upper bounds}).
Thus far we have been unable to find a family of example spaces which come close to the general upper bounds.
However, restricting our attention to spaces whose differential is quadratic (i.e., coformal spaces), we are able to futher refine the upper bounds (\Cref{cor:coformal isolip upper bounds}) and demonstrate a family of examples which nearly meet these bounds (\Cref{thm:coformal example}).

To be more concrete, we briefly introduce some notation.
First, following \cite{gromov1999metric}, we introduce the idea of the volume of a Lipschitz map $g : X \to Y$ between compact (piecewise) smooth Riemannian manifolds with boundary.
After embedding $Y \into \RR^{N}$ in a (piecewise) smooth manner, Rademacher's Theorem implies that $g$ is differentiable almost everywhere (since $g$ is Lipschitz), and hence one may define the \textbf{volume} of $g$ by
\[
    \vol{g} \defeq \int_{X}|\operatorname{Jac}(g)|\dvol{}.
\]

Given $f : X \to Y$ an $L$-Lipschitz, nullhomotopic function we define the \term{(homotopic) filling volume of $f$} to be
\begin{equation*}
    \filvol{f} \defeq \inf \set{\vol{H} \sothat H : X \times I \to Y,\, H\res{X \times \set{0}} = f,\, H\res{X \times \set{1}} = *}
\end{equation*}
where $*$ is the constant map to the basepoint of $Y$ and $I = [0,1]$ is the unit interval, and where the infimum is taken with respect to Lipschitz maps $H$ so that the notion of volume is well-defined.
The term ``filling volume'' comes from when $X = \SS^n$ and we view the nullhomotopy $H$ as a map $H : \DD^{n+1} \to Y$.
We may then define a higher order Dehn function $\isolipverbose{X}{Y}{L}$ as the maximal filling volume among $L$-Lipschitz nullhomotopic maps $Y \to X$. 
That is,
\begin{equation*}
    \isolipverbose{X}{Y}{L} \defeq \sup \set{\filvol{f} \sothat f : X \to Y,\, \Lip{f} \le L,\, [f] = [*] \in [X,Y]}.
\end{equation*}
In the event $X = \SS^n$ we write $\isolipverbose{n}{Y}{L} \defeq \isolipverbose{\SS^n}{Y}{L}$.
We summarize our results in \Cref{tab:results} below for comparison.

\begin{table}[htb!]
    \begin{center}
    \renewcommand{\arraystretch}{2}
    \begin{tabular}[htb!]{c||c|c|c|c}
        \textbf{\textit{Type of Space}} $Y$
            & \shortstack[c]{\textit{Simply} \\ \textit{Connected}}
            & \textit{Simple}
            & \shortstack[c]{\textit{$c$-Step} \\ \textit{Nilpotent}}
            & \shortstack[c]{\textit{Coformal} \\ \textit{+ $c$-Step} \\ \textit{Nilpotent}}
            \\ \hline \hline
            &&&\\[-1.5em]
        \shortstack[c]{\textbf{Lower Bounds} \\ \textbf{on} $\isolipverbose{n}{Y}{L}$}
            & $\Omega(L^{2(n-1)})$ 
            & ?
            & ?
            & $\Omega(L^{(c-1)(n-1)})$
            \\
        \textbf{Reference}
            & {\cite[\textsection 7]{Chambers_DMW_2018_nullcobordism}} 
            & \text{}
            & \text{}
            & \Cref{thm:coformal example}
            \\ \hline
            &&&\\[-1.5em]
        \shortstack[c]{\textbf{Upper Bounds} \\ \textbf{on} $\isolipverbose{n}{Y}{L}$}
            & $O(L^{2n})$ 
            & $O(L^{2n+1})$
            & $O(L^{(4c-1)n})$
            & $O(L^{(c-1)n})$
            \\
        \textbf{Reference}
            & {\cite[Theorem A]{Manin_2019}} 
            & \Cref{cor:simple isolip upper bounds}
            & \Cref{cor:nilpotent isolip upper bounds}
            & \Cref{cor:coformal isolip upper bounds}
    \end{tabular}
    \caption{Upper Bounds on $\isolip{L}$}
    \label{tab:results}
    \end{center}
\end{table}

We expect the general $c$-step nilpotent upper bounds can be refined:
\begin{conjecture}\label{conj:better nilpotent bounds}
    If $Y$ is a $c$-step nilpotent space then $\isolipverbose{n}{Y}{L} = O(L^{3nc})$.
\end{conjecture}

\subsection{Structure of Paper}
In \Cref{sec:background} we orient the reader with the notation, conventions, and background concerning fibrations, nilpotent spaces, minimal Sullivan extensions/models, and the category equivalences induced by the localization and realization functors. 
In \Cref{sec:shadows} we generalize the Shadowing Principle to principal $K(G,n)$ fibrations, which specifically results in a corresponding principle for nilpotent spaces.
As an application of the principle, in \Cref{sec:upper bounds} we produce upper bounds on the volumes of nullhomotopies filling in most of the bottom row of \Cref{tab:results}.
Finally, we demonstrate the (near) sharpness of the results for coformal spaces in \Cref{sec:coformal-example}.

\subsection{Acknowledgements}
The work in this paper was done as part of my doctoral thesis at the University of California, Santa Barbara; consequently there are too many people to name who made this work possible.
Principally, I thank Fedya Manin for drawing my attention to many of these questions in the first place.
His constant guidance and inspiration have allowed me to discover the answers herein.
Furthermore, I am deeply indebted to Daryl Cooper and Guofang Wei for their thorough reading of early drafts of this work. 
Their questions and advice have led to clearer notation, more precise statements, and more accurate proofs.
All errors that remain---fatal or otherwise---are solely my own contribution.

\section{Background and Preliminaries}\label{sec:background}
We provide here our notation and conventions for nilpotent spaces, as well as Sullivan DGAs.
All spaces throughout the paper, unless otherwise stated, are assumed to be compact Riemannian manifolds or finite simplicial complexes.
Let $Y$ be such a space, and recall that there is an action $\pi_1(Y) \actson \pi_n(Y)$ for all $n \ge 1$ induced by the Whitehead bracket.
Namely, we have $\gamma \bullet \beta \defeq [\gamma,\beta] + \beta$ where $[\cdot,\cdot] : \pi_{1}(Y) \tensor \pi_{n}(Y) \to \pi_{n}(Y)$ is the Whitehead product.

Before providing all relevant definitions, we note here that we are primarily interested in spaces up to \term{rational equivalence} which we define a little later on.
Suffice to say, we only focus on the non-torsion aspects of algebraic invariants when possible. 
For instance we might identify $\SS^{2n+1}$ and $K(\ZZ,2n+1)$.
We will be primarily interested in spaces which admit triangulations, and the algebra of their differential forms, conflating the notation in this setting.
We let $\Forms{Y}$ denote the smooth differential forms when $Y$ is a Riemannian manifold, as well as the algebra of piecewise smooth forms when $Y$ is a simplicial complex.

\subsection{Nilpotent Spaces}\label{ssec:nilpotent spaces}
Let $G$ be a nilpotent group acting on $M$, an abelian group.
Recall that an action $G \actson M$ is \term{nilpotent} if the lower central series of the action terminates.
That is, $\lcsa{M}{c+1} \aisom \set{0}$ and $\lcsa{M}{c} \not\aisom \set{0}$ for some $c$, where
    \begin{equation*}
        \lcsa{M}{j} = \lcsaverbose{M}{G}{j} \defeq \begin{cases}
            M & j = 1, \\
            \gen{g \bullet m - m \sothat g \in G, m \in \lcsaverbose{M}{G}{j-1}} & j > 1.
        \end{cases}
    \end{equation*}
    
    \begin{definition}
        If $\pi_1(Y)$ is nilpotent and acts nilpotently on $\pi_n(Y)$ for all $n$ we say that $Y$ is \term{nilpotent}. 
        If $\pi_1(Y)$ is $c$-step nilpotent and the action is $c$-step nilpotent for all $n$, we say that $Y$ is \term{$c$-step nilpotent}.
        If $c = 1$, then we say that $Y$ is \term{simple} or \term{abelian}.
        If one of the above properties holds for all $n \le N$ for some $N$ we say that the property holds \term{up through degree $N$}.
    \end{definition}
    Clearly, simply connected spaces are simple, and simple spaces are nilpotent.
    Note that nilpotent spaces don't always enjoy the same properties as simply connected spaces.
    This is especially relevant for questions of how to create new nilpotent spaces out of old ones.
    For instance, nilpotent spaces are closed under direct products; however, nilpotent spaces are \textbf{not} closed under wedge sums (a property which \textbf{is} enjoyed by simply connected spaces). 
    One can see nearly immediately that $\nsphere[1] \wedgesum \nsphere[1]$ is not nilpotent, even though both factors are. 
    
    Rather than constructing examples of nilpotent spaces using products or wedges, we instead build nilpotent spaces up inductively. This is, nilpotent spaces may be defined as those spaces whose Postnikov tower admits a principal refinement; in particular, a nilpotent space may be constructed as the inverse limit of principal $K(G,n)$ fibrations.

    \begin{definition}
        Let $K = K(G,n)$ be an Eilenberg-MacLane space with $G$ abelian.
        Define the space $K^{+} = K(G,n+1)$ with fixed basepoint $* \in K^{+}$.
        Recall that the \term{pathspace fibration} is the obvious fibration $PK^{+} \to K^{+}$ with fiber $p^{-1}(*) = \Omega K^{+} \homotopic K$.
        Now let $\map{k}{B}{K^{+}}$ for some topological space $B$.
        The \term{pullback of $PK^{+}$ by $k$} is the fibration $Y = B \times_{k} PK^{+} \to B$ (unique up to weak homotopy equivalence) under the pullback 
        \begin{diagram}
            Y 
                \arrow[dr, phantom, "\scalebox{1.5}{$\lrcorner$}" , very near start, color=black]
                \rar{}
                \dar{p}
            \& PK^{+}
                \dar{}
            \\
            B 
                \rar{k}
            \& K^{+}
        \end{diagram}
        Such a fibration $\map{p}{Y}{B}$ has fiber $K$ and is called a \term{principal $K(G,n)$ fibration}.
        In particular, a fibration $\map{p}{Y}{B}$ is a principal $K(G,n)$ fibration if there is some map $\map{k}{B}{K^{+}}$ for which $p$ is the pullback.
    \end{definition}

    \begin{proposition}
        Given a principal $K(G,n)$-fibration $Y \to B$ induced by a map $\map{k}{B}{K(G,n+1)}$, there is a class $[k] \in H^{n+1}(B;G)$ which classifies $\map{p}{Y}{B}$ up to weak homotopy equivalence.
        That is there are bijections
        \begin{equation*}
            \set{
                \renewcommand{\arraystretch}{1}
                \begin{array}{c} 
                    \text{principal} \\ \text{$K(G,n)$ fibrations}
                \end{array}}
                \onetoone
                \set{
                    \renewcommand{\arraystretch}{1}
                    \begin{array}{c}
                    \text{homotopy classes of} \\ \text{maps $B \to K(G,n)$}
                \end{array}}
                \onetoone
                H^{n+1}(B;G).
        \end{equation*}
        \begin{proof}
            $Y$ is up to weak equivalence the homotopy fiber of $k$, and therefore it is enough to show that the choice of map $k$ is equivalent to a class $[k] \in H^{n+1}(B;G)$. 
            But by the universal coefficient theorem (or rather, by an interpretation of the theorem) the homotopy class of $k : B \to K(G,n+1)$ is precisely determined by $[k] \in H^{n+1}(B;G)$.
        \end{proof}
    \end{proposition}

    \begin{definition}
        The class $[k] \in H^{n+1}(B;G)$ is called the \term{$k$-invariant} of the fibration.
        We call $\map{k}{B}{K(G,n+1)}$ itself the \term{classifying map} of the fibration.
    \end{definition}

    Fix a principal $K(G,n)$ fibration $Y \rto{p} B$ of simplicial complexes $Y$ and $B$.
    \begin{definition}
        Let the following diagram commute:
        \begin{diagram}
            A \dar[hook,swap]{\iota} \rar{\alpha} 
                \& Y \dar[two heads]{p} \\
            X \rar[swap]{\beta}
                \& B
        \end{diagram}
        Let $k : B \to K(G,n+1)$ be the classifying map of the fibration.
        Let $CX$ and $CA$ denote the cones of $X$ and $A$ respectively.
        The map $\alpha$ induces a nullhomotopic map $A \to K(G,n+1)$, and the \term{obstruction class} for this diagram is the class $\mc{O} \in H^{n+1}(X,A;G)$ induced by the composition 
        \[
            X \cup CA \to B \rto{k} K(G,n+1).
        \]
    \end{definition}

    \begin{proposition}[{\cite[Proposition 4.72]{hatcher2002algebraic}}]\label{prop:obstruction to relative extension}
        A map $X \to Y$ fitting into the diagonal of the above diagram exists if and only if the obstruction of the diagram vanishes.
    \end{proposition}

    In a similar manner (see \Cref{prop:topological obstruction} below), the obstruction to a relative homotopy between two potential fillings defines an element of $H^n(X,A;G)$. 
        Let $\ms{D} = \ms{D}(X,A,Y,B)$ be the set of all diagrams of the above form whose obstruction class vanishes.
    We will often be interested in the situation $(X,A) = (\simplex, \boundary \simplex)$, the standard simplex.
    When $X = \simplex$, we will assume that $A = \boundary \simplex$, often without specifying so.
    \begin{definition}
        The \term{external diagram} of a map $\map{f}{X}{Y}$ is the diagram $\diag{f}$ given by 
        \begin{diagram}
            A \rar{f \circ \iota} \dar[swap,hook]{\iota}
                \& Y \dar[two heads]{p} \\
            X \rar[swap]{p \circ f} 
                \& B
        \end{diagram}
        By \Cref{prop:obstruction to relative extension} we have $\diag{f} \in \ms{D}$ since $f$ itself provides the diagonal map of interest.
        If $D_{f} = D_{g}$ for two maps $\map{f,g}{X}{Y}$ we call $f$ and $g$ \term{externally equivalent}.
        Moreover, suppose that $f$ and $g$ are externally equivalent via a common diagram 
        \begin{diagram}
            A \arrow[d,swap,hook,"\iota"] \arrow[r,"\alpha"]
                \& Y \dar[two heads]{p}
                \\
            X \arrow[ur, bend right=-7,"f"]\arrow[ur,bend right=7,swap,"g"] \arrow[r,swap,"\beta"]
                \& B
        \end{diagram}
        Suppose that there is a homotopy $H$ from $f$ to $g$ for which $H\res{t}$ is externally equivalent to $f$ (hence to $g$) for all $t \in I$.
        Then $f$ and $g$ are said to be \term{internally homotopic through $H$} and we call $H$ an \term{internal homotopy} from $f$ to $g$.
        If $f$ and $g$ are internally homotopic through some (unspecified) homotopy, we say that $f$ and $g$ are \term{internally homotopic} or \term{internally equivalent}.
    \end{definition}
    Under this terminology we have:
    \begin{proposition}[{\cite[Exercise 4.24]{hatcher2002algebraic}}]\label{prop:topological obstruction}
        Let $f,g: X \to Y$ be externally equivalent maps.
        The obstruction to an internal homotopy from $f$ to $g$ is a class in $H^{n}(X,A;G)$.
    \end{proposition}

    \begin{corollary}\label{cor:homotopy torsor}
        Internal homotopy classes form a torsor for $H^{n}(X,A;G)$.
        In particular, the collection of internal homotopy classes for $(\simplex[n],\boundary \simplex[n])$ is a torsor for $G \aisom H^n(\simplex[n],\boundary \simplex[n];G)$. 
    \end{corollary}

    Recall that a Postnikov system for $Y$ is an inverse system 
    \begin{equation*}
        \cdots \to \postspace{Y}{n} \rto{\postmap{n}} \postspace{Y}{n-1} \to \cdots
    \end{equation*}
    where each $\postmap{n}$ is a $K(\pi_n(Y),n)$ fibration.
    Moreover, the inclusion $Y \to Y_n$ induces an isomorphism $\pi_k(Y) \aisom \pi_k(\postspace{Y}{n})$ for all $k \le n$.
    In general this fibration is \textit{not} principal, but when $Y$ is nilpotent, $p_n$ decomposes into a \term{principal refinement} of principal $K(G,n)$ fibrations.
    That is, if the action of $\pi_1(Y)$ on $\pi_n(Y)$ is $c$-step nilpotent, \cite[Thm. I.2.9]{hilton2016localization} shows that there are abelian groups $\pinstep{n}{j}$ and a sequence of principal $K(\pinstep{n}{j},n)$ fibrations $\postmapstep{n}{j} : \postspacestep{Y}{n}{j} \to \postspacestep{Y}{n}{j-1}$ with $\postspacestep{Y}{n}{c} = \postspacestep{Y}{n}{}$ and 
    $\postmap{n} = \postmapstep{n}{1} \circ \postmapstep{n}{2} \circ \cdots \circ \postmapstep{n}{c}$.
    In the event that $Y$ is a Riemannian manifold, the construction of the Postnikov tower in general creates a fibration of CW complexes whose fibration maps $\postmapstep{n}{j}$ are \textit{a priori} known only to be continuous.
    More specifically, there is no guarantee we have any reasonable notion of ``local trivializations'' around each fiber as one might hope.
    However, one may replace $Y$ with a rationally equivalent simplicial complex whose geometric properties are not far from those of $Y$.
    Working in the simplicial category instead allows one to produce fibrations in the simplicial category; in particular, the maps $\postmapstep{n}{j}$ may be assumed to be simplicial maps (see for example \cite[Theorem 11.4]{goerss2009simplicial} and \cite[\textsection 1.6]{felix2015rational}).
    
\subsection{Sullivan Models}
    We now turn to producing algebraic models of the above ideas.
    This will allow use to produce a more broad collection of examples of nilpotent spaces fairly easily, without needing to describe the Postnikov tower of a space quite so tediously.

    Much of this chapter is a conglomeration of definitions and results from \cite{felix2012rational,felix2015rational,griffiths2013rational}.
    We again start by providing the relevant definitions and notation for Sullivan models, and then describe when it is possible to build up such models as elementary extensions.
    We conclude with a few examples of such models.

    For our purposes, a \term{commutative differential graded algebra} (a \term{DGA}) is a cochain complex $(\mc{A},d)$ of vector spaces $\mc{A} = \set{\adeg{\mc{A}}{n}}$ over $\QQ$ or $\RR$, together with a graded commutative multiplication which satisfies the (graded) Leibniz rule.
    We assume from here on out that every DGA $\mc{A}$ is \term{connected}, meaning that $H^0(\mc{A}) = \QQ$.
    We define $\aplus{\mc{A}} \defeq \set{\adeg{\mc{A}}{n}}_{n > 0}$.

    The alternating algebra $\wedge Z$ on a graded vector space $Z = \set{\vsdeg{Z}{n}}_{n \ge 1}$ is a graded algebra, which can be made into a DGA with an appropriate choice of differential (a trivial example being $d = 0$).
    We define 
    \begin{equation*}
        \aprod{k} Z = \sspan{z_1 \wedge \cdots \wedge z_k \sothat z_i \in Z}
    \end{equation*}
    and define $\aprod{\ge k} Z$ and $\aprod{\le k} Z$ and $\aprod{> k} Z$ and $\aprod{< k} Z$ in the obvious ways.

    When $\mc{A} = \wedge Z$ the \term{space of indecomposables} is the graded vector space
    \begin{equation*}
        \indecomp{\mc{A}} = \frac{\aplus{\mc{A}}}{\aplus{\mc{A}} \wedge \aplus{\mc{A}}} \aisom Z
    \end{equation*}
    and in degree $n$ the map $Z \to \wedge Z$ inducdes an identification 
    \begin{equation*}
        \adeg{(\indecomp{\mc{A}})}{n} \aisom \vsdeg{Z}{n}.
    \end{equation*}

    Define a commutative cochain algebra morphism $(\mc{A},d) \to (\mc{A} \tensor \wedge Z,d)$ by ${a \mapsto a \tensor 1}$.
    We identify $\mc{A} = \mc{A} \tensor 1$ and $Z = 1 \tensor Z$.
    This morphism is a \term{Sullivan extension} if
    \begin{enumerate}
        \item $H^{0}(\mc{A},d) = \QQ$ (that is, $\mc{A}$ is connected),
        \item $Z = \set{\vsdeg{Z}{n}}_{n \ge 1}$ is concentrated in positive degrees, and
        \item $Z$ satisfies the \term{nilpotence condition}. 
        That is, $Z$ is the union of an increasing family of subspaces 
        \begin{equation*}
            0 = \vsstep{Z}{0} \sse \cdots \sse \vsstep{Z}{r} \sse \cdots    
        \end{equation*}
        where the differential satisfies
        \begin{equation*}
            \im{d} \res{\vsstep{Z}{r}} \sse \begin{cases}
                \mc{A} & r = 1 \\
                \mc{A} \tensor \wedge \vsstep{Z}{r-1} & r > 1.
            \end{cases}
        \end{equation*}
    \end{enumerate}
    The extension is \term{minimal} if in addition
    \begin{enumerate}[resume]
        \item $\im{d} \sse \aplus{\mc{A}} \tensor \wedge Z + \mc{A} \tensor \aprod{\ge 2} Z$.
    \end{enumerate} 
    The \term{fiber} of the extension is the DGA $(\wedge Z,\bar{d})$ where $\bar{d}$ is induced by the quotient
    \begin{equation*}
        \wedge Z \aisom \frac{\mc{A} \tensor \wedge Z}{\aplus{\mc{A}} \tensor \wedge Z}.
    \end{equation*}

    We will drop the adjective ``Sullivan'' from here on out, except if the extension in question is not minimal.
    Abusing terminology, we may call $\mc{A} \tensor \wedge Z$ itself a minimal extension of $\mc{A}$.
    We may write $\sulext{\mc{A}}{Z}$ for $\mc{A} \tensor \wedge Z$ in the event that $\im{d\res{Z}} \sse \mc{A}$ (i.e., under the nilpotence condition $Z = Z(1)$), and we call such a minimal extension \term{elementary}.

    We are primarily interested in studying connected DGAs of the form $\mc{A}= \wedge V$ for some graded vector space $V = \set{V_n}_{n \in \NN}$ that can be built by minimal extensions of the form $\wedge V_{\le n} \tensor \wedge V_{n+1}$.
    Here, $\vsdeg{V}{\le n}$ denotes $\set{\vsdeg{V}{k}}_{k \le n}$.
    
    Suppose for instance that the extension $\wedge V_{\le n} \tensor \wedge \vsdeg{V}{n+1}$ is minimal.
    If $\vsdeg{V}{n+1} = \vsstep{\vsdeg{V}{n+1}}{c}$, then this extension decomposes into a sequence of $c$ elementary extensions.
    If this is true for every $n$, we say that $\wedge V$ is \term{$c$-step nilpotent}.
    If $c = 1$ then every extension is elementary; hence for every $n$ the differential satisfies the more strict requirement
    \begin{equation}
        \im{d\res{\vsdeg{V}{n}}} \sse \aprod{\ge 2}\vsdeg{V}{< n}
    \end{equation}
    and we call $\wedge V$ \term{simple} or \term{abelian}.
    In the event $\vsdeg{V}{1} = 0$ the above expression is automatically satisfied, and we call $\wedge V$ \term{simply connected}.  

    Recall that a map of DGAs is a \term{quasi-isomorphism} if it induces an isomorphism on cohomogology.
    A \term{minimal (Sullivan) model} of a DGA morphism $\map{\phi}{(\mc{A},d)}{(\mc{B},d)}$ is a quasi-isomorphism $\map{m_{\phi}}{(\mc{A} \tensor \wedge Z,d)}{(\mc{B},d)}$ where $\mc{A} \tensor \wedge Z$ is a minimal extension of $\mc{A}$, and where $m_\phi(a \tensor 1) = \phi(a)$ for all $a \in \mc{A}$.
    Again, we will drop the adjective ``Sullivan'' and simply call such models \term{minimal}.

    Note that our terminology disagrees slightly with that in \cite{felix2015rational}. We require that $Z = \set{\vsdeg{Z}{n}}_{n \ge 1}$, but they relax this condition to allow for nontrivial subspaces of degree $n = 0$ in the fiber.
    In this way, their definition of a minimal Sullivan model is slightly more broad; however, we will soon restrict our attention to a context where these definitions agree. 
    (See \cite[\textsection 3.7]{felix2015rational}.)

    The existence of minimal models of fibrations is a well-established fact:

    \begin{theorem}[{\cite[Thm. 14.12]{felix2012rational}}]
        Let $\map{\phi}{(\mc{A},d)}{(\mc{B},d)}$ be a morphism of connected DGAs for which $H^1(\phi)$ is an injection. 
        Then $\phi$ has a (unique) minimal model. 
    \end{theorem}

    In particular, the minimal model of a connected DGA always exists:

    \begin{proposition}[{\cite[Prop. 1.8]{felix2015rational}}]
        If $(\mc{B},d)$ is a connected DGA, there is a minimal Sullivan algebra $\wedge V$ and a quasi-isomorphism $\phi: (\wedge V,d) \to (\mc{B},d)$.
    \end{proposition}
    Abusing terminology, we call $\wedge V$ itself the \term{minimal Sullivan model} of $\mc{B}$ in this situation.
    In particular, we can now define the minimal model of a topological space.
    As before, we will frequently drop the adjective ``Sullivan.''

    \begin{definition}
        If $Y$ is a Riemannian manifold or simplicial complex, the \term{minimal model of $Y$} is $\minimal{Y}$, the minimal model of $\Forms{Y}$.
        If $p : Y \to B$ is a principal $K(G,n)$ fibration for $G$ abelian, the \term{minimal model of $p$} is the minimal model of $\map{p^*}{\Forms{B}}{\Forms{Y}}$.
    \end{definition}
    
    It is worth noting that the differential $d$ of a DGA decomposes into $d = d_1 + d_2 + \cdots$ where $\im{d_k\res{V}} \subseteq \wedge^{k+1}V$ so that $d_k$ maps an indecomposable to a linear combination of products of $k+1$ indecomposables.
    Of particular note is the \textbf{quadratic} component $d_1$. 
    It is straightforward to see that if $(\wedge V,d)$ is a DGA, restricting to the quadratic component $(\wedge V, d_1)$ yields another DGA, which we call the \textbf{associated quadratic} DGA. 
    The remaining components $d - d_1$ of the differential are a \textbf{deformation}, which are studied for instance in \cite{schlessinger2012deformation}.
    We will eventually be interested in studying the setting when the differential itself $d = d_1$ is itself quadratic. 
    When $(\minimal{Y},d) = (\minimal{Y},d_1)$ we will call $Y$ \textbf{coformal} (see \cite{neisendorfer1978formal}).

\subsection{Realization and Rationalization}
    In \cite{hilton2016localization} it is shown that if $G$ is finitely generated nilpotent, there is a \term{rationalization} (equivalently, $\set{0}$-localization) of $G$, by which we mean that there is a group denoted $G \tensor \QQ$ in which $g \mapsto g^p$ is a bijection for all $p \neq 0$.
    Moreover, the rationalization functor $G \mapsto G \tensor \QQ$, thought of as the Malcev completion of $G$, is exact in the category of ($c$-step) nilpotent groups.
    When $M$ is abelian, this rationalization $M \tensor \QQ$ agrees with the standard tensor product with $\QQ$.

    \begin{definition}
        Let $X$ and $Y$ be nilpotent spaces.
        Suppose $f : X \to Y$ is such that $f_k : \pi_k(X) \tensor \QQ \to \pi_k(Y) \tensor \QQ$ is an isomorphism for all $k$. 
        In this case we call $f$ a \term{rational (homotopy) equivalence}.
        Two spaces $X$ and $X'$ are \term{rationally equivalent} if there is a chain of rational equivalences connecting them:
        \[
            X \rto{} Y_0 \lto{} Y_1 \rto{} \cdots \lto{} Y_N \rto{} X'
        \]
    \end{definition}

    Minimal models provide a functor from topological spaces to differential graded algebras $\minimal{\square} : \mathbf{Top} \to \mathbf{DGA}$. 
    Dually there is a ``realization'' functor $\realize{\;\cdot\;} : \mathbf{DGA} \to \mathbf{Top}$ for producing a topological space out of the given algebraic data.
    
    \begin{definition}
        Given a space $X$ the \term{rationalization of $X$} is $X_{\QQ} = \realize{\minimal{X}}$.
    \end{definition}
    A central power of rational homotopy theory is that for nilpotent spaces, these functors are inverses up to rational equivalence.
    That is
    \begin{theorem}[{\cite[Proposition 8.2]{felix2015rational}}]
        Given a nilpotent space $X$, the geometric realization of the minimal model $\map{m_X}{\minimal{X}}{\Forms{X}}$ induces a rational equivalence $$X \to \realize{\Forms{X}} \rto{\realize{m_X}} \realize{\minimal{X}} = X_{\QQ}.$$
    \end{theorem}
    
    Since we are interested in nilpotent spaces up to rationalization, we will always assume that we are working with any rationally equivalent space when possible.
            Since duality between principal $K(G,n)$ fibrations and elementary extensions allows us to make natural identifications between certain rational invariants, we are able to consolidate some of our current notation in the nilpotent setting.

    When $M = \pi_n(Y)$ for some $n \ge 2$ we define
    \begin{equation*}
        \pindual{n} = \pindualverbose{n}{Y} \defeq 
            \Hom(M, \QQ) = \Hom(\pi_n(Y), \QQ).
    \end{equation*}
    In the situation $G = \pi_1(Y)$, it is fitting to be slightly more delicate. 
    If we assume that $G$ is finitely generated, there is a Lie group $\hat{G}$ called the \term{Malcev completion} of $G$ for which $G \subseteq \hat{G}$ is a cocompact lattice.
    In fact, $\hat{G}$ has rational structure constants, meaning that there is a basis $\set{X_1,\dots,X_m}$ of the Lie algebra $\mf{g} = \mf{lie}(\hat{G})$ for which all bracket relations have rational coefficients:
    \begin{equation*}
        [X_{i},X_{j}] \in \mathrm{span}_{\QQ}\set{X_1,\dots,X_m},
    \end{equation*}
    and one views $\hat{G}$ as a rational Lie group.
    Following \cite[\textsection 3.2.1]{felix2008algebraic} one may define $\pindual{1} = \mf{g}^{*} = \Hom(\mf{g},\QQ)$.
    In \cite{griffiths2013rational} they denote $G \tensor \QQ = \mf{g}$, and we will often adopt this notation as well.
    We say that $Y$ is of \term{finite type} if $\pindual{n}$ is a finite dimensional $\QQ$ vector space for every $n$.
    Unless otherwise stated, we assume that all nilpotent spaces are of finite type.
    Moreover define
    \begin{equation*}
        \pindualstep{n}{j} = \pindualstepverbose{n}{j}{Y} \defeq 
            \Hom(\pinstep{n}{j},\QQ).
    \end{equation*}
    (Recall the definition of $\pinstep{n}{j}$ from the end of \Cref{ssec:nilpotent spaces}.)
    Furthermore define
    \begin{equation*}
        \pindualupto{n}{J} \defeq \Dirsum_{j = 1}^{J} \pindualstep{n}{j}.
    \end{equation*}
    In particular, observe that $\pindualupto{n}{J} \aisom \pi_n(\postspacestep{Y}{n}{J}) \tensor \QQ$, and that $\pindualupto{n}{c} \aisom V_n \aisom \pi_n(\postspace{Y}{n}) \tensor \QQ$.
    A key consequence in this context is that the minimal model of $Y$ can be described by $\minimal{Y} = \wedge V$ where $V = \set{\vsdeg{V}{n}}_{n \ge 1}$ with $\pindual{n} \aisom \Hom(\pi_n(Y), \QQ) \aisom \pi_n(Y) \tensor \QQ$ for all $n$.
    We make this identification without confusion.
    In fact, the principal refinement of the Postnikov tower of $Y$ induces a sequence of minimal extensions whose limit is $\minimal{Y}$.
    That is, the principal $K(\pinstep{n}{j},n)$ fibration 
    \begin{equation*}
        \map{\postmapstep{n}{j}}{\postspacestep{Y}{n}{j}}{\postspacestep{Y}{n}{j-1}}
    \end{equation*}
    is precisely dual to an elementary extension 
    \begin{equation*}
        \sulext{\vsstep{V_n}{j-1}}{\pindualstep{n}{j}} = \vsstep{V_n}{j}.
    \end{equation*}

    The following powerful theorem connects some of the dualities we have been discussing in a bit more generality.

    \begin{theorem}[{\cite[eq. 3.8, and Thm. 5.1]{felix2015rational}}]\label{thm:model of fiber}
        Let $F \to Y \to B$ be a fibration of path connected Riemannian manifolds of finite type.
        Then there is a minimal extension $\minimal{B} \to \minimal{B} \tensor \wedge V$ such that the following diagram of DGAs commutes
        \begin{diagram}
            (\minimal{B},d) \rar \arrow[d,"m_B","\quasiisom"']
                \& (\minimal{B} \tensor \wedge V, d) \rar \arrow[d,"m_Y","\quasiisom"']
                \& (\wedge V,\bar{d}) \arrow[d,"m_F"] \\
            \Forms{B} \rar 
                \& \Forms{Y} \rar 
                \& \Forms{F}
        \end{diagram}
        Moreover, if $\pi_1(B)$ acts locally nilpotently on $H^*(F;\QQ)$ then $m_F : \wedge V \rto{\quasiisom} \Forms{F}$ is a quasi-isomorphism, and hence a minimal model for $F$.
    \end{theorem}

    If the fibration $F \to Y \to B$ is a principal $K(G,n)$ fibration, it is shown in \cite[\textsection 12.2]{griffiths2013rational} that the $k$-invariant $[k] \in H^n(B;G)$ is dual to the differential $d$ of $\minimal{B} \tensor \wedge V$. 
    That is, $[d] \in H^n(\minimal{B};V)$ corresponds to $[k] \tensor 1 \in H^n(B;G) \tensor \QQ$ under the following, which follows by applications of the Universal Coefficient Theorem:
    \begin{lemma}
        $H^n(\minimal{B};V) \aisom H^n(B;G) \tensor \QQ$.
    \end{lemma} 

    In particular, then, as in the proof of \cite[Proposition 15.3]{griffiths2013rational}, we have that the following diagram commutes for any (finite) complex $X$:
    \begin{equation}\label{eq:obstructions agree}
        \begin{tikzcd}[ampersand replacement = \&]
            {[X,B]} \arrow[rr,"\text{obstruction}"]\arrow[rr,swap,"\text{to lifting}"] \arrow[dd]
                \& {} \& {H^n(X;G)} \dar{\tensor \QQ} \\
            {}
                \& {} \& {H^n(X;V)} \dar{\aisom} \\
            {[\minimal{B},\minimal{X}]} \arrow[rr,"\text{obstruction}"]\arrow[rr,swap,"\text{to extending}"]
                \& {} \& {H^n(\minimal{X};V)}
        \end{tikzcd}
    \end{equation}

\subsection{Algebraic Obstruction Theory} 
Following now \cite[\textsection 11]{griffiths2013rational}, we present some homotopy theory for DGAs, with a quantitative focus as in \cite[\textsection 3.2]{Manin_2019}. 
    Let $\mc{A}$ and $\mc{B}$ be DGAs with morphisms $\phi,\psi : \mc{A} \to \mc{B}$.
    An \term{(algebraic) homotopy} from $\phi$ to $\psi$ is then a morphism
    \begin{equation*}
        \Phi : \mc{A} \to \mc{B} \trtdt
    \end{equation*}
    for which $\Phi\res{\substack{t=0 \\ dt=0}} = \phi$ and $\Phi\res{\substack{t=1 \\ dt=0}} = \psi$.
    For such an algebraic homotopy $\Phi$, we will write $\Phi^i_j$ for the $t^i (dt)^j$ term of $\Phi$.
    Moreover, we think of $\rtdt \defeq \RR\gen{t^{(0)},dt^{(1)}}$ as an algebraic model for the unit interval $I \defeq [0,1]$.
    In the event that $\mc{B} = \Forms{X}$ for a sufficiently nice space $X$, we may realize this algebraic model as a genuine interval by a morphism 
    \begin{equation*}
        \rho : \Forms{X} \trtdt \to \Forms{X \times I}.
    \end{equation*}
    We may further utilize this algebraic model of the unit interval as a tool for integration of homotopies.
    That is, we define a linear operator $\aint_0^1 : \mc{B} \trtdt \to \mc{B}$ by
    \begin{equation*}
        \aint_0^1 b \tensor t^i = 0 \qquad \text{and} \qquad  \aint_0^1 b \tensor t^i dt = (-1)^{\deg{b}}\frac{1}{i+1} b,
    \end{equation*}
    as well as an operator $\aint_0^t : \mc{B} \trtdt \to \mc{B} \trtdt$ by
    \begin{equation*}
        \aint_0^t b \tensor t^i = 0 \qquad \text{and} \qquad  \aint_0^t b \tensor t^i dt = (-1)^{\deg{b}}\frac{1}{i+1}b \tensor t^{i+1}.
    \end{equation*}

    These operators will allow us to discuss obstructions to extending homotopies, and serve as formal analogues of fiberwise integration (see \cite{bott2013differential}).
    Before this discussion, however, we state the following facts which are easily proved from the given definitions.

    \begin{proposition}
        These operators satisfy the conditions 
        \begin{align}
            d\l(\aint_0^1 u\r) + \aint_0^1 du  & = u\res{\substack{t = 1 \\ dt = 0}} - u\res{\substack{t = 0 \\ dt = 0}}  \label{eq:fundamental theorem 1} \\ 
            d\l(\aint_0^t u\r) + \aint_0^t du  & = u - u\res{\substack{t = 0 \\ dt = 0}} \tensor 1. \label{eq:fundamental theorem 2}
        \end{align}
    \end{proposition}

    \begin{corollary}
        If $\Phi : \mc{A} \to \mc{B} \trtdt$ is a homotopy from $\phi$ to $\psi$ then 
        \begin{equation*}
            d\l(\aint_0^1 \Phi(a)\r) + \aint_0^1 d\Phi(a) = \psi(a) - \phi(a).
        \end{equation*}
    \end{corollary}

    These tools allow us to develop algebraic obstruction theory, which provides algorithmic processes for computing algebraic obstructions. These obstructions are used to determine whether a partial algebraic homotopy (a map which is a homotopy when restricted to a sub-DGA) may be extended to the entire DGA.

        Given a morphism $\phi : \mc{A} \to \mc{B}$ we may define the \term{relative cohomology} of $\phi$ as the cohomology of the complex $(C^n(\phi),d) = (\mc{A}^{n} \dirsum \mc{B}^{n-1},d)$ where $d(a,b) = (da, \phi(a) - db)$. 
    We denote this cohomology ring by any of the following means, sometimes with vector valued coefficients 
    \begin{equation*}
        H^*(\phi : \mc{A} \to \mc{B}) = H^*(\phi) = H^*(\mc{A} \to \mc{B}).
    \end{equation*}
    Moreover (see \cite[\textsection 10.3]{griffiths2013rational}) there is an exact sequence 
    \begin{equation*}
        \cdots \to H^{n-1}(\mc{A}) \to H^{n-1}(\mc{B}) \to H^{n}(\phi) \to H^{n}(\mc{A}) \to \cdots
    \end{equation*}

    \begin{proposition}[{\cite[Proposition 11.1]{griffiths2013rational}}]\label{}
        Let $Z = \vsdeg{Z}{n}$ and let $\iota : \mc{A} \to \mc{A}\gen{Z}$ be an elementary Sullivan extension of a DGA $\mc{A}$.
        Suppose we are given the following (non-commutative) DGA diagram
        \begin{diagram}
            \mc{A} \rar{f} \arrow[d,hook,"\iota"]
                \& \mc{B} \dar{\eta} \\
            \sulext{\mc{A}}{Z} \rar{g} 
                \& \mc{C}
        \end{diagram}
        as well as a homotopy $\Phi$ from $\eta \circ f$ to $g\res{\mc{A}} \defeq g \circ \iota$. 

        Then the obstruction class $[\mc{O}] \in H^{n+1}(\eta; Z)$ to producing a map $\tilde{f} : \mc{A}\gen{Z} \to \mc{B}$ such that $\tilde{f} \circ \iota = f$, as well as a homotopy $\tilde{\Phi}$ from $\eta \circ \tilde{f}$ to $g$ for which $\tilde{\Phi} \circ \iota = \Phi$, is represented by 
        \begin{align*}
            \mc{O} : Z & \to \mc{B}^{n+1} \dirsum \mc{C}^{n} \\
            z & \mapsto \l(f(dz),\, g(z) + \aint_0^1 \Phi(dz)\r).
        \end{align*}
        
        Moreover, in the event that the obstruction vanishes, there are maps $(b,c) : Z \to \mc{B}^n \dirsum \mc{C}^{n-1}$ such that $d(b,c) = \mc{O}$, which means
        \begin{align*}
            db(z) & = f(dz), \\
            dc(z) & = \eta \circ b(z) - g(z) - \aint_0^1H(dz).
        \end{align*}
        The map $\tilde{f}$ can then be defined by $\tilde{f}(z) = b(z)$ and the homotopy extension $\tilde{\Phi}$ can be defined by 
        \[
            \tilde{\Phi}(z) = g(z) + d(c(z) \tensor t) + \aint_0^t \Phi(dz).
        \]
    \end{proposition}

    As a relative version of this, we have 
    \begin{proposition}[{\cite[Proposition 11.2]{griffiths2013rational}}]\label{prop:relative obstruction}
        Let $Z = \vsdeg{Z}{n}$ and let $\iota : \mc{A} \to \mc{A}\gen{Z}$ be an elementary Sullivan extension of a DGA $\mc{A}$.
        Suppose we are given the following (non-commutative) DGA diagram
        \begin{diagram}
            \mc{A} \rar{f} \arrow[d,hook,"\iota"]
                \& \mc{B} \dar{\eta} \rar{\mu}
                \& \mc{D} \\
            \sulext{\mc{A}}{Z} \rar{g} 
                \& \mc{C} \urar{\nu}
                \& 
        \end{diagram}
        such that 
        \begin{enumerate}
            \item $g\res{A} \homotopic \eta f$ by a homotopy $\Phi : \mc{A} \to \mc{C} \tensor \RR\gen{t,dt}$ with $\nu \circ \Phi$ constant
            \item $\mu$ is surjective
            \item the righthand triangle commutes exactly (i.e., $\nu \circ \eta = \mu$ exactly)
            \item the outer portion of the diagram commutes exactly (i.e., $\mu \circ f = \nu \circ g\res{\mc{A}}$)
        \end{enumerate}
        Then the obstruction class $[\mc{O}] \in H^{n+1}(\eta; Z)$ to producing a map $\tilde{f} : \mc{A}\gen{Z} \to \mc{B}$ such that $\tilde{f} \circ \iota = f$, as well as a homotopy $\tilde{\Phi}$ from $\eta \circ \tilde{f}$ to $g$ for which $\tilde{\Phi} \circ \iota = \Phi$, and for which $(\nu \circ \tilde{\Phi} = (\mu \circ \tilde{f}) \tensor 1)$ is constant, is represented by a map 
        \begin{align*}
            \mc{O} : Z & \to \mc{B}^{n+1} \dirsum \mc{C}^{n} \\
            z & \mapsto \l(f(dz),\, g(z) + \aint_0^1 \Phi(dz)\r).
        \end{align*}
    \end{proposition}

    We will often make use of the following specific relative version of this fact:

    \begin{proposition}\label{prop:relative qualitative obstruction}
        Let $\sulext{\mc{A}}{Z}$ be an elementary extension of a DGA $\mc{A}$. 
        Suppose there are maps
        \begin{diagram}
            \sulext{\mc{A}}{Z} \rar{\phi,\psi}
                \& \mc{B} \arrow[r, two heads, "\mu"]
                \& \mc{C}
        \end{diagram}
        with $\mu$ surjective, for which $\phi\res{\mc{A}} \homotopic \psi\res{\mc{A}}$ by a homotopy $\Phi$, and such that $\mu \circ \phi \homotopic \mu \circ \psi$ by a homotopy $\chi : \sulext{\mc{A}}{Z} \to \mc{C} \trtdt$ extending $\mu \circ \Phi$. 
        The obstruction $[\mc{O}] \in H^n(\mu;Z)$ to extending $\Phi$ to a homotopy $\tilde{\Phi} : \mc{A} \to \mc{C} \trtdt$ which simultaneously lifts $\chi$ and extends $\Phi$ is given by 
        \[ 
            \mc{O}(z) = \l(\psi(z) - \phi(z) - \aint_0^1\Phi(dz); \chi(z)\r) .
        \]
    \end{proposition}

    The proofs of these facts can be found in \cite[\textsection 3.1]{Manin_2019}.
        Now, fix $X$ a finite piecewise Riemannian simplicial complex, and a principal $K(G,n)$ fibration $Y \to B$ in which $G$ is abelian and $V = \vsdeg{V}{n} = \Hom(G, \QQ)$ is a finite dimensional vector space. 
    As in \Cref{thm:model of fiber}, there is a minimal model $\minimal{Y} \defeq \minimal{B} \tensor \wedge V \rto{m_Y} \Forms{Y}$ relative to the fibration.
        Furthermore, equip $\Forms{X}$ with the $L^{\infty}$ norm.
    Fix some norms on $\minimal{B}$ and on $V = \vsdeg{V}{n}$.
    Assume moreover that $n \le \dim{X}$. 
    Since our spaces are of finite type, any dependence on this choice of norm for $V$ will be scaled only up to a constant. 
    Let $\phi : \minimal{Y} \to \Forms{X}$.
    Assume that there is some fixed definition of $\Dil{\phi\res{\minimal{B}}} \in \RR^{\ge 0}$.
    For instance, if $\minimal{B} = \wedge W$ is a minimal model (say for a nilpotent space) we set 
    \begin{equation*}
        \Dil{\phi\res{\minimal{B}}} \defeq 
        \max_{1 \le k \le \dim X}
        \set{\|\phi\res{W_{k}}\|^{1/k}_{op}}.
    \end{equation*}
    Having fixed a definition for $\Dil{\phi\res{\minimal{B}}}$ define the \term{dilatation} of $\phi$ to be
    \begin{equation*}
        \Dil{\phi} \defeq 
        \max\set{\Dil{\phi\res{\minimal{B}}},\|\phi\res{V}\|_{op}^{1/n}}.
    \end{equation*}
    For example, suppose that $\bar{f} : X \to B$ is $L$-Lipschitz.
    Then $\bar{f}^* : \Forms[k]{B} \to \Forms[k]{X}$ satisfies 
    \begin{equation*}
        \|\bar{f}^*\omega\|_{\infty} \le L^k\|\omega\|_{\infty}.
    \end{equation*}
    Then under an appropriate definition we have $\Dil{\bar{f}^*m_B} \le CL$.
    Suppose further that $f : X \to Y$ is a lift of $\bar{f}$ for which $\|f^*m_Y\res{V}\|^{1/n}_{op} \le L$. 
    Since $(f^*m_Y)\res{\minimal{B}} = \bar{f}^*m_B$, we have $\Dil{f^* m_Y} \le C'L$ as well, where $C'$ depends only on $m_B$ and the norms of $\Forms{B}$ and $V$.
    More specifically if $f : X \to Y$ is itself $L$-Lipschitz, then $\Dil{f^*m_Y} \le C'L$.

    Given a DGA homotopy $\Phi : \minimal{Y} \to \Forms{X} \trtdt$, we can define the dilatation once we fix a particular realization of the algebraic interval into its geometric counterpart. 
    For each $T \in \RR^{> 0}$ we define a realization $\rho_{T} : \Forms{X} \trtdt \to \Forms{X \times I}$ by $t \mapsto t/T$, using $T$ as a heuristic for the ``formal length'' of the formal interval.
    With this, we may define a family of dilatations for $\Phi$ by 
    \begin{equation*}
        \Dil[T]{\Phi} \defeq \Dil{\rho_T \Phi}.
    \end{equation*}

    The following proposition provides quantitative bounds on extending homotopies beteween maps, quantizing the algebraic information from \Cref{prop:relative qualitative obstruction}.

    \begin{proposition}[{\cite[Proposition 3.9]{Manin_2019}}]\label{prop:quantitative relative obstruction}
        Let $\phi, \psi : \minimal{Y} \to \Forms{X}$, and suppose there is a homotopy $\Phi : \minimal{B} \to \Forms{X} \trtdt$ from $\phi\res{\minimal{B}}$ to $\psi\res{\minimal{B}}$.
        \begin{enumerate}
            \item The obstruction to extending $\Phi$ to a homotopy $\tilde{\Phi} : \minimal{Y} \to \Forms{X} \trtdt$ from $\phi$ to $\psi$ is a class in $H^{n}(X;V)$ represented by a cochain $\sigma$ for which (for any $\tau > 0$) we have
            \begin{equation*}
                \|\sigma\|_{op} \le \tau C(n,d\res{V})\Dil[\tau]{\Phi}^{n+1} + \Dil{\phi}^{n} + \Dil{\psi}^n.
            \end{equation*}
            \item In the event the obstruction vanishes, we can choose $\tilde{\Phi}$ so that 
            \begin{equation*}
                \|\tilde{\Phi}_i^j\res{V}\|_{op} \le (C_{IP} + 2)\l(\tau C(n,d\res{V})\Dil[\tau]{\Phi}^{n+1} + \Dil{\phi}^{n} + \Dil{\psi}^{n}\r)
            \end{equation*}
            where $C_{IP}$ is the isoperimetric constant for $(n+1)$-forms in $X$, and $\tau > 0$ is arbitrary.
        \end{enumerate}
        Furthermore, consider the relative situation (i.e., there is a subcomplex  $A \sse X$ for which we have an existing homotopy $\chi : \minimal{Y} \to \Forms{A} \trtdt$ between $\phi\res{A}$ and $\psi\res{A}$).  If the obstruciton from \Cref{prop:relative qualitative obstruction} vanishes, we can get similar bounds using a relative isoperimetric constant together with an additional $O(\tau \Dil[\tau]{\chi})$ term.
    \end{proposition}

    In \cite[\textsection 3.2]{Manin_2019}, it is shown that if $\Dil{\phi},\Dil{\psi} \le L$, and a homotopy $\Phi$ between them can be built formally up through degree $n$ without encountering any nonzero obstructions, then \textit{as long as $\minimal{Y}$ is simply connected} with minimal model $\minimal{Y} = \wedge V$ for $V = \set{\vsdeg{V}{n}}_{n \ge 2}$, we can create more specific bounds on the dilatation of $\Phi$.
    Namely, it is shown that for $j \in \set{0,1}$ we have
    \begin{equation*}
        \|\Phi^j_i\res{\vsdeg{V}{k}}\|_{op} \le CL^{2k - 2}
    \end{equation*}
    for some constant $C = C(k,X,Y)$ depending additionally on the choice of norms for each $\vsdeg{V}{n}$.
    We revisit this in \Cref{sec:upper bounds} to produce upper bounds on dilatation when $Y$ is a nilpotent space.

\section{The Shadowing Principles}\label{sec:shadows}
In this, the heart of the paper, we generalize the Shadowing Principle of \cite{Manin_2019} to targets which appear as principal $K(G,n)$ fibrations, and hence to nilpotent spaces.
We would ultimately like to get an upper bound on the volume of nullhomotopies in nilpotent spaces, and the original principle provided a useful tool for such analysis in the simply connected setting (see \cite[Theorem 5.5]{Manin_2019}).
By the help of \Cref{thm:lifting shadows}, we attain our goal in \Cref{sec:upper bounds}.

Before turning to our versions of the Shadowing Principle, we provide some intuition for interpreting the principle outside of the fibered setting.
Given a target $Y$ and an $n$-dimensional space $X$, we assume that we have an algebraic map $\phi : \minimal{Y} \to \Forms{X}$ which has dilatation bounded by a quantity $L$.
We might wonder how far away we must look to find a genuine map which has similar bounds.
That is, we put some reasonable metric on $\Hom(\minimal{Y},\Forms{X})$, such as 
\begin{equation*}
    d(\phi,\psi) = \inf\set{\operatorname{length}(\Phi) \sothat \phi \stackrel{\Phi}{\homotopic} \psi}
\end{equation*}
where the \term{formal length} of such a homotopy $\Phi$ is
\begin{equation*}
    \operatorname{length}(\Phi) \defeq \Dil{\aint_0^1 \Phi}.
\end{equation*}
This gives constant homotopies trivial length, is nondegenerate according to applications of \eqref{eq:fundamental theorem 2}, and satisfies the triangle inequality by \cite[Prop. 3.11]{Manin_2019}.
With this metric in mind, the Shadowing Principle states that if $\phi$ is itself in the homotopy class of some genuine map $f : X \to Y$ and $\Dil{\phi} \le L$, then there is $O(L)$-Lipschitz map $g \homotopic f$ whose pullback $g^*m_Y$ is at most $O(L)$ far from $\phi$.
If it is too hard to construct a bounded geometric map with certain properties, it may be easier to construct one whose pullback is bounded.
The principle then yields a genuine map homotopic to the original map, but with bouded Lipschitz constant.

In \Cref{ssec:Simplicial Geometry} we introduce establish necessary quantitative results in simplicial geometry, which are crucial to the proof of the fibered and nilpotent versions of the principle. 
In \Cref{ssec:Shadowing Principle Proofs} we prove the relevant principles.

\subsection{Simplicial Geometry}\label{ssec:Simplicial Geometry}
We work primarily in the simplicial category as our model for homotopy theory. 
See \cite{goerss2009simplicial} for details on this setting, and its equivalence to the model category of CW complexes.
In particular, fibrations will all be (geometric realizations of) fibrations of finite simplicial sets.
A key ingredient in this work (as in \cite{Chambers_DMW_2018_nullcobordism, Chambers_MW_2018_quantitative_nullhomotopy, Manin_2019}) is that of simplicial approximation. 
Given a simplicial complex $X$ with standard metric, one of the most common subdivision schemes of $X$ is barycentric subdivision. 
This scheme has certain geometric disadvantages though; simplices may become arbitrarily skinny, causing massive distortion in comparison to the standard simplex.

We say that a simplicial subdivision of $X$ is \term{$L$-regular} if there is a constant $r > 0$ (depending on the dimension of $X$) for which simplices are $r$-bilipschitz to a standard simplex with edge length $\frac{1}{L}$.
Although barycentric subdivision is not regular, there are a handful of schemes (see \cite{ferry2013quantitative,edelsbrunner1999edgewise} for examples) which are.
We assume from here on out that all subdivisions will be regular.

Lipschitz maps between simplicial complexes can be approximated by simplicial maps without increasing the Lipschitz constant too much.
First, recall the following definitions:
\begin{definition}
    Let $H : X \times I \to Y$ be a Lipschitz homotopy. 
    After embedding $Y \into \RR^{N}$ in a (piecewise) smooth manner, Rademacher's Theorem implies that $H$ is differentiable almost everywhere (since $H$ is Lipschitz), and hence one may define the \textbf{volume} of $H$ by
    \[
        \vol{H} \defeq \int_{X \times I}|\operatorname{Jac}(H)|\dvol{}.
    \]
    The \term{length} of $H$ is 
    \begin{equation*}
        \operatorname{length}(H) \defeq \max \set{\Lip(H\res{\set{x} \times I}) \sothat x \in X}.
    \end{equation*}
    The \term{thickness} of $H$ is 
    \begin{equation*}
        \operatorname{thickness}(H) \defeq \max \set{\Lip(H\res{X \times \set{t}}) \sothat t \in I}.
    \end{equation*}
    Observe that for $X$ an $n$-dimensional space, $\vol(H) = O(\operatorname{length}(H)) \cdot O(\operatorname{thickness}(H)^n)$.
\end{definition}

\begin{proposition}[Quantitative Simplicial Approximation Theorem {\cite[Prop. 2.1]{Chambers_DMW_2018_nullcobordism}}]\label{prop:simplicial approximation}
    Let $X$ be a finite-dimensional simplicial complex with standard simplices, and let $Y$ be a finite simplicial complex with piecewise linear metric. 
    There is a constant $C = C(X,Y)$ such that any $L$-Lipschitz map $f : X \to Y$ has a simplicial approximation $f_{\Delta}$ with $\Lip{f_{\Delta}} \le C(L+1)$.
    Moreover, $f$ is homotopic to $f_{\Delta}$ by a $C$-length homotopy $H$ of thickness $C(L+1)$.
    In particular, $\Lip{H} \le C'(L+1)$ for some $C' = C'(X,Y)$.
    \begin{proof}[Sketch of proof]
        Paying quantitative attention to the proof of the (non-quantiative) simplicial approximation theorem from \cite{hatcher2002algebraic}, one takes a regular subdivision of $X$ and uses the Lebesgue covering lemma to carefully produce nearly-linear homotopies between the uncontrolled map and its simplicial approximation. 
        By induction up the skeleta of $X$, this process is made to be globally well-controlled.
    \end{proof}
\end{proposition}

Recalling notation from \Cref{sec:background}, we are primarily interested in diagrams $\diag{f}$ of the form 
\begin{diagram}
    A \rar{\alpha}
    \arrow[d,"\iota", hook]
        \& Y\dar[two heads]{p} \\
    X \arrow[ur,"f"]
        \arrow[r,"\beta"]
        \& B
\end{diagram}
where $A \into X$ is a simplicial pair and $Y \onto B$ is a principal $K(G,n)$ fibration.
We will often write $f\res{A} = \alpha$ and $\bar{f} = \beta$.

\begin{definition}
    Let $\map{f}{X}{Y}$.
    If each of the maps $\bar{f}$ and $f\res{A}$ in $\diag{f}$ are $L$-Lipschitz, we say that $\diag{f}$ (or $f$) is \term{externally $L$-Lipschitz}. 
    If $f$ is externally equivalent to some $L$-Lipschitz map, we say that $\diag{f}$ (or $f$) is \term{$L$-Lipschitz at heart}. 
    If $f$ is internally equivalent to some $L$-Lipschitz map, we say that $f$ is \term{internally $L$-Lipschitz}. 
\end{definition}

If we have a map $\map{f}{X}{Y}$ which is externally $L$-Lipschitz, we would like to bound the best-behaved fillings of the corresponding diagram $\diag{f}$.
In the setting of maps from the standard simplex, we can be a bit more precise.
Let $\map{f}{\simplex}{Y}$.
We would like to come up with a map $f_0$ such that, for each $\xi \in G$  there is a map $f_{\xi}$ representing $\xi$ in the $G$-torsor relative to $f_0$.
Moreover, we would like to bound $\Lip(f_{\xi})$ for each $\xi$.

One useful tool for quantifying the behavior of the maps we build will be the language of shard complexes and mosaic maps.
\begin{definition}
    For each $k$, let $\mc{F}_{k}$ be a finite set of maps $\simplex[k] \to Y$.
    The collection $\mc{F} = \bigcup_{k} \mc{F}_{k}$ is called a \term{shard complex} on $Y$, and is naturally a finite subcomplex of the singular simplicial complex on $Y$. 
    In particular, $\mc{F}$ is a semi-simplicial set whose $k$-simplices are $\mc{F}_{k}$.
    We call a map $f : X \to Y$ from a simplicial complex $X$ an \term{$\mc{F}$-mosaic} map if $f\res{\simplex[k]} \in \mc{F}_{k}$ for all $k$-simplices $\simplex[k] \sse X$, for all $k$.
\end{definition}

Note that being mosaic is preserved under postcomposition. 
For example, if $\phi : Y \to Y'$ is a homotopy equivalence, then $\phi \circ f$ will be mosaic with respect to a shard complex $\mc{F}'$ on $Y'$ if $f$ is mosaic with respect to some $\mc{F}$ on $Y$.

For the sake of brevity and clarity, it will be convenient to prove results for $1$-Lipschitz maps.
However, such results will often times generalize to $L$-Lipschitz maps.
In particular, if $X \to Y$ is $L$-Lipschitz, one subdivides $X$ at scale $\frac{1}{L}$ and rescales so that simplices are unit size. 
The same map then becomes $1$-Lipschitz, and when convenient we will assume without loss of generality that $L = 1$.

Fix a principal $K(G,n)$-fibration of finite simplicial sets $\map{p}{Y}{B}$.
A (minimal) fibration of simplicial sets has a geometric realization which is locally trivial over simplices in the base. 
In other words, if the fiber of $p : Y \to B$ is $F$, then for any simplex $\sigma \subseteq B$ we may identify $p^{-1}(\sigma)$ with $F \times \sigma$. 
(See, for example, the proof of \cite[Theorem I.10.9]{goerss2009simplicial}.) 
We make this identification freely in what transpires and after sufficient subdivision we may assume that for $b \in \skel{B}{0}$ we have $p^{-1}(\openstar{b}) \isom \openstar{b} \times F$, where $\openstar{b}$ denotes the star of $b$.
By rescaling $Y$ and $B$ as necessary, we assume that if $s : \simplex \to B$ is $1$-Lipschitz, then $\im{s} \sse \openstar{b}$ for some $b \in B^{(0)}$.

\begin{proposition}\label{prop:elect rep}
    Let $L > 0$ be given.
    There is an $L_0 = L_0(L,p)$ such that the following holds:
    Suppose the diagram 
    \begin{diagram}
        A
            \arrow[r,"\alpha"]
            \dar{}
        \& Y 
            \dar{}
    \\
        X
            \arrow[r,"\beta"']
            \arrow[ur,"f"]
        \& B
    \end{diagram}
    is externally $L$-Lipschitz.
    Assume moreover that $\alpha$ is mosaic with respect to some fixed shard complex $\mc{Z}$ on $Y$ for which $\mc{Z}^{(0)}$ consists of basepoints over fibers of vertices in $B$.
    Then $f$ is $L_0$-Lipschitz at heart through some $f_0$.
    Moreover, for each $\xi \in H^n(X,A;G)$, there is an $L_{\xi}$ such that $f$ is $L_{\xi}$-Lipschitz at heart through some $f_{\xi}$ for which the obstruction to performing an internal homotopy from $f_{\xi}$ to $f_{0}$ is $\xi$.
    \begin{proof}
        Let $\mf{D}$ be the collection of diagrams 
        \begin{diagram}
            A
                \arrow[r,"\alpha"]
                \dar{}
            \& Y 
                \dar{}
        \\
            X
                \arrow[r,"\beta"']
                \arrow[ur,dashed]
            \& B
        \end{diagram}
        for which $\beta$ is simplicial, for which $\alpha$ is mosaic with respect to $\mc{Z}$, and whose obstruction to finding a lift in the diagonal vanishes.
        Given such a diagram $\mc{D} \in \mf{D}$, let $f_{\mc{D},0}$ be a filling of minimal Lipschitz constant.
        Because $\alpha$ is assumed to be mosaic and $\beta$ is assumed to be simplicial, the set $\mf{D}$ is finite, and hence the quantity $C_{\simplex,0} \defeq \max_{\mc{D} \in \mf{D}}\set{f_{\mc{D},0}}$ is finite.

        Furthermore, given $\xi \in H^{n}(X,A;G)$, let $f_{\mc{D},\xi}$ be of minimal Lipschitz constant such that $[f_{\mc{D},\xi}] = \xi$ in the torsor relative to $f_{\mc{D},0}$. (That is, the obstruction to performing an internal homotopy from $f_{\mc{D},\xi}$ to $f_{\mc{D},0}$ is $\xi$.) 
        Again, the quantity $C_{\mc{D},\xi} = \max_{\mc{D} \in \mf{D}}\set{\Lip(f_{\mc{D},\xi})}$ is finite.
        Therefore, if we had first assumed that $\beta$ were simplicial, we would be done.

        It suffices then to show that, given an arbitrary $L$-Lipschitz diagram $\mc{D}$, one can find a ``nearby'' simplicial diagram $\mc{D}_{\simplex}$ of the above form whose fillings $f_{\mc{D}_{\simplex},0}$ and $f_{\mc{D}_{\simplex},\xi}$ induce $L_0$-Lipschitz fillings of $\mc{D}$.
        Indeed, this is easily attained by similar methods to \Cref{prop:simplicial approximation}. 
        At the risk of increasing the Lipschitz constant slightly we attain a simplicial approximation $\beta'$ of $\beta$, as well as a linear homotopy $H$ from $\beta$ to $\beta'$.
        By the homotopy lifting property, we get the following commutative diagram
        \begin{diagram}
                A \times I
                    \dar{}
                    \rar{H_{\alpha}}
                \& Y \dar{p}
            \\
                X \times I 
                    \rar[swap]{H_{\beta}}
                \& B 
            \end{diagram}
            where at $t = 1$ the diagram is filled as 
            \begin{diagram}
                A
                    \dar{}
                    \rar{\alpha'}
                \& Y \dar{p}
            \\
                X 
                    \rar[swap]{\beta'}
                \& B 
            \end{diagram}
            This is our diagram $\mc{D}_{\simplex}$.
            If $\mc{Z}$ were the shard complex with respect to which $\alpha$ was mosaic, this induces a new shard complex $\mc{Z}'$ for which $\alpha'$ is $O(L+1)$-Lipschitz and mosaic.
            By our previous work, since $\beta'$ is simplicial, there is a map $f_{\mc{D}_{\simplex},\xi}$ fitting into the diagonal of $\mc{D}_{\simplex}$ with the desired properties.
            By the homotopy lifting property, we attain a homotopy from $f_{\mc{D}_{\simplex},\xi}$ to some map $f_{\xi}$ at $t = 0$ for which the following diagram commutes:
            \begin{diagram}
                A
                    \dar{}
                    \rar{\alpha}
                \& Y \dar{p}
            \\
                X 
                    \urar{f_{\xi}}
                    \rar[swap]{\beta}
                \& B 
            \end{diagram}
            In fact, by the identification $p^{-1}(\openstar{b}) \isom \openstar{b} \times F$ for each $b \in \skel{B}{0}$ and all $t \in I$ we may define the homotopy $H_{\alpha}(t)$ to be piecewise linear. 
            Therefore $f_{\xi} \defeq H_{\alpha}(0)$ has Lipschitz constant $L_{\xi}$ depending only on the geometry of the local trivializations and on $L$, as well as on the initial shard complex $\mc{Z}$.
        \end{proof}
\end{proposition}

\begin{corollary}\label{cor:no lower obstructions}
    Let $\mc{Z}_{0}$ be a fixed shard complex on $Y$ for which $\mc{Z}_{0}^{(0)}$ consists of basepoints over fibers of vertices in $B$.
    There is a $C = C(\mc{Z}_{0},n-1,Y)$ such that, if $f : X \to Y$ is externally $1$-Lipschitz, and $f\res{A}$ is mosaic, then $f\res{X^{(n-1)}}$ is internally equivalent to a $C$-Lipschitz map.
    In particular there is a shard complex $\mc{Z}_{n-1} = \mc{Z}_{n-1}(\mc{Z}_{0},n-1,Y)$ with respect to which $f\res{\skel{X}{n-1}}$ is mosaic. 
    \begin{proof}
        We proceed by induction.
        Because $Y$ is finite, there is a $C_{0} = \diam{Y}$ for which $f\res{\skel{X}{0}}$ is $C_{0}$-Lipschitz.
        Inductively, for $0 < k < n$ assume that $f$ is internally homotopic to a map $\tilde{f}_{k-1}$ for which $\tilde{f}_{k-1}\res{\skel{X}{k-1}}$ is $C_{k-1}$-Lipschitz.
        Apply \Cref{prop:elect rep} to the diagram
        \begin{diagram}
            A \cup \skel{X}{k-1}
                \dar{}
                \arrow[rrr,"\tilde{f}_{k-1}\res{A \cup \skel{X}{k-1}}"]
            \&\&\& Y \dar{p}
        \\
            A \cup \skel{X}{k}
                \arrow[urrr,"\tilde{f}_{k-1}\res{\skel{X}{k}}",pos=0.35]
                \arrow[rrr,swap,"p \circ f\res{A \cup \skel{X}{k}}"]
            \&\&\& B 
        \end{diagram}
        to attain an (\textit{a priori}) externally equivalent $C_{k}$-Lipschitz map $\tilde{f}_{k}\res{\skel{X}{k}} : \skel{X}{k} \to Y$ for which $C_{k}$ depends only on $C_{k-1}$ and on $p : Y \to B$.
        The obstruction to performing an internal homotopy between $\tilde{f}_{k}\res{\skel{X}{k}}$ and $\tilde{f}_{k-1}\res{\skel{X}{k}}$ is an element of $H^{n}(A \cup \skel{X}{k},A\cup \skel{X}{k-1};G)$. 
        This obstruction is trivial, however, since $k \neq n$.
        In particular, $\tilde{f}_{k}\res{\skel{X}{k}}$ must in fact be \textit{internally} equivalent to $f_{k-1}\res{\skel{X}{k}}$ over this diagram, and we extend $\tilde{f}_{k}$ over higher dimensional cells to a map $\tilde{f}_{k}$ on all of $X$ which is internally equivalent (relative to $A$) to $f$ itself.
        Take $\mc{Z}_{k} \defeq \mc{Z}_{k-1} \cup \set{\tilde{f}_{k}\res{\simplex} \sothat \simplex \sse X \text{ is a $k$-simplex}}$. 
        Induction ceases at $k = n-1$.
    \end{proof}
\end{corollary}

Of primary importance for our purposes is the situation where the domain is the standard $n$-simplex, $(\simplex[n],\boundary \simplex[n])$.

\begin{corollary}\label{cor:elect rep}
    Let $p : Y \to B$ be a principal $K(G,n)$ fibration of simplicial sets.
    Let $\xi \in G \aisom H^{n}(\simplex, \boundary \simplex; G)$.
    There is a constant $C_{\xi} = C_{\xi}(p)$ such that if $f : \simplex \to Y$ is externally $1$-Lipschitz, then there is a $C_{\xi} $-Lipschitz map $f_{\xi}$ externally equivalent to $f$.
    Moreover, the obstruction to performing an internal homotopy from $f_0$ to $f_{\xi}$ is $\xi$ (where $f_0$ is the map applied to $\xi = 0$).
\end{corollary}

The following lemma from \cite[\textsection 2.2]{Manin_2019} will be necessary for producing quantitative bounds on integration over simplices. 
As in that paper, we adapt the notation $\Forms{X}$ to denote the algebra of simplex-wise smooth forms on a simplicial complex $X$. 
Given a simplicial pair $(X,A)$, we use $\Forms{X,A}$ to denote the forms whose restriction to $A$ is trivial.
We refer the reader to \cite{Manin_2019} for the proof.

\begin{lemma}[Second Quantitative Poincar\'{e} Lemma]\label{lem:second quantitative poincare}
    For every $0 < k < n$ there is a constant $\tilde{C}_{n,k}$ such that, whenever $\omega \in \Forms[k]{\Delta^n}$ is a closed $k$-form, and $\alpha_{\boundary} \in \Forms[k-1]{\boundary\Delta^{n}}$ is a $(k-1)$-form with $d\alpha_{\boundary} = \omega\res{\boundary\Delta^n}$, then there is a $(k-1)$-form $\alpha \in \Forms[k-1]{\Delta^{n}}$ extending $\alpha_{\boundary}$, such that $d\alpha = \omega$, and for which $\|\alpha\|_{\infty} \le \tilde{C}_{n,k}\l(\|\omega\|_{\infty} + \|\alpha_{\boundary}\|_{\infty}\r)$.
    In the event $k = n$, the same holds if the pair additionally $(\omega, \alpha_{\boundary})$ satisfies Stokes' Theorem, so that $\int_{\Delta^{n}} \omega = \int_{\boundary\Delta^{n}} \alpha_{\boundary}$.
\end{lemma}

\subsection{The Shadowing Principles}\label{ssec:Shadowing Principle Proofs}
For the rest of this section we fix $G$ a finitely generated abelian group, we fix a principal $K(G,n)$ fibration $Y \rto{p} B$ of simplicial sets $Y$ and $B$.
Moreover, fix norms on $\minimal{B}$ and $V = G \tensor \QQ$. 
Define $F_{b}$ to be the fiber over a given point $b \in B$, and let $F$ denote a generic fiber.
Assume moreover that the skeleta of $Y$ and $B$ are finite in every dimension.
For example this is the case at each stage of the Postnikov tower of a finite type nilpotent space.

Recall that $\ms{D} = \ms{D}(X,A,Y,B)$ is the set of all diagrams of the form 
\begin{diagram}
    A \rar{} \dar[hook]{\iota}
        \& Y \dar[two heads]{p} \\
    X \rar{} 
        \& B
\end{diagram}
whose obstruction class vanishes, for some simplicial pair $(X,A)$.
We also set $\ms{D}_{\simplex} \defeq \ms{D}(\simplex,\boundary \simplex,Y,B)$.
Recall furthermore that $\map{f}{X}{Y}$ is \term{internally $L$-Lipschitz} if $f$ is internally equivalent to an $L$-Lipschitz map, relative to $\diag{f} \in \ms{D}(X,A,Y,B)$. 
That is, there is an $L$-Lipschitz map $\map{g}{X}{Y}$ and an internal homotopy from $f$ to $g$ relative to $\diag{f}$.

Our goal, in some sense, is to produce a constant $C$ depending only on $\map{p}{Y}{B}$ and $N = \dim{X}$ such that every map $f$ which is externally $1$-Lipschitz will be internally $C$-Lipschitz.
However, this will only be possible if the corresponding algebraic data is bounded in the right ways.
Recall that the following diagram commutes up to homotopy, where the vertical maps are quasi-isomorphisms:
\begin{diagram}
    (\minimal{B}, d) \rar{\Id \tensor 1} \arrow[d,"m_B"]
        \& (\minimal{B} \tensor \wedge V,d)\rar  \arrow[d,"m_Y"]
        \& (\wedge V,\bar{d}) \arrow[d,"m_F"]
        \\
    \Forms{B} \rar{p^*} 
        \& \Forms{Y} \rar
        \& \Forms{F}
\end{diagram}
The top row is a minimal Sullivan model for the morphism $p^* : \Forms{B} \to \Forms{Y}$.
Fix now a map $\map{f}{X}{Y}$, and define $\bar{f} \defeq p \circ f$.

    A \term{quadruple} for $(\bar{f},f)$ as above consists of DGA morphisms $(\bar{\phi},\phi,\bar{\Phi},\Phi)$ where 
    \begin{enumerate}[label=(Q\arabic*)]
        \item $\bar{\phi} : \minimal{B} \to \Forms{X}$ 
        \item $\phi : \minimal{B} \tensor \wedge V \to \Forms{X}$ 
        \item $\bar{\Phi} : \minimal{B} \to \Forms{X} \trtdt$ 
        \item $\Phi : \minimal{B} \tensor \wedge V \to \Forms{X} \trtdt$
    \end{enumerate}
    such that 
    \begin{enumerate}[resume,label=(Q\arabic*)]
        \item $\bar{\phi} \homotopic \bar{f}^*m_B$ via $\bar{\Phi}$ 
        \item $\phi \homotopic f^*m_Y$ via $\Phi$ 
        \item $\Phi\res{A}$ is constant (in particular, $\phi\res{A} = f^*m_Y\res{A}$).
    \end{enumerate}
    We say that $(\phi,\Phi)$ \term{extends} $(\bar{\phi},\bar{\Phi})$ if in addition 
    \begin{enumerate}[label=(E\arabic*)]
        \item $\phi\res{\minimal{B}} \defeq \phi \circ (\Id \tensor 1) = \bar{\phi}$ 
        \item $\Phi\res{\minimal{B}} \defeq \Phi \circ (\Id \tensor 1) = \bar{\Phi}$.
    \end{enumerate}

\begin{figure}[H]
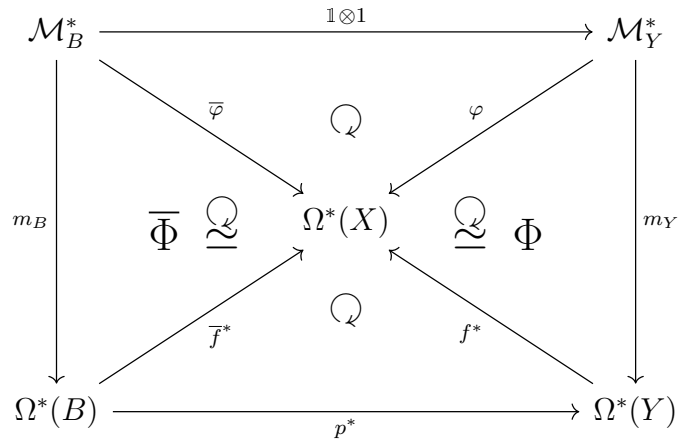

        \begin{diagram}
            \minimal{B} \arrow[rrrr,"\Id \tensor 1"] \arrow[dddd,swap,"m_B"] \arrow[ddrr, "\bar{\phi}"]
                \& 
                \& 
                \& 
                \& \minimal{Y} \arrow[dddd,"m_Y"] \arrow[ddll, swap, "\phi"]
                \\
    
                \& 
                \& \commutes
                \& 
                \& 
                \\
    
                \& {\arrow[mysymbol]{}[description]{{\text{\Large$\bar{\Phi}$}}\;\;\;\;\homotopicCommuteSymbol\;\;}}
                \& \Forms{X}
                \& {\arrow[mysymbol]{}[description]{\;\; \homotopicCommuteSymbol \;\;\;\; {\text{\Large$\Phi$}}}}
                \& 
                \\
    
                \& 
                \& \commutes
                \& 
                \& 
                \\
            \Forms{B} \arrow[rrrr,swap,"p^*"] \arrow[uurr, swap, "\bar{f}^*"]
                \& 
                \& 
                \& 
                \& \Forms{Y} \arrow[uull, "f^*"]
        \end{diagram}
        \caption{Interaction of Quadruples}
        \label{fig:quadruples}
    \end{figure}

    Unless otherwise stated, all quadruples will be assumed to be two pairs that extend each other.
    Ignoring the complex $A$, quadruples fit into the diagram of \Cref{fig:quadruples}, where the top and bottom triangles commute exactly, while the left and right triangles commute up to homotopy through $\bar{\Phi}$ and $\Phi$ respectively.

    We say that this quadruple $(\bar{\phi},\phi,\bar{\Phi},\Phi)$ is \term{bounded by $L$ downstairs} if 
    \begin{enumerate}[label=(B\arabic*)]
        \item $\Dil{\bar{\phi}} \le L$ 
        \item $\Dil{\phi} \le L$
        \item $\Dil[1/L]{\bar{\Phi}} \le L$.
    \end{enumerate}
    The quadruple is \term{totally bounded by $L$} (or a \term{totally $L$-bounded quadruple}) if in addition to being bounded by $L$ downstairs we have 
    \begin{enumerate}[resume,label=(B\arabic*)]
        \item $\Dil[1/L]{\Phi} \le L$.
    \end{enumerate}

For a subcomplex $X' \sse X$ we say that $(\bar{\phi},\phi,\bar{\Phi},\Phi)\res{X'}$ is a totally $L$-bounded quadruple for $(\bar{f},f)\res{X'}$ if $(\bar{\phi},\phi,\bar{\Phi},\Phi)$ is a quadruple for $(\bar{f},f)$ whose appropriate restrictions to $X'$ are all $L$-bounded.
In this situation we may say more concisely that \term{the restriction of the quadruple to $X'$} is bounded.

The remainder of this section has strong roots in the original paper \cite{Manin_2019}, but the current approach (induction up the Postnikov tower) allows the subsequent ideas to be presented somewhat atomically.
Before moving into the relevant lemmata, we set up some final notation which will help streamline some of the proofs.

\begin{notation}
    Given maps $\map{f,g}{\simplex}{Y}$ for which $f\res{\boundary \simplex} = g\res{\boundary \simplex}$, define $\stackmap{g}{f} : \SS^{n} \to Y$ to be the map given by $\stackmap{g}{f}\res{+} = g$ on the northern hemisphere and $\stackmap{g}{f}\res{-} = f$ on the southern hemisphere, glued along their common boundary.
\end{notation}

Recall from \Cref{cor:homotopy torsor} that since $\map{p}{Y}{B}$ is principal, if $f$ and $g$ are externally equivalent, the obstruction to performing an internal homotopy from $f$ to $g$ defines a class 
\begin{equation*}
    \l[\stackmap{g}{f}\r] \in H^{k}(\simplex,\boundary \simplex;\pi_n(F)) = H^{k}(\simplex,\boundary \simplex;G) \aisom \begin{cases}
        G & k = n, \\
        0 & k \neq n.
    \end{cases}
\end{equation*}
We may write simply $[f]$ for $\l[\stackmap{g}{f}\r]$ if $g$ is fixed ahead of time.

\begin{definition}
    Let $\diag{} \in \ms{D}_{\simplex}$ be an externally $1$-Lipschitz diagram.
    The map $f_{\xi}$ as given by \Cref{cor:elect rep} is called the \term{elect representative of $\xi \in G$ in $\diag{}$}. 
    This representative is defined relative to the diagram, as long as the diagram allows some diagonal filling.
    If $\map{f}{\simplex}{Y}$ is a specific map for which $\diag{f} = \diag{}$, we may also say that $f_{\xi}$ is the \term{elect representative of $\xi$ for $f$}, though in reality $f_{\xi}$ is completely agnostic about the map $f$.
\end{definition}

In what follows we let $\square \tensor 1: G \to G \tensor \QQ \aisom V^*$ be the obvious map, noting that in case $n = 1$ this is a slight abuse of notation.
We may similarly denote $(\square \tensor 1)^* : G \to V$.

\begin{lemma}\label{lem:bounded obstruction}
    There is a $\tilde{C} > 0$ so that the following holds.
    Suppose that $\map{f}{\simplex}{Y}$ is externally $1$-Lipschitz, and let $\bar{f} \defeq f \circ p$.
    Let $(\bar{\phi},\phi,\bar{\Phi},\Phi)$ be a quadruple for $(\bar{f},f)$ which is $1$-bounded downstairs, and whose restriction to $\boundary \simplex$ is totally $1$-bounded.
    
    Let $f_0$ be the elect representative of $0 \in G$ relative to $f$, and let $\mc{R} \defeq \stackmap{f_0}{f}$.
    Then $[\mc{R}] \tensor 1$ is contained in a $\tilde{C}$-ball around 
    \begin{equation*}
        \int_{\simplex} \l(-\phi\res{V} + f^*m_Y\res{V} + \aint_0^1\Phi d \res{V}\r)\in V^*.
    \end{equation*}
    \begin{proof}
        Let $C_0$ be  given by \Cref{cor:elect rep} applied to $\xi = 0$, let $C(n-1,d)$ be the constant from \Cref{prop:quantitative relative obstruction}, and define $\tilde{C} = C(n-1,d) + (C_0)^{n} + 1$.        
        This constant depends only on $n$ and the norms of $V$.

        The class $[\mc{R}] \tensor 1$ is the obstruction in $\minimal{B} \tensor \wedge V$ to homotoping $\mc{R}^*m_Y$ to $0$, or equivalently, to another nullhomotopic map, such as $\stackmap{\phi}{\phi}$, the map which is $\phi$ on each hemisphere.
        First, define $\Psi\res{\minimal{B}} \defeq \Phi\res{\minimal{B}}$ on each hemisphere.
        We now evaluate the obstruction to extending $\Psi$ to $V$.
        According to \Cref{prop:relative qualitative obstruction} this is given by 
        \begin{equation*}
            \l[-\stackmap{\phi}{\phi} + \mc{R}^*m_Y + \aint_0^1 \Psi d\res{V}\r] \in V^*.
        \end{equation*}
        Now because $\Dil[1]{\Phi\res{\minimal{B}}} = \Dil[1]{\bar{\Phi}} \le 1$, by \Cref{prop:quantitative relative obstruction} the obstruction on the northern hemisphere is represented by a class $\mc{O}$ of norm
        \begin{equation*} 
            \|\mc{O}\|_{op}
                \le C(n-1,d)\Dil[1]{\Psi\res{\minimal{B}}}^{n+1} + \Dil{f_0^*m_Y}^{n} + \Dil{\phi}^{n} = C(n-1,d) + (C_0)^n + 1  = \tilde{C},
        \end{equation*}
        while on the southern hemisphere this class sends $v \in V$ to 
        \begin{equation*}
            \int_{\simplex} (-\phi(v) + f^*m_Y(v) + \aint_0^1 \Phi(dv)).
        \end{equation*}
        Therefore, $[\mc{R}] \tensor 1$ is within distance $\tilde{C}$ of the term as claimed.
    \end{proof}
\end{lemma}

\begin{lemma}\label{lem:local fibered shadows}
    There is a constant $C$ such that the following holds.
    Suppose the following diagram commutes:
    \begin{diagram}
        A \rar{f\res{A}}\dar[hook]{\iota}
            \& Y \dar[two heads]{p}\\
        X \rar{\bar{f}} \urar{f}
            \& B
    \end{diagram}
    Suppose that $f\res{A}$ and $\bar{f}$ are $1$-Lipschitz.
    Let $(\bar{\phi},\phi,\bar{\Phi},\Phi)$ be a quadruple for $(\bar{f},f)$ which is $1$-bounded downstairs, and whose restriction to $A$ is totally $1$-bounded.
    Suppose moreover that $\Lip{f\res{\skel{X}{n-1}}} \le 1$ and $\Dil[1]{\Phi\res{\skel{X}{n-1}}} \le 1$.
    
    Then there is a homotopy $H$ from $f$ to a map $\map{\tilde{f}}{X}{Y}$ such that 
    \begin{enumerate}
        \item $H\res{(\skel{X}{n-2} \cup A) \times \set{t}} = f\res{\skel{X}{n-2} \cup A}$ for all $t \in I$ (i.e., $H$ is a homotopy $\rel{\skel{X}{n-2} \cup A}$)
        \item $p \circ H_t = \bar{f}$ for all $t \in I$ (i.e., $H$ is a fiberwise homotopy)
        \item $\tilde{f}\res{\skel{X}{n-1}} = f\res{\skel{X}{n-1}}$ 
        \item $\tilde{f}$ is $C$-Lipschitz.
    \end{enumerate}
    \begin{proof}
        We assume without loss of generality that the star of $A$ retracts to $A$ after further subdivision if necessary.
        Let $\simplex \sse X$ be an $n$-simplex.
        For each $\zeta \in G$, let $C_{\zeta}$ be as in the consequence of \Cref{cor:elect rep}.
        Let $\tilde{C}$ be as in the result of \Cref{lem:bounded obstruction}. 
        Define 
        \begin{equation*}
            G_{\tilde{C}} \defeq \set{\zeta \in G \sothat \|\zeta \tensor 1\|_{V^*} \le \tilde{C} + 1}.
        \end{equation*}
        Since $\ker{\square \tensor 1}$ is torsion, $\#G_{\tilde{C}} < \infty$ is a constant depending only on the same data as $\tilde{C}$, and so we may define a constant $C_{n}$ depending again only on the same data as $\tilde{C}$ by 
        \begin{equation*}
            C_{n} \defeq \max \set{C_{\zeta} \sothat \zeta \in G_{\tilde{C}}}.
        \end{equation*}

        Now, let the map $f$ be given.
        For each $\zeta \in G_{\tilde{C}}$, let $\map{f_{\zeta}}{\simplex}{Y}$ be the elect representative of $\zeta$ relative to $\diag{f}$, as given by \Cref{cor:elect rep}.
        We may define a shard complex on $Y$ by
        \[
            \mc{Z} = \set{f_{\zeta}\res{\simplex} \;\bigg\vert\; f_{\zeta}\res{\boundary \simplex} = f\res{\boundary \simplex},\, \zeta \in G_{\tilde{C}}} \cup \mc{A}
        \]
        where $\mc{A}$ is a given shard complex for which $f \res{A \cup \skel{X}{n-1}}$ is mosaic.
        
        We show that there is an element $\xi \in G_{\tilde{C}}$ such that (the restriction to $\simplex$ of) the elect representative $\tilde{f} \defeq f_{\xi}$ of $\xi$ relative to $f$ is mosaic with respect to this complex, and that this shard complex has bounded Lipschitz constant.
        The process goes by building a sufficiently nice homotopy on the boundary of each $\simplex[n] \sse \skel{X}{n}$, so that the obstruction to extending to $\simplex$ is as desired.
        
        Let $\beta = \aint_0^1\Phi$, and observe that for $v \in V$ we have $d\beta(v) = \phi(v) - f^*m_Y(v) - \aint_0^1\Phi(dv)$.
        Let $b \in \cochains{X,A;G}{n-1}$ be such that the image of $\inner{b,\sigma} \in G$ inside $V^*$ is as close in norm as possible (say, within distance $1$) to $\l\|\aint_{\sigma} \beta\res{V}\r\|_{V^*}$ for any $(n-1)$-simplex $\sigma \sse \skel{X}{n-1}$. 
        Let $\Sigma \sse \skel{X}{n}$ be a simplex.
        Set $\xi \defeq [f] + \inner{\coboundary b,\Sigma}$. 
        Let $\tilde{f}\res{\Sigma} \defeq (f\res{\Sigma})_{\xi}$ be the elect representative of $\xi$ relative to $f\res{\Sigma}$.
        In particular, $f_{\xi}$ is externally equivalent to $f\res{\Sigma}$.
        
        \begin{claim}
            The map $\tilde{f}\res{\Sigma}$ has $\Lip{\tilde{f}\res{\Sigma}} \le C_{n}$.
            \begin{proofofclaim}
                Since $\tilde{C}$ was the result of \Cref{lem:bounded obstruction}, observe that $\torsor{f\res{\Sigma}} \tensor 1$ is contained in a $\tilde{C}$-ball around
                
                \begin{equation*}
                    \int_{\Sigma} (-\phi\res{V} + f^*m_Y\res{V} + \aint_0^1\Phi d \res{V}) = -\displaystyle\int_{\Sigma} d\beta(v).
                \end{equation*}
                
                As a result, $\torsor{f\res{\Sigma}} \tensor 1$ is contained within a $\tilde{C} + 1$ ball around $-\inner{\coboundary b,\Sigma} \tensor 1$ by Stokes' Theorem.
                In particular, 
                \[  
                    \|\xi \tensor 1\| = \l\|\Big([f\res{\Sigma}] + \inner{\coboundary b, \Sigma}\Big) \tensor 1\r\| \le \tilde{C} + 1.
                \]
                Therefore, $\xi \in G_{\tilde{C}}$.
                Since $\tilde{f}\res{\Sigma}$ was elect, we have $\Lip{\tilde{f}} \le C_{\xi} \le C_{n}$.        
            \end{proofofclaim}
        \end{claim} 
        
        It remains to construct the desired homotopy $H$ from $f$ to $\tilde{f}$.
        For all $(n-2)$-cells $\tau$ Define $H\res{\tau \times \set{t}} = f\res{\tau}$ for all $t \in I \defeq [0,1]$.
        Let $\sigma \sse \Sigma$ be an $(n-1)$-face. 
        Define $H\res{\sigma \times \set{1}} \defeq \ H\res{\sigma \times \set{0}} = f\res{\sigma}$.
        This provides a complete definition for $H\res{\boundary(\sigma \times I)}$, which makes the following diagram commute for each $(n-1)$-face $\sigma \sse \Sigma$:
        \begin{diagram}
            \boundary(\sigma \times I) \arrow[d,hook] \arrow[r,"H\res{\boundary(\sigma \times I)}"]
                \& Y \arrow[d]
                \\
            \sigma \times I \arrow[r,"(\bar{f} \circ \pi)\res{\sigma}",swap] \arrow[ur,dashed]
                \& B
        \end{diagram}
        By \Cref{cor:homotopy torsor} the obstruction to extending this to a map which fits in the dashed diagonal defines an element of $G$.
        There is well-defined extension $H\res{\sigma \times I}$ whose degree is $\inner{b,\sigma} \in G$.
        That is, precomposing $H$ with a linear retract of $\sigma$ induces a canonical nullhomotopy of this map, and one can find a filling $H\res{\sigma \times I}$ whose degree is as desired.
        
        This gives a complete definition of $H$ on $\boundary(\Sigma \times I) = (\Sigma \times \set{0}) \cup (\boundary \Sigma \times I) \cup (\Sigma \times \set{1})$.
        It remains to show that $H$ extends to all of $\Sigma \times I$, which is the result of the following claim:

        \begin{claim}
            There is an extension of $H\res{\Sigma}$ which defines a homotopy from $f\res{\Sigma}$ to $\tilde{f}\res{\Sigma}$.
            \begin{proofofclaim}
                Consider the following diagram:
                \begin{diagram}
                    \SS^{n} \isom (\Sigma \times \set{0}) \cup (\boundary \Sigma \times I) \cup (\Sigma \times \set{1}) \arrow[d,swap,hook] \arrow[rrr,"f \cup H\res{\boundary \Sigma \times I} \cup \tilde{f}"]
                    \&
                        \&
                        \& Y \dar[two heads]{p}
                        \\
                    \DD^{n+1} \isom \Sigma \times I \arrow[urrr,dashed] \arrow[rrr,swap,"\bar{f} \circ \pi"]
                        \&
                        \&
                        \& B
                \end{diagram}
                By \Cref{prop:obstruction to relative extension}, an extension of $H$ will fit in this diagram if and only if the obstruction class $\torsor{H} \in G$ vanishes.
                By Stokes' Theorem this means that such a homotopy exists if and only if 
                \begin{equation*}
                    \torsor{H\res{\Sigma \times \set{1}}} - \torsor{H\res{\Sigma \times \set{0}}} = \torsor{H\res{(\boundary\Sigma) \times I}}.
                \end{equation*}
                In particular, this happens if and only if
                \begin{equation*}
                    [\tilde{f}] - [f] = \inner{\coboundary b,\Sigma} 
                \end{equation*}
                which is true by definition of $[\tilde{f}] = \xi = [f] + \inner{\coboundary b,\Sigma}$.
                Therefore $H$ extends to an internal homotopy from $f$ to $\tilde{f}$.
            \end{proofofclaim}
        \end{claim}

        Therefore, $\tilde{f}\res{\skel{X}{n}}$ is $C_{n}$-Lipschitz. 
        In fact, $\tilde{f}\res{\skel{X}{n}}$ is mosaic with respect to $\mc{Z}_{n}$.
        We redefine $\tilde{f}_{n} = \tilde{f}$ on higher skeleta of $X$ using \Cref{prop:elect rep} as in the proof of \Cref{cor:no lower obstructions}.
        Namely, inductively for $n < k \le N = \dim{X}$ the diagram
        \begin{diagram}
            A \cup \skel{X}{k-1}
                \dar{}
                \arrow[rrr,"\tilde{f}_{k-1}\res{A \cup \skel{X}{k-1}}"]
            \&\&\& Y \dar{p}
        \\
            A \cup \skel{X}{k}
                \arrow[urrr,"\tilde{f}_{k-1}\res{\skel{X}{k}}"']
                \arrow[rrr,swap,"p \circ f\res{A \cup \skel{X}{k}}"]
            \&\&\& B 
        \end{diagram}
        is externally $C_{k-1}$-Lipschitz.
        We attain an (\textit{a priori}) externally equivalent $C_{k}$-Lipschitz map $\tilde{f}_{k}\res{\skel{X}{k}} : \skel{X}{k} \to Y$ for which $C_{k}$ depends only on $C_{k-1}$ and on $p : Y \to B$.
        As before, the obstruction to performing an internal homotopy between these maps is an element of $H^{n}(A \cup \skel{X}{k},A\cup \skel{X}{k-1};G)$, which is trivial since $k \neq n$.
        In particular, $\tilde{f}_{k}\res{\skel{X}{k}}$ must in fact be \textit{internally} equivalent to $f_{k-1}\res{\skel{X}{k}}$ over this diagram, and we extend $\tilde{f}_{k}$ over higher dimensional cells to a map $\tilde{f}_{k}$ on all of $X$ which is internally equivalent to $f$ itself, and whose restriction to lower skeleta is as desired.
        This extends the shard complex as well.
        When $k = N$ we attain (noting slight equivocation on the symbol $\tilde{f}$) the desired map $\tilde{f} \defeq \tilde{f}_{N}$. 
    \end{proof}
\end{lemma}

 \begin{lemma}\label{lem:local algebraic shadows}
    Let $(\bar{\phi},\phi,\bar{\Phi},\Phi)$ and $(\bar{f},f)$ be as in the setup to \Cref{lem:local fibered shadows} and let $\tilde{f}$ and $H$ be as in the conclusion.
    In particular, assume that $\Dil[1]{\Phi\res{\skel{X}{n-1}}} \le 1$.

    Let $\pi : X \times I \to X$ be the obvious projection. 
    Then there is a homotopy 
    \begin{equation*}
        \Psi : \minimal{B} \tensor \wedge V \to \Forms{X \times I} \tensor \RR\gen{s,ds}
    \end{equation*} 
    from $\pi^*\phi$ to $H^*m_Y$ for which $\Psi\res{t = 0} = \Phi$, and such that $\tilde{\Phi} \defeq \Psi\res{t = 1}$ is a homotopy from $\phi$ to $\tilde{f}^*m_Y$ $\rel{A}$ for which $\Dil[1]{\tilde{\Phi}} \lesssim 1$.
    In particular, $(\bar{\phi},\phi,\bar{\Phi},\tilde{\Phi})$ is a quadruple for $(\bar{f},\tilde{f})$ which  is totally $(\lesssim 1)$-bounded.
    \begin{proof}
        First, define $\Psi\res{\minimal{B}} \defeq \pi^*\Phi\res{\minimal{B}} = \pi^*(\bar{\Phi} \tensor 1)$.
        Observe that $\Psi\res{\minimal{B}}$ is a homotopy from $\pi^{*}\phi\res{\minimal{B}}$ to $H^*m_Y\res{\minimal{B}}$, since 
        \begin{equation*}
            (\Psi\res{\minimal{B}})\res{X\times \set{0}}
                = \pi^*(\Phi\res{t=0})\res{\minimal{B}} 
                = \pi^*\phi\res{\minimal{B}}
        \end{equation*}
        and (because $p \circ H = p \circ f \circ \pi$) we furthermore have
        \begin{equation*}
            (\Psi\res{\minimal{B}})\res{X\times \set{1}} 
                = \pi^*(\bar{\Phi}\res{t=1} \tensor 1)\res{\minimal{B}} 
                = \pi^*\bar{f}^*m_B 
                = (p \circ f \circ \pi)^*m_B 
                = H^*p^*m_B 
                = H^*m_Y\res{\minimal{B}}.
        \end{equation*}
        This defines $\Psi$ (hence, $\tilde{\Phi}$) on $\minimal{B}$, and in particular on $\im{d\res{V}}$.
        We seek to extend $\Psi$ to $V$.

        Recall that $V = \Hom(G,\QQ)$ is of degree $n$, and that $\beta = \aint_0^1 \Phi$.
        The primary goal is to define a bounded $\tilde{\beta}\res{\skel{X}{n}}$ on $V$ for which $\tilde{\beta}(v) = \phi(v) - \tilde{f}^*m_Y(v) - \aint_0^1 \Phi(dv)$,
        and hence to define
        \[
            \tilde{\Phi}(v) \defeq \tilde{f}^*m_Y(v) + d(\tilde{\beta}(v) \tensor s) + \aint_0^{s} \Phi(dv).
        \]
        On $\skel{X}{n-2}$ we may define $\Psi(v)$ constantly equal to $\Phi(v)$, but on $\skel{X}{n-1}$ we will want to be somewhat more careful.
        The obstruction to defining such a homotopy $\Psi$ from $\Phi$ to such a map $\tilde{\Phi}$ is (by \Cref{prop:relative obstruction}) a class $\mc{O}(v) \in H^{n}(X \times I, X \times \set{0,1}; V)$, defined by 
        \begin{equation*}
            \mc{O}(v) = (\pi^*\phi(v) - H^*m_Y(v) - \aint_0^1\Psi(dv) ; \beta(v) \dirsum \tilde{\beta}(v)).
        \end{equation*}
        We analyze this obstruction class in order to determine a proper definition for $\tilde{\beta}$.
        To aid in this, we denote $\mf{O}(v) \defeq \pi^*\phi(v) - H^*m_Y(v) - \aint_0^1\Psi(dv)$.
        First, note that $\mf{O}(v)\res{\skel{X}{n-2} \times I} = 0$, and therefore for any $(n-1)$-simplex $\sigma \sse \skel{X}{n-1}$ we have $\mf{O}(v)\res{(\boundary \sigma) \times I} = 0$ as well.
        Hence, by the Poincar\'{e} Lemma, there is a form $\mc{B}(v)$ for which $d\mc{B}(v) = \mc{O}(v)$, and whose restriction to $\skel{X}{n-2} \times I$ vanishes.
        
        Therefore, it is enough to analyze $\mf{B}(v) \res{X \times \set{0,1}}$ and thus to determine the class of $\tilde{\beta}(v)$.
        Let $b \in C^{n-1}(X,A;G)$ be as in the previous lemma. 
        In particular, recall that $H$ was defined so that $[H\res{\sigma \times I}] = \inner{b,\sigma}$.
        
        The following claim will allow us to apply Stokes' Theorem in order to evaluate $\mf{B}(v)$ in the desired manner:

        \begin{claim}\label{claim:integral}
            For any $(n-1)$-simplex $\sigma$ we have \[ \int_{\sigma \times I} d\mf{B}(v) = -\inner{b,\sigma}^*(v).\]
            \begin{proofofclaim}
                Identifying $\sigma \times \set{0} \sim \sigma \times \set{1}$ and flattening $\boundary \sigma \times I$ to a single copy $\boundary \sigma$, we form a quotient $Q : \sigma \times I \to \SS^{n}$.
                Observe that $\mf{O}(v)\res{\sigma \times \set{0}} = \mf{O}(v)\res{\sigma \times \set{1}}$, and that $\mf{O}(v)\res{\boundary \sigma \times I}$ is constant. 
                Similarly, by \Cref{lem:local fibered shadows} we have $H\res{\sigma \times \set{0}} = H\res{\sigma \times \set{1}}$ and constant $H\res{\boundary\sigma \times I}$.
                Therefore $\mf{O}(v)$ and $H$ each factor through $Q$. 
                In particular, we view $\int_{\sigma \times I}\mf{O}(v)$ as a map $V \to \Forms[n]{\SS^{n}} \aisom \RR$ (and therefore as an element of $V^{*} \aisom V$ itself) which is the obstruction to performing a homotopy from $H^*m_Y$ to $\pi^*(\phi\res{\sigma})$ within $\SS^{n}$.
                However, $\pi^*(\phi\res{\sigma})$ factors through $\Forms[*]{\sigma} \aisom \Forms[*]{\DD^{n-1}}$ and is therefore algebraically nullhomotopic.
                
                Hence, $\int_{\sigma \times I}d\mf{B}$ is equivalently seen as the obstruction to nullhomoping $H^*m_Y$ in the quotient.
                By \eqref{eq:obstructions agree}, the obstruction to an algebraic nullhomotopy of $H^*m_Y$ is the same as the rationalization of the obstruction to a geometric nullhomotopy of $H$ in this quotient. 
                But this obstruction is exactly measured by the homotopy class of the torsor element $-[H\res{\sigma \times I}] \in G$, and so as required we get 
                \begin{equation*}
                    \int_{\sigma \times I} d\mf{B}(v) = -[H\res{\sigma \times I}]^*(v) = -\inner{b,\sigma}^*(v).
                    \qedhere
                \end{equation*}
            \end{proofofclaim}
        \end{claim}
        
        Therefore, it is enough that $\tilde{\beta}$ satisfy 
        \begin{equation}\label{eq:definition of new beta}
            \int_{\sigma} \tilde{\beta}(v) = \int_{\sigma}\beta(v) - \inner{b,\sigma}^*(v)
        \end{equation}
        since then
        \begin{align*}
            \int_{\sigma \times \set{1}}\mf{B}(v) - \int_{\sigma \times \set{0}}\mf{B}(v) 
                & = \int_{\boundary(\sigma \times I)}\mf{B}(v) \tag{since $\mf{B}(v)\res{(\boundary\sigma) \times I} = 0$} \\
                & = \int_{\sigma \times I}d\mf{B}(v) \tag{Stokes' Theorem} \\
                & = -\inner{b,\sigma}^*(v) \tag{\autoref{claim:integral}}\\
                & = \int_{\sigma}\tilde{\beta}(v) - \int_{\sigma}\beta(v). \tag{Property \eqref{eq:definition of new beta}}
        \end{align*}
        
        Concretely, define $\tilde{\beta}\res{\sigma}$ be the volume form of $\sigma$ times a bump function, scaled so that \eqref{eq:definition of new beta} is satisfied.
        This completes the definition of $\Psi(v)\res{\skel{X}{n-1} \times I}$, and it remains to define $\Psi(v)\res{\skel{X}{n}}$.
        Let $\Sigma \sse \skel{X}{n}$ be an $n$-simplex.
        By \Cref{lem:second quantitative poincare}, there is an extension $\tilde{\beta}(v)\res{\Sigma}$ of $\tilde{\beta}(v)\res{\boundary \Sigma}$ to all of $\Sigma$ for which $\|\tilde{\beta}\|_{op} \le O(1)$, and such that $d\tilde{\beta} = \phi(v) - \tilde{f}^*m_Y(v) - \aint_0^1 \Phi(dv)$.
        Define
        \begin{equation*}
            \tilde{\Phi}(v) \defeq \tilde{f}^*m_Y(v) + d(\tilde{\beta}(v) \tensor s) + \aint_0^s \Phi(dv).
        \end{equation*}
        Therefore, we have
        \[
            \Dil[1]{\tilde{\Phi}\res{\skel{X}{n}}} \le \Lip{\tilde{f}} + O(1) + 1 \lesssim 1.
        \]
        
        Finally, we determine $\tilde{\Phi}$ on higher dimensional skeleta.
        Let $\simplex$ be a $k$-simplex for some $k > n$. 
        We inductively apply \Cref{lem:second quantitative poincare} (swapping the roles of $n$ and $k$) to $\alpha_{\boundary} = \tilde\beta(v)\res{\boundary \simplex}$ and $\omega = \phi(v) - \tilde{f}^*m_Y(v) - \aint_0^1 \Phi(dv)$ to produce an extension $\tilde{\beta}(v)\res{\simplex}$ of $\tilde{\beta}(v)\res{\boundary \simplex}$, for which $d\tilde{\beta}(v) = \omega(v)$ with $\|\tilde{\beta}(v)\res{\simplex}\|_{\infty} \le O(1)$.
        Since $v$ was an arbitrary unit vector this induces a definition of $\tilde{\beta} : V \to \Forms{\simplex}$ for which $\|\tilde{\beta}\|_{op} \le O(1)$.
        Define $\tilde{\Phi}(v) = \tilde{f}^*m_Y(v) + d(\tilde{\beta}(v) \tensor s) + \aint_0^{s} \Phi(dv)$ as before. 
        Since $k > n$ there are no obstructions to extending $\Psi$ to a homotopy from $\Phi$ to $\tilde{\Phi}$ over these cells, completing the proof.
    \end{proof} 
\end{lemma}

We put all of these pieces together to get the following theorem, which we state for simplicial sets, and in particular can be applied to compact Riemannian manifolds:

\begin{theorem}[The Shadowing Principle for Lifts of Lipschitz Maps]\label{thm:lifting shadows}
    Let the following diagram commute:
    \begin{equation*}
        \begin{tikzcd}[ampersand replacement = \&]
        A \dar[hook]{\iota} \rar{f\res{A}} 
            \& Y \dar[two heads]{p} \\
        X \urar{f} \rar{\bar{f}}
            \& B
        \end{tikzcd}
    \end{equation*}
    where $p$ is a principal $F = K(G,n)$ fibration of compact Riemannian manifolds with boundary, and where $(X,A)$ is a finite $N$-dimensional simplicial pair.
    Assume that $\Lip{f\res{A}} \le L$ and that $\Lip{\bar{f}} \le L$.
    Moreover, suppose that there is a DGA map $\phi : \minimal{B} \tensor \wedge V \to \Forms{X}$ with $\phi\res{\minimal{B}} = \bar{f}^*m_B$, and $\Dil{\phi} \le L$ such that $\phi \homotopic f^* m_Y$ by a homotopy $\Phi$ with $\Phi\res{\minimal{B}} = \bar{f}^*m_B \tensor 1$.
    
    Then $f$ is homotopic $\rel{A}$ to a $C(L+1)$-Lipschitz map $\tilde{f}$ for which $p \circ \tilde{f} = \bar{f}$, and such that $\phi \homotopic \tilde{f}^*m_Y$ by a homotopy $\tilde{\Phi} : \minimal{B} \tensor \wedge V \to \Forms{X} \trtdt$ with $\Dil[1/L]{\tilde{\Phi}} \le C(L+1)$.
    Here, $C$ is a constant depending only on $\map{p}{Y}{B}$, on the minimal model of the fibration, on the norms of $V$, and on $N$ (but otherwise is independent of $X$).
    
    \begin{proof}
        Let $\bar{\phi} = \phi\res{\minimal{B}} = \bar{f}^*m_B$ and $\bar{\Phi} = \bar{\phi} \tensor 1$.
        After taking a (regular) subdivision of $(X,A)$ at scale $1/L$, we may assume without loss of generality that $L = 1$, and that each simplex of $X$ is uniformly bilipschitz to the standard simplex. 
        At the cost of slightly increasing the Lipschitz constant, and a slight adjustment to the DGA map $\phi\res{\openstar{A}}$, we may assume without loss of generality (see the beginning of \cite[Theorem 4.1]{Manin_2019}) that $\alpha$ is mosaic with respect to some fixed shard complex $\mc{Z}$ on $Y$ for which $\mc{Z}^{(0)}$ consists of basepoints over fibers of vertices in $B$.
        By \Cref{cor:no lower obstructions}, $f$ is internally equivalent to a map $f'$ for which $f'\res{\skel{X}{n-1}}$ is $C_{n-1}$-Lipschitz.
        Moreover, by an inductive application of \Cref{lem:second quantitative poincare} (as in the end of \Cref{lem:local algebraic shadows}) we can assume without loss of generality that $\Dil[1]{\Phi\res{\skel{X}{n-1}}} \le C_{n-1}$.
        By \Cref{lem:local fibered shadows}, there is an internal homotopy $H$ to a $C_{N}$-Lipschitz map $\tilde{f}$ for which $H$ and $\tilde{f}$ meet the assumptions of \Cref{lem:local algebraic shadows}.
        This yields a quadruple $(\bar{\phi},\phi,\bar{\Phi},\tilde{\Phi})$ for $(\bar{f},\tilde{f})$ which is totally $\lesssim C \defeq C_N$ bounded, as desired.
    \end{proof}
\end{theorem}

An inductive argument allows us to replace this single fibration with fibrations which admits a refinement of principal $K(G,n)$ fibrations.
In particular, we recover the original shadowing principle for simply connected targets, but more generally this applies to nilpotent spaces.

\begin{theorem}[The Shadowing Principle for Nilpotent Targets]\label{thm:nilpotent shadows}
    Let $(X,A)$ be an $N$-dimensional finite simplicial pair, and let $Y$ be a compact nilpotent space of finite type.
    Let $f : X \to Y$ be a map and let $\phi : \minimal{Y} \to \Forms{X}$ be such that 
    \begin{enumerate}
        \item $f^* m_Y \res{A} = \phi\res{A}$,
        \item $f^* m_Y \homotopic \phi \rel{A}$,
        \item $\Lip{f\res{A}} \le L$, and
        \item $\Dil{\phi} \le L$.
    \end{enumerate}
    Then $f$ is homotopic $\rel{A}$ to a $(\lesssim L)$-Lipschitz map $\tilde{f} : X \to Y$ such that $\tilde{f}^*m_Y \homotopic \phi$ by a homotopy $\tilde{\Phi} : \minimal{Y} \to \Forms{X} \trtdt$ such that $\tilde{\Phi}\res{A}$ is constant, and $\Dil[1/L]{\tilde{\Phi}} \lesssim L$.
    \begin{proof}
        One simply applies \Cref{thm:lifting shadows} inductively up the Postnikov tower of $Y$.
        The only serious adjustment we must make is that, in order to continue induction at any stage, we require the homotopies at each stage in the Postnikov tower to extend into all of $Y$ both geometrically, and algebraically.

        At any stage in the inductive process, we end up with a map $\mapdegstep{\tilde{g}}{K}{J+1} : X \to \postspacestep{Y}{K}{J+1}$ of bounded Lipschitz constant.
        There is no topological obstruction to lifting this map up the remainder of the tower to produce a genuine map $g_{K}^{J+1} : X \to Y$.
        We similarly end up with a bounded homotopy $\mapdegstep{\tilde{\Psi}}{K}{J+1} : \minimal{\postspacestep{Y}{K}{J+1}} \to \Forms[]{X} \tensor \RR\gen{t,dt}$ from $\mapdegstep{\phi}{K}{J+1}$ to $(\mapdegstep{\tilde{g}}{K}{J+1})^*m_{\postspacestep{Y}{K}{J+1}}$.
        To extend this to a homotopy $\Psi_{K}^{(J+1)} : \minimal{Y} \to \Forms[]{X} \tensor \RR\gen{t,dt}$ we apply \Cref{prop:relative qualitative obstruction} to the diagram 
        \begin{diagram}
            \minimal{\postspacestep{Y}{K}{J+1}} \arrow[rr, "\mapdegstep{\tilde{\Psi}}{K}{J+1}"] \dar[hook]
            \&
            \& \Forms{X \times I} \tensor \RR\gen{s,ds} \dar{t = 0} \drar{}
            \& \\
            \minimal{Y} \arrow[rr,swap,"\Psi_{K}^{J}"] \arrow[urr,dashed]
            \&
            \& \Forms{X} \tensor \RR\gen{s,ds} \rar[equal] 
            \& \Forms{X} \tensor \RR\gen{s,ds} 
        \end{diagram}
        in which the middle vertical map is a quasi-isomorphism.
        Once $n = N + 1 = \dim{X} + 1$, the map $\mapdegstep{g}{N+1}{0}$ may therefore be identified with $\map{g}{X}{Y}$ itself, and so we take $\tilde{f} = \mapdegstep{g}{N+1}{0}$ and $\tilde{\Phi} = \mapdegstep{\Psi}{N+1}{0}$, and the result follows.
    \end{proof}
\end{theorem}

\section{Upper Bounds on Volumes of Nullhomotopies}\label{sec:upper bounds}
Throughout this section, we implicitly work with the minimal model $\minimal{Y} = \wedge \pidual$ of a space $Y$ which is $c$-step nilpotent (up through some finite stage), where $\pidual = \set{\pindual{n}}_{n \ge 1}$.
Since $Y$ is $c$-step nilpotent, we use our notation in \Cref{sec:background}. 
Recall then that $\pindual{n} = \pindualupto{n}{c}$, where $\pindualupto{n}{J} = \Dirsum_{1 \le j \le J} \pindualstep{n}{j}$.
Moreover, we recall that we have a sequence of elementary extensions
\begin{equation*}
    \minimal{\postspacestep{Y}{n}{J-1}} \subseteq \sulext{\minimal{\postspacestep{Y}{n}{J-1}}}{\pindualstep{n}{J}} = \minimal{\postspacestep{Y}{n}{J}}
\end{equation*}
dual to $\postmapstep{n}{J} : \postspacestep{Y}{n}{J} \to \postspacestep{Y}{n}{J-1}$, and therefore that in particular the extension $\minimal{\postspace{Y}{n-1}} \subseteq \minimal{\postspace{Y}{n}}$ passes through a refinement of $c$-many elementary extensions.

Equipped with \Cref{thm:nilpotent shadows}, the same proof as in \cite[\textsection 5]{Manin_2019} leads to the following theorem. 
We omit the proof for brevity, since the only novelty is the application of the Shadowing Principle in the nilpotent setting.

\begin{theorem}[{\cite[Thm. 5.7]{Manin_2019}} + \Cref{thm:nilpotent shadows}]\label{thm:homotopy length and thickness}
    Let $X$ be a finite $n$-dimensional simplicial complex with the standard metric and let $Y$ be a nilpotent space of finite type.
    Let $f : X \to Y$ be a nullhomotopic Lipschitz map and let 
    \begin{equation*}
        \Phi : \minimal{Y} \to \Forms{X} \trtdt 
    \end{equation*}
    be an algebraic nullhomotopy of $f^*m_Y$ with $\Dil[\tau]{\Phi} \le \sigma$ for some $\tau,\sigma \ge 0$.
    There are constants $\kappa = \kappa(X,Y)$ depending only on $X$ and $Y$, and $C = C(Y,m_Y,n)$ depending in addition on $m_Y$ such that there is a nullhomotopy of $f$ with length $C(L_{\lambda} + 1)$ where 
    \begin{equation*}
        L_{\lambda} = \sigma \tau + 1
    \end{equation*}
    and thickness $C(L_{\theta} + 1)$ where 
    \begin{equation*}
        L_{\theta} = \max\set{\sigma,\kappa}.
    \end{equation*}
\end{theorem}

Using this, we will get upper bounds on the volumes of nullhomotopies in nilpotent spaces.
In particular, the goal of \Cref{ssec:nilpotent bounds} will be to prove the following proposition:

\begin{restatable}{proposition}{propnilpotentbounds}\label{prop:nilpotent lipschitz bounds}
    Let $X$ be a finite $N$-dimensional simplicial complex with the standard metric.
    Let $\phi : \minimal{Y} \to \Forms{X} \trtdt$ have $\Dil{\phi} \le L$, and let $\Phi$ be a nullhomotopy which can be built formally up through degree $N$ without encountering any nonzero obstructions that make extendability dependent on choices made in lower degrees. 
    Then writing $\Phi_i^j$ for the $t^idt^j$ term of $\Phi$, for all $1 < n \le N$ we have 
    \begin{equation*}
        \|\Phi_i^j|_{\pindual{n}}\|_{op} \le C(c,n,Y) L^{(n-1)(4c-1) - c + 1} \lesssim L^{n(4c-1) - 4c + 1}
    \end{equation*}
    while for $n = 1$ we have 
    \begin{equation*}
        \|\Phi_i^j|_{\pindual{1}}\|_{op} \le C(c,Y) L^{c}
    \end{equation*}
    In particular, $\Dil[L^{-4c+1}]{\Phi} \le CL^{4c-1}$ for $N > 1$ and $\Dil[L^{c-1}]{\Phi} \le CL$ for $N = 1$.
\end{restatable}

This fact is the analogue of \cite[eq. 3.10]{Manin_2019}, and therefore once we have proved this fact the proof of \cite[Thm. 5.5]{Manin_2019} yields the following:

\begin{theorem}\label{thm:nilpotent upper bounds}
    Let $Y$ be a $c$-step nilpotent space up through dimension $n$, and let be $X$ a finite simply connected complex of dimension $n$. 
    Then there is a constant $C = C(X,Y)$ such that any nullhomotopic $L$-Lipschitz map $f : X \to Y$ is nullhomotopic via a $C$-length nullhomotopy of thickness $O(L^{4c-1})$. 

    \begin{proof}[Proof (Assuming \Cref{prop:nilpotent lipschitz bounds})]
        It suffices to show that we can build a nullhomotopy $\Phi$ by formal extension in each degree without encountering any nontrivial obstructions, thus putting us in the setting of \Cref{prop:nilpotent lipschitz bounds}.
        The proof of this is identical to that of \cite[Thm. 5.5]{Manin_2019}, and we omit the details for brevity. 
        The crux of the argument is that at each stage of the minimal model the morphism
        \begin{equation*}
            [\minimal{Y}(r+1),\Forms{X} \tensor\gen{e}]_0 \to [\minimal{Y}(r),\Forms{X} \tensor\gen{e}]_0 
        \end{equation*}
        is surjective so that obstructions to formally extending nullhomotopies vanish.
        That is, we can take $\Phi : \RR = \minimal{Y}(0) \to \Forms{X} \trtdt$ to be the trivial map, and inductively apply \Cref{prop:quantitative relative obstruction} to formally extend $\Phi$ from $\minimal{Y}(r)$ to $\minimal{Y}(r+1)$ without obstruction.
        We do this in such a way as to construct $\Phi$ as in \Cref{prop:nilpotent lipschitz bounds} for which 
        \begin{equation*}
            \Dil[L^{-4c + 1}]{\Phi} = O(L^{4c-1}).
        \end{equation*} 
        Applying \Cref{thm:homotopy length and thickness} the length of such a nullhomotopy is constant and has thickness bounded by $O(L^{4c-1})$.
    \end{proof}
\end{theorem}

\begin{corollary}\label{cor:nilpotent isolip upper bounds}
    If $Y$ is $c$-step nilpotent up through dimension $n > 1$, then $\isolipverbose{n}{Y}{L} = O(L^{(4c-1)n})$.
\end{corollary}
Note in the application of \Cref{prop:nilpotent lipschitz bounds}, we may have actually reduced the length of the homotopy by considering the slightly stronger bound $\|\Phi_i^{j}\res{\pindual{n}}\|_{op} \lesssim L^{(n-1)(4c-1) - c + 1}$, but this would lead to a less symmetric result.
This lack of symmetry suggests these bounds might be improved under more careful techniques.
For example, a more careful analysis of the situation when $c = 1$ (see \Cref{ssec:simple bounds}) yields the following:
\begin{restatable}{proposition}{propsimplebounds}\label{prop:simple lipschitz bounds}
    Let $X$ be a finite $N$-dimensional simplicial complex with the standard metric.
    Let $c = 1$ so that $Y$ is a simple space up through degree $n$.
    Let $\phi : \minimal{Y} \to \Forms{X} \trtdt$ have $\Dil{\phi} \le L$, and let $\Phi$ be a nullhomotopy which can be built formally up through degree $n$ without encountering any nonzero obstructions that make extendability dependent on choices made in lower degrees. 
    Then writing $\Phi_i^j$ for the $t^idt^j$ term of $\Phi$, for all $1 < n \le N$ we have 
    \begin{equation*}
        \|\Phi_i^j|_{\pindual{n}}\|_{op} \le C(n,Y) L^{2n-1}.
    \end{equation*}
    In particular, $\Dil[L^{-1}]{\Phi} \le CL^{2}$.
\end{restatable}

Applying this in the same manner as above, we can construct a homotopy with length $O(L)$ and thickness $O(L^2)$, and therefore attain the following:
\begin{corollary}\label{cor:simple isolip upper bounds}
    If $Y$ is simple up through dimension $n$, then $\isolipverbose{n}{Y}{L} = O(L^{2n+1})$.
\end{corollary}

If the differential $d$ of $(\minimal{Y},d)$ is quadratic, $Y$ is said to be coformal, and we can further refine our algebraic analysis. 
Consequently, we have the following:
\begin{restatable}{proposition}{propcoformalbounds}\label{prop:coformal lipschitz bounds}
    Let $X$ be a finite $N$-dimensional simplicial complex with the standard metric.
    Let $Y$ be a coformal, $c$-step nilpotent space up through degree $n$.
    Let $\phi : \minimal{Y} \to \Forms{X} \trtdt$ have $\Dil{\phi} \le L$, and let $\Phi$ be a nullhomotopy which can be built formally up through degree $n$ without encountering any nonzero obstructions that make extendability dependent on choices made in lower degrees. 
    Then writing $\Phi_i^j$ for the $t^idt^j$ term of $\Phi$, for all $1 < n \le N$ we have 
    \begin{equation*}
        \|\Phi_i^j|_{\pindual{n}}\|_{op} \le C(n,Y) L^{(c+1)(n-1)}.
    \end{equation*}
    In particular, $\Dil[L^{-(c+1)}]{\Phi} \le CL^{c+1}$.
\end{restatable}

As above, this allows us to conclude the following, which we demonstrate in \Cref{sec:coformal-example} to be (nearly) sharp:

\begin{corollary}\label{cor:coformal isolip upper bounds}
    If $Y$ is $c$-step nilpotent and coformal then $\isolipverbose{n}{Y}{L} = O(L^{(c+1)n})$.
\end{corollary}

The following subsections are largely a technical analysis of the structure of the differential on Sullivan models. 
This is purely algebraic information about the structure of Sullivan models which will provide the foundation of the quantitative results we are interested in.

The rest of this section is organized as follows: in \Cref{ssec:filtrations} we define a filtration on each $\pindual{n}$ from which we will deduce that any nilpotence in the Sullivan model is quantitatively recognized by the action of the abelianization of the fundamental group. 
In \Cref{ssec:nilpotent bounds}, we will analyze the quantitative bounds to prove \Cref{thm:nilpotent upper bounds}. 
In \Cref{ssec:simple bounds} we refine our approach when the target $Y$ is a simple space, with more careful analysis of the differential. 
We use this to prove \Cref{prop:simple lipschitz bounds} which results in \Cref{cor:simple isolip upper bounds}.
Similarly, in \Cref{ssec:coformal bounds} we demonstrate the bounds of \Cref{cor:coformal weights} via \Cref{prop:coformal lipschitz bounds} for coformal spaces.

\subsection{Filtrations}\label{ssec:filtrations}
The main quantity we analyze is the weight of indecomposables of a minimal model.
This will give us a measure of how the differential interacts with certain indecomposables.

\begin{definition}
    Let $\omega \in \pindual{n}$ be indecomposable of degree $n$ for a minimal Sullivan DGA $(\wedge V,d)$.
    Fix a basis $\set{\basis{a}_{i}}$ on $V$.
    The \term{weight} of $\omega$ is defined by $\wt(\omega) \defeq |\omega| = n$ if $d\omega = 0$, and otherwise, writing 
    \begin{equation*}
        d\omega = \sum_{i} r_i (\basis{a}_{i_1} \wedge \basis{a}_{i_2} \wedge \cdots \wedge \basis{a}_{i_{R_i}})
    \end{equation*}
    as a linear combination of products of indecomposables, we inductively define
    \begin{equation*}
        \wt(\omega) = \max_{i}\set{\sum_{k}\wt(\basis{a}_{i_k})}.
    \end{equation*}
\end{definition}
    
We may further define the weight of a product of indecomposables to be the maximal sum of weights of the individual factors.
Note that though we started with a basis for $V$, the weight is well-defined independent of the specific basis we started with.
Later we will set particular bases and conventions of notation.

In order to perform analysis on the weights of indecomposables, we present some helpful filtrations on the Sullivan DGAs in terms of the differential $d$.
The indexing of \cite[\textsection 2]{felix2015rational} (from whom we borrow our notation) begins at $J = 0$.
Unfortunately, this conflicts with the convention of initializing the indices of the central series at $1$, to which we adhere.
This collision of conventions is somewhat inevitable.
We stick with convention---at the cost of this collision---and index starting at $J = 1$.

\begin{definition}
    The \term{cautious filtration} on $\pindual{n}$ is $\cautious{n}{J} \subseteq \cautious{n}{J+1}$ defined by 
    \begin{equation*}
        \cautious{n}{J} \defeq \begin{cases}
            \pindual{n} 
                \cap d^{-1}(\wedge^{\ge 2} \pindual{< n}) & J = 1 \\
            \pindual{n} 
                \cap d^{-1}\l(\wedge^{\ge 2} \pindual{< n} 
                    \dirsum \l(\cautious{1}{J-1} \wedge \cautious{n}{J-1}\r)\r) 
                & J > 1.
        \end{cases}
    \end{equation*}
    The \term{naive filtration} on $\pindual{n}$ is $\naive{n}{J} \subseteq \naive{n}{J+1}$ defined by 
    \begin{equation*}
        \naive{n}{J} \defeq \begin{cases}
            \pindual{n} 
                \cap d^{-1}(\wedge^{\ge 2} \pindual{< n}) & J = 1 \\
            \pindual{n} 
                \cap d^{-1}\l(\wedge^{\ge 2} \pindual{< n} 
                    \dirsum \l(\pindual{1} \wedge \naive{n}{J-1}\r)\r) 
                & J > 1.
        \end{cases}
    \end{equation*}
    That is, one may always pick $\pindual{n}^{\le J} = \naive{n}{J}$, but it may be possible to find a more refined definition of $\pindual{n}^{\le J}$.
    Moreover, having fixed a norm on $\pindual{n}$ we define the $\cautiousstep{n}{J}$ (resp. $\naivestep{n}{J}$) the \term{$J$th cautious step} of the cautious filtration by
    \begin{equation*}
        \cautiousstep{n}{J} 
            \defeq \begin{cases} 
                \cautious{n}{1} & J = 1 \\
                \l(\cautious{n}{J-1}\r)^{\perp} \cap \cautious{n}{J} & J > 1.
            \end{cases}
    \end{equation*}
    and the \term{$J$th naive step} of the naive filtration by
    \begin{equation*}
        \naivestep{n}{J} 
            \defeq \begin{cases} 
                \naive{n}{1} & J = 1 \\
                \l(\naive{n}{J-1}\r)^{\perp} \cap \naive{n}{J} & J > 1.
            \end{cases}
    \end{equation*}
\end{definition}
In particular observe that
\begin{equation*}
    \cautiousstep{n}{J+1} \dirsum \cautious{n}{J} = \cautious{n}{J+1}
\end{equation*}
and 
\begin{equation*}
    \naivestep{n}{J+1} \dirsum \naive{n}{J} = \naive{n}{J+1}.
\end{equation*}

\begin{note}
    The distinction in the naive and cautious filtrations is subtle, but important. 
    In the naive filtration, we essentially say that ``we know $\pi_1$ acts on our space, and we don't care \textit{which} elements are doing the action'', and so we take all of the action all at once.
    In the cautious filtration, we move more slowly, considering only those pieces in the $J$th step of $\pi_1$ itself, when building the $(J+1)$st stage of the filtration.
    Because of this, the cautious filtration might take twice as long as the naive one to terminate.
    The main result of this section is that in fact the two filtrations agree.
\end{note}

We begin with the easy direction:
\begin{proposition}\label{prop:naive is bigger}
    $\pindualupto{n}{J} \sse \naive{n}{J}$ and $\cautious{n}{J} \sse \naive{n}{J}$ for all $n$ and all $J$.
    \begin{proof}
        For $J = 1$, by definition we have $\im{d\res{\pindualupto{n}{1}}} \sse \wedge^{\ge 2} \pindual{< n}$, and so obviously 
        \begin{equation*}
            \pindualupto{n}{1} \sse \pindual{n} \cap d^{-1}\l(\wedge^{\ge 2}\pindual{< n}\r) = \cautious{n}{1}.
        \end{equation*}
        Moreover, clearly by definition we have $\cautious{n}{1} = \naive{n}{1}$.
        Inductively then we have 
        \begin{equation*}
            \im{d\res{\pindualupto{n}{J}}} 
                \sse \wedge^{\ge 2} \pindual{< n} \dirsum \l(\pindual{1} \wedge \pindualupto{n}{J-1}\r)
                \sse \wedge^{\ge 2} \pindual{< n} \dirsum \l(\pindual{1} \wedge \naive{n}{J-1}\r)
        \end{equation*}
        so that 
        \begin{equation*}
            \pindualupto{n}{J} 
                \sse \pindual{n} \cap d^{-1}\l(\wedge^{\ge 2} \pindual{< n}
                    \dirsum \l(\pindual{1} \wedge \naive{n}{J-1}\r)\r) 
                = \tilde{W}_n^{J}.
        \end{equation*}
        Moreover, since $\cautious{1}{J} \sse \pindual{1}$ for all $J$ we have 
        \begin{align*}
            \cautious{n}{J} 
                & = \pindual{n} 
                    \cap d^{-1}\l(\wedge^{\ge 2} \pindual{< n} 
                        \dirsum \l(\cautious{1}{J-1} \wedge \cautious{n}{J-1}\r)\r) \\
                & \sse \pindual{n} 
                    \cap d^{-1}\l(\wedge^{\ge 2} \pindual{< n} 
                        \dirsum \l(\pindual{1} \wedge \cautious{n}{J-1}\r)\r) \\
                & \sse \pindual{n} 
                    \cap d^{-1}\l(\wedge^{\ge 2} \pindual{< n} 
                        \dirsum \l(\pindual{1} \wedge \naive{n}{J-1}\r)\r) \\
                & = \naive{n}{J}.
                \qedhere
        \end{align*}
    \end{proof}
\end{proposition}

\begin{corollary}
    $\naive{n}{c} = \pindual{n}$.
    \begin{proof}
        By the previous proposition
            $\pindual{n} = \pindualupto{n}{c} \sse \naive{n}{c} \sse \pindual{n}$.
    \end{proof}
\end{corollary}

Now let $\omega \in \CS{n}{J+1}$.
Then $d\omega \in \wedge^{\ge 2}\pindual{< n} \dirsum (\cautious{1}{J} \wedge \cautious{n}{J})$. 
Let $d = \scomp + \ncomp$ where
\begin{equation*}
    \scomp 
        : \CS{n}{J+1} \to \wedge^{\ge 2}\pindual{< n}
\end{equation*}
and
\begin{equation*}
    \ncomp 
        : \CS{n}{J+1} \to \cautious{1}{J} \wedge \cautious{n}{J}.
\end{equation*}
We call $\scomp(\omega)$ the \term{simple component}, and $\ncomp(\omega)$ the \term{nilpotent component} of $d\omega$. 
Moreover, for $n > 1$ we may write 
\begin{equation*}
    \cautious{1}{J} \wedge \cautious{n}{J} = \Dirsum_{1 \le i,j \le J} \CS{1}{i} \wedge \CS{n}{j}
\end{equation*}
and we write 
\begin{equation*}
    \ncompstep{i,j} : \CS{n}{J+1} \to \CS{1}{i} \wedge \CS{n}{j}
\end{equation*}
as the $(i,j)$-nilpotent component of the differential.
We may similarly define $\scomp$ and $\ncomp$ on $\pindualupto{n}{J}$ in the obvious way.

Let $\omega \in \CS{n}{J}$. We define $\abcomp\omega = \ncompstep{1,J-1}(\omega)$, and write $\abcomp = \ncompstep{1,J-1}$, which we call the \term{abelian component} of the differential.
Context should always make clear the domain of the abelian component.
This terminology is inspired by \cite{felix2015rational}, where the orthogonal subspaces of $\im{d}$, and summands therein, are called ``components.''

A number of the subsequent proofs will be based on certain maximality arguments.
To aid with this, we define a lexicographical order on tuples of natural numbers.
Many of the methods used here are essentially those found in \Cref{lem:cautious step higher degrees}, but perhaps a bit clarified by this notation.

\begin{definition}
    For $K \in \NN$, let $U_K = \set{(m_1,\dots,m_K) \in \NN^{K} \sothat m_i \le m_{i+1}}$. 
    For $K = 1$, define $m_1' \len{1} m_1 \iff m_1' \le m_1$. 
    Inductively, define $(m_1',\dots,m_K') \len{K} (m_1,\dots,m_K)$ if and only if either $m_K' < m_K$ or 
    \begin{equation*}
        m_K' = m_K \;\text{ and }\; (m_1',\dots,m_{K-1}') \len{K-1} (m_1,\dots,m_{K-1}).
    \end{equation*}
    In the event that $(m_1',\dots,m_K') \len{K} (m_1,\dots,m_K)$ and $(m_1',\dots,m_K') \neq (m_1,\dots,m_K)$ we write $(m_1',\dots,m_K') \lessn{K} (m_1,\dots,m_K)$.
\end{definition}

The following is a result of \cite{felix2015rational}, but note that by our reindexing, the exact statements look slightly different.
\begin{lemma}[{\cite[Lem 2.3]{felix2015rational}}]\label{lem:cautious step degree 1}
    The nilpotent component of the differential on $\CS{1}{J}$ satisfies 
    \begin{equation*}
        \ncomp = \sum_{i + j \le J} \ncompstep{i,j}.
    \end{equation*}
    That is, if $\omega \in \CS{1}{J}$ then $d\omega = \sum_{i + j \le J}\Omega_{i,j}$ where $\Omega_{i,j} \in \CS{1}{i} \wedge \CS{1}{j}$. Moreover, $\abcomp$ is injective.  
\end{lemma}

We generalize this to the following:

\begin{lemma}\label{lem:cautious step higher degrees}
    For all $n$ and $J$, the nilpotent component of the differential on $\CS{n}{J}$ satisfies \begin{equation*}
        \ncomp = \sum_{i + j \le J} \ncompstep{i,j}.
    \end{equation*}
    In particular, if $\omega \in \CS{n}{J}$ is nonzero then  
    \begin{equation*}
        d\omega 
            = \scomp\omega + \ncomp\omega 
            = \scomp\omega + \sum_{i + j \le J}\Omega_{ij}
    \end{equation*}
    where each $\Omega_{ij} \in \CS{1}{i} \wedge \CS{n}{j}$.
    Moreover, the morphism $\abcomp : \CS{n}{J} \to \CS{1}{1} \wedge \CS{n}{J-1}$ is injective for all $J \ge 2$.
    \begin{proof}
        When $n = 1$ this is the result of \Cref{lem:cautious step degree 1}, so let $n > 1$.
        First consider $J = 2$ and let $\omega \in \CS{n}{2}$. 
        Note that in this case $\abcomp\omega = \ncomp\omega$.
        Thus if $\abcomp \omega = 0$ then $\ncomp\omega \in \CS{1}{1} \wedge \CS{n}{1}$ is trivial.
        But by definition $\omega \in d^{-1}(\wedge^{\ge 2} \pindual{< n}) = \CS{n}{1} = \cautious{n}{1}$, and since $\CS{n}{2}$ is orthogonal to $\cautious{n}{1}$ we have
        \begin{equation*}
            \omega \in \cautious{n}{1} \cap \CS{n}{2} = \set{0}.
        \end{equation*}
        Therefore the result holds for $J = 2$.

        Now let $J > 2$ and assume the result for all $2 \le j < J$.
        Let $\omega \in \CS{n}{J}$.
        Note that $0 = d^2\omega = (d\circ\scomp)\omega + (d\circ\ncomp)\omega$.
        Moreover, note that 
        \begin{equation*}
            \im{d\circ\scomp} \sse \wedge^{\ge 3}\pindual{< n}
        \end{equation*}
        and that 
        \begin{equation*}
            \im{d \circ \ncomp} \sse 
                (\wedge^{\ge 3}\pindual{< n}) 
                \dirsum \left((\cautious{1}{J-2} \wedge \cautious{1}{J-2} \wedge \cautious{n}{J-1}) 
                + (\cautious{1}{J-1} \wedge \cautious{1}{J-2} \wedge \cautious{n}{J-2})\right).
        \end{equation*}
        In particular, the components of $(d \circ \ncomp)\omega$ outside of $\wedge^{\ge 3}\pindual{< n}$ must all be trivial, as they cannot cancel with components of $(d \circ \scomp)\omega$.
        To this end, write $\ncomp \omega = \sum_{i,j \le J-1}\Omega_{i,j}$, such that $\Omega_{i,j} \in \CS{1}{i} \wedge \CS{n}{j}$.
        A priori, $i + j$ is unbounded.
        Let $\ell \defeq \max\set{i + j \sothat \Omega_{i,j} \neq 0}$.
        Moreover, fix $(i_0,j_0) \defeq \max_{\len{2}} \set{(i,j) \sothat i+j = \ell ,\; \Omega_{i,j} \neq 0}$.
        In particular, $j_0$ is maximal over such pairs $(i,j)$ for which $\Omega_{i,j} \neq 0$.
        
        Suppose to the contrary that $d\Omega_{i_0,j_0} \neq 0$.
        Write $\Omega_{i_0,j_0} = \sum_{r} \basis{t}_r \wedge \inde{v}_r \in \CS{1}{i_0} \wedge \CS{n}{j_0}$ such that the $\basis{t}_r$ terms are linearly independent.
        By \cite[Lem 2.3]{felix2015rational} (equivalently, by the base case $n = 1$) we know that $\abcomp\res{\CS{1}{i_0}}$ is injective. 
        Hence the collection $\set{\abcomp \basis{t}_r} \sse \CS{1}{1} \wedge \CS{1}{i_0-1}$ is linearly independent.
        Consequently $\sum_r (\abcomp \basis{t}_r) \wedge \inde{v}_r \in \CS{1}{1} \wedge \CS{1}{i_0-1}\wedge\CS{n}{j_0}$ is itself nonzero.

        Therefore, there is some term $\Omega_{i',j'} \in \CS{1}{i'}\wedge\CS{n}{j'}$ of $\ncomp \omega$ for which $d\Omega_{i',j'}$ has a summand $\gamma \in \CS{1}{1} \wedge \CS{1}{i_0-1} \wedge \CS{n}{j_0}$. 
        Moreover, by induction (applying the differential to $\Omega_{i',j'}$) we must have 
        \[
            \gamma \in \left(\CS{1}{a_i} \wedge \CS{1}{b_i} \wedge \CS{n}{j'}\right) \dirsum \left(\CS{1}{i'} \wedge \CS{1}{a_j} \wedge \CS{n}{b_j}\right)
        \]
        for which $a_i + b_i \le i'$ and $a_j + b_j \le j'$.
         
        First, if
        \[
            \gamma \in \left(\CS{1}{a_i} \wedge \CS{1}{b_i} \wedge \CS{n}{j'}\right) \cap \left(\CS{1}{1} \wedge \CS{1}{i_0-1} \wedge \CS{n}{j_0}\right)
        \]
        then (without loss of generality) we have $a_i = 1$ and $b_i = i_0 -1$ and $j' = j_0$. 
        But $i_0 = a_i + b_i \le i'$ and $j' = j_0$.
        By maximality of $(i_0,j_0)$ this implies that in fact $(i',j') = (i_0,j_0)$ and hence that $\Omega_{i',j'} = s\Omega_{i_0,j_0}$ for some $s \in \RR$.
        
        On the other hand, if 
        \[
            \gamma \in \left(\CS{1}{i'} \wedge \CS{1}{a_j} \wedge \CS{n}{b_j}\right) \cap \left(\CS{1}{1} \wedge \CS{1}{i_0-1} \wedge \CS{n}{j_0}\right)
        \]
        then $b_j = j_0$.
        But then, since $j_0$ is maximal and since $a_j > 0$, we have
        \[
            j_0 < j_0 + a_j = b_j + a_j \le j'
        \]
        a contradiction to the maximality of $j_0$.
        
        Thus, $d\Omega_{i_0,j_0} = 0$.
        But this implies that $\sum_r (\abcomp \basis{t}_r) \wedge \inde{v}_r = 0$. 
        By linear independence of the $\abcomp \basis{t}_r$ terms we must have $\abcomp \basis{t}_r = 0$ for all $r$, which by injectivity of $\abcomp$ in degree $1$ implies that $\basis{t}_r = 0$ for all $r$.
        In particular we have $i_0 = 1$, and since $j_0 \le J-1$ we have $\ell = i_0 + j_0 = 1 + j_0 \le J$.
        By the maximality of $\ell$ we see that $\ncomp$ is as desired.

        Finally, if $\ell \neq J$, then each $\Omega_{i,j}$ actually belongs to $\cautious{1}{J-2} \wedge \cautious{n}{J-2}$ so that $\omega \in \cautious{n}{J-1}$. 
        Thus if $\omega \in \CS{n}{J}$ as we assumed, we would have had $\omega \in \CS{n}{J} \cap \cautious{n}{J-1} = \set{0}$, a contradiction. 
        Therefore, $\ell = J$.
        The above shows that the component $\Omega_{1,J-1}$ is nontrivial, and thus that $\abcomp$ is injective as claimed.
    \end{proof}
\end{lemma}

This allows us to prove:

\begin{lemma}\label{lem:naive is smaller}
    $\naive{n}{J} = \cautious{n}{J}$ for all $n$ and $J$.
    \begin{proof}
        For $J = 1$, this is clearly true by definition.
        By \Cref{prop:naive is bigger} we know that $\cautious{n}{J} \sse \naive{n}{J}$. 
        Inductively we can write 
        \begin{equation*}
            \cautious{n}{J} = \cautious{n}{J-1} \dirsum \CS{n}{J} = \naive{n}{J-1} \dirsum \CS{n}{J} \sse \naive{n}{J-1} \dirsum \NS{n}{J} = \naive{n}{J}
        \end{equation*}
        where in fact $\CS{n}{J} \sse \NS{n}{J}$ as they are both orthogonal to the same $\cautious{n}{J-1} = \naive{n}{J-1}$ by induction.
        More generally, we have $\CS{n}{K}$ orthogonal to $\cautious{n}{J-1} = \naive{n}{J-1}$ for all $K > J-1$. 
        Let $K$ be maximal so that there is $\omega \in \NS{n}{J}$ for which the natural projection $p : \NS{n}{J} \to \CS{n}{K}$ is nonzero.
        By injectivity of $\abcomp$ we have $\abcomp(p(\omega)) \neq 0$ as a component of $d\omega$ in $\CS{1}{1} \wedge \CS{n}{K-1}$.
        But also, by definition, the differential sends the naive filtration down a degree in the filtration.
        Therefore, since $p(\omega) \in \NS{n}{J}$, we have
        \begin{equation*}
            \abcomp(p(\omega)) \in \pindual{1} \wedge \naive{n}{J-1} = \pindual{1} \wedge \cautious{n}{J-1}
        \end{equation*}
        by induction. 
        Therefore 
        \begin{equation*}
            \abcomp(p(\omega)) \in \l(\pindual{1} \wedge \cautious{n}{J-1}\r) \cap \l(\CS{1}{1} \wedge \CS{n}{K-1}\r)
        \end{equation*}
        is nonzero which in particular implies that $K \le J$.
        Since $K$ was maximal for which this projection was nontrivial, every element of $\naive{n}{J}$ projects trivially onto $\CS{n}{j}$ for $j > J$ so that $\cautious{n}{J} \sse \naive{n}{J} \sse \Dirsum_{j=1}^{J}\CS{n}{j} = \cautious{n}{J}$ as required.
    \end{proof}
\end{lemma}

\begin{corollary}\label{cor:naive is cautious}
    For all $J$ we have $\pindualupto{n}{J} \sse \cautious{n}{J}$.
    In particular
    \begin{equation*}
        \im{\ncomp\res{\pindualupto{n}{J}}} \sse \im{\ncomp\res{\cautious{n}{J}}} \sse \Dirsum_{i + j \le J}\CS{1}{i} \wedge \CS{n}{j}.
    \end{equation*}
\end{corollary}

We work exclusively with the cautious filtration from now on.
Here we establish conventions and terminology which we frequently use throughout the rest of the paper.
In particular, defining a basis of $\CS{n}{J}$ for each $n$ and $J$ induces a basis of $V$ and thus on $\wedge V$ in the obvious ways.
Given such a basis we induce an inner product $\inner{\cdot,\cdot}$ for which this basis is orthonormal.
Given such a basis and inner product, we try to adhere to the following conventions in what follows:

\begin{notation}\text{}\label{notation:indecomposables}
    \begin{itemize}
        \item Standard \textit{italic} Roman letters such as $r,s$ and $i,N$ denote scalars and indices.
        \item $\basis{Bold\text{ }and\text{ }barred}$ lower case Roman letters such as $\basis{x}$ denote basic indecomposables, and in particular are elements of $\pindual{n}$ for some degree $n$.
        \item The lower case Greek letter $\omega$, and \textbf{bold} Roman lower case letters such as $\inde{x}$ denote generic indecomposables, and in particular are elements of $\pindual{n}$ for some degree $n$.
        \item Greek letters such as $\Omega, \alpha, \beta, \xi$ (but \textbf{not} $\omega$; see above) denote scalar multiples of basis elements of $\wedge V$. In particular, these denote basis elements of $\wedge^{k}V$ for some $k$, and are products of basic indecomposables.
        \item Capital Roman letters such as $A,B$ denote generic elements of $\wedge V$.
    \end{itemize}
\end{notation}
For example we might write $\inde{x} \wedge \basis{y} = \left(\sum_{i} r_{i}\basis{x}_{i} \wedge \basis{y}\right)$ to denote an element of $V_{|x|} \wedge \sspan{\basis{y}}$.
We will sometimes be concerned with specifying that certain summands live in complementary subspaces.
Generally then if $A \in \wedge V$ (not necessarily basic) for some $k$, and $B \in \wedge V$ write $A + B = A \pp B$ if $A$ and $B$ live in complementary subspaces.

Given $A \in \wedge V$ we will typically say that $\Omega \in \wedge^{k}V$ is a \term{(basic) summand} of $A$ if $\Omega \in \sspan{\tilde{\Omega}}$ for some basic element $\tilde{\Omega}$ of $\wedge^{k}V$ for which $\inner{A, \Omega} > 0$.
In particular, $\Omega$ is a nonzero scalar multiple of a basic element.
(Note, that $\Omega$ is a basic summand of $d\omega$ if and only if $d\omega = \Omega \pp B$ for some possibly trivial $B$.)
A \term{summand} more generally is a linear combination of basic summands, whose coefficients are typically assumed to be nonzero, but we will often use the word ``summand'' when we really mean a basic summand. 
Context should make the meaning of the term clear.
Finally, suppose that $\xi$ is basic and $A \in \wedge V$.
We say that $\xi$ \term{cancels} with $A$ if $\inner{A,\xi} < 0$; that is, $A$ contains $r\xi$ as a basic summand for some $r < 0$.
The following is a straightforward consequence of the definitions:

\begin{lemma}\label{lem:cancelling property}
    Let $\omega$ be indecomposable, let $\inner{d\omega,A} \neq 0$, and suppose $dA$ has a (nonzero) basic summand $\xi$.
    Then there is a (basic) summand $\Omega$ of $d\omega$ complementary to $A$ for which $\inner{d\Omega,\xi} < 0$.
    In particular, $d\Omega$ has a (basic) summand which cancels with $\xi$.
    \begin{proof}
        Write $d\omega = A \pp B$.
        Since $dA + dB = d^2\omega = 0$ we have 
        \[
            0 = \inner{0,\xi} = \inner{d^2\omega,\xi} = \inner{dA,\xi} + \inner{dB,\xi} 
        \]
        so that $\inner{dB,\xi} = -\inner{dA,\xi}$.
        Since $\xi$ is a basic summand of $dA$ we have $\inner{dA,\xi} > 0$ by definition.
        Since $B$ is complementary to $A$, there must be a (basic) summand ${\Omega}$ of $B$ for which $\inner{d{\Omega},\xi} < 0$.
    \end{proof}
\end{lemma}

\begin{definition}
    If $\xi$ is as in the above lemma we say that $\xi$ is \term{cancelled externally} relative to $A$.
\end{definition}

\begin{definition}
    Let $\Omega$ be a basic summand of $d\omega$, and write $\Omega = r\basis{x}_1 \wedge \cdots \wedge \basis{x}_{k}$.
    Let $\xi$ be a (possibly trivial) summand of $r\basis{x}_{1} \wedge \cdots \wedge (d\basis{x}_{i}) \wedge \cdots \wedge \basis{x}_{k}$ for some $i$.
    We say that $\xi$ is \term{cancelled internally} relative to $\Omega$ if $\xi = 0$ or if $\xi$ is a basic summand which cancels with a summand of $r\basis{x}_1 \wedge \cdots \wedge (d\basis{x}_{j}) \wedge \cdots \wedge \basis{x}_{k}$ for some $j$.
\end{definition}

\begin{note}
    Let $\Omega$ be as above and let $\xi$ be a basic summand of $r\basis{x}_{1} \wedge \cdots \wedge (d\basis{x}_{i}) \wedge \cdots \wedge \basis{x}_{k}$ for some $i$.
    If $\xi$ is not cancelled externally (relative to $\Omega$), then by contraposition of the previous lemma, $\xi$ is \textbf{not} a nonzero summand of $d\Omega$. 
    Hence, $\xi$ must be cancelled internally (relative to $\Omega$).
\end{note}

\begin{lemma}\label{lem:internal}
    Let $\Omega$ be as above and let $\xi$ be a summand of $r\basis{x}_{1} \wedge \cdots \wedge (d\basis{x}_{i}) \wedge \cdots \wedge \basis{x}_{k}$ for some $i$.
    If $\xi$ is cancelled internally (relative to $\Omega$) by a summand of $r\basis{x}_1 \wedge \cdots \wedge (d\basis{x}_{j}) \wedge \cdots \wedge \basis{x}_{k}$ of $d\Omega$, then $\xi = 0$ or $\basis{x}_i \neq \basis{x}_j$.
    \begin{proof}
        For contraposition, suppose $\basis{x}_{i} = \basis{x}_{j}$ and $\xi \neq 0$.
        If $i = j$ then $\xi = 0$, and so we may assume that $i \neq j$.
        If $|\basis{x}_{i}|$ is odd, then $\Omega = 0$ since $\basis{x}_{i} \wedge \basis{x}_{j} = \basis{x}_{i}^{2} = 0$.
        Assume therefore that $|\basis{x}_{i}|$ is even, and without loss of generality $\basis{x}_{i} = \basis{x}_{1} = \basis{x}_{2} = \basis{x}_{j}$.
        Then by graded commutativity we have 
        \[
            (d\basis{x}_{1}) \wedge \basis{x}_{2} \wedge \cdots \wedge \basis{x}_{k} = \basis{x}_{2} \wedge (d\basis{x}_{1}) \wedge \cdots \wedge \basis{x}_{k} = \basis{x}_{1} \wedge (d\basis{x}_{2}) \wedge \cdots \wedge \basis{x}_{k}
        \]
        since $\basis{x}_1 = \basis{x}_{2}$.
        Therefore, by definition of $\xi$ being a basic summand we have
        \[
            0 < \inner{(d\basis{x}_{1}) \wedge \basis{x}_{2} \wedge \cdots \wedge \basis{x}_{k},\xi} = \inner{\basis{x}_{1} \wedge (d\basis{x}_{2}) \wedge \cdots \wedge \basis{x}_{k},\xi}
        \]
        so that $\xi$ is not cancelled by such a summand.
    \end{proof}
\end{lemma}

\begin{lemma}\label{lem:trivial or external}
    Let $N$ be given and let $\Omega = \alpha \wedge \basis{x}^{\ell} \in( \wedge^{\ge 0}V) \wedge (\wedge^{\ell}\CS{N}{J})$ be a basic summand of $d\omega$.
    In particular, $\basis{x}$ is basic, and $\alpha$ is a product of basic indecomposables.
    Moreover, assume that $\ell$ is maximal in the sense that $\alpha$ does not contain $\basis{x}$ as a factor.
    Finally, assume that $J > 1$.
    Let $\basis{t} \wedge \zeta \in \CS{1}{p} \wedge (\wedge^{k}V)$ be a basic summand of $d \basis{x}$, and let $\xi \defeq \alpha \wedge \basis{x}^{\ell-1} \wedge \basis{t} \wedge \zeta$.
    Then either $\xi = 0$ (in which case $\basis{t}$ is a factor of $\alpha$, or $\zeta$ and $\alpha$ share a common factor) or $\xi$ is cancelled externally relative to $\Omega$.
    \begin{proof}
        If $\xi$ is not cancelled externally, then $\xi$ is cancelled internally relative to $\Omega$.
        Suppose to the contrary that $\xi \neq 0$.
        By the above lemma, internal cancellation does not happen within any of the $\basis{x}$ factors, so we must have some basic summand $\beta$ of $d\alpha$ for which $\inner{\beta \wedge \basis{x}^{\ell},\alpha \wedge \basis{x}^{\ell-1} \wedge \basis{t} \wedge \zeta} < 0$.
        There are $\ge \ell$ factors of $\basis{x}$ in $\beta \wedge \basis{x}^{\ell}$, but only $\ell-1$ factors of $\basis{x}$ in $\alpha \wedge \basis{x}^{\ell-1} \wedge \basis{t} \wedge \zeta$ by the maximality of $\ell$ and the fact that $\basis{x}$ is not a factor of $\basis{t} \wedge \zeta$ by degree/step considerations.
        Since $\inner{\cdot,\cdot}$ is defined to make our basis orthogonal, we must have $\inner{\beta \wedge \basis{x}^{\ell},\alpha \wedge \basis{x}^{\ell-1} \wedge \basis{t} \wedge \zeta} = 0$, a contradiction.

        Hence, $0 = \xi = \alpha \wedge \basis{x}^{\ell-1} \wedge \basis{t} \wedge \zeta$. 
        Since each factor is basic indecomposable, we must have a repeated factor of odd degree.
        Since $\Omega$ itself is nonzero, this repeated factor must be a factor of $\basis{t} \wedge \zeta$.
        Since $\basis{t} \wedge \zeta$ is a summand of $d\basis{x}$ this repeated factor cannot be $\basis{x}$ itself, which completes the proof.
    \end{proof}
\end{lemma}

\subsection{Bounds on Nilpotent Spaces}\label{ssec:nilpotent bounds}
We continue using our prior notation from \Cref{ssec:filtrations}.
The main goal of this section is the proof of \Cref{prop:nilpotent lipschitz bounds}.
We now put some quantitative bounds in place on the weights of indecomposables.
First we make an observation about the indecomposables in degree $1$.

\begin{proposition}\label{prop:degree 1 weights}
    If $\omega \in \cautious{1}{J}$ then $\wt(\omega) \le J$.
    In particular, $\wt(\omega) = J$ if $\omega \in \CS{1}{J}$.
    \begin{proof}
        When $J = 1$ by definition $d\omega = 0$ and so $\wt(\omega) = 1$. 
        Now \Cref{lem:cautious step degree 1} implies that any basic summand of $d\omega$ is of the form $\Omega_{i,j} = \basis{x_{i}} \wedge \basis{y_{j}} \in \cautiousstep{1}{i} \wedge \cautiousstep{1}{j}$ where $i + j \le J$, and inductively $\wt(\basis{x}_i) \le i$ and $\wt(\basis{y}_j) \le j$.
        Thus if $\Omega_{i,j}$ maximizes the weight of $\omega$, then $\wt(\omega) \le i + j \le J$.
        Because $\abcomp$ is injective $\abcomp \omega \in \CS{1}{1} \wedge \CS{1}{J-1}$ is nonzero and therefore inductively has weight exactly $J$.
    \end{proof}
\end{proposition}

Naively, one might expect that the indecomposables in $\CS{2}{1}$ of a $c$-step nilpotent space might have weight at most $3c$, but this bound can be sharpened. 
In fact, we see in \Cref{cor:degree 2 weights} that $3c$ ends up being the largest weight of $\pindual{2}$ overall.
We expect the techniques in degree $2$ carry through higher dimensions to produce sharper upper bounds on the weight of indecomposables.
Here, we examine one such technique, with an eye towards \Cref{conj:better nilpotent bounds}.

\begin{definition}
    If $\Omega$ is a summand of $d\omega$ and $\wt(\omega) = \wt(\Omega)$, we say that $\Omega$ is a \term{full weight}, or \term{weight maximizing} summand of $\omega$.
    On the other hand, if $\inde{x} \in \CS{n}{1}$ for some $n$, we say that $\inde{x}$ is \term{light}.
\end{definition}

\begin{lemma}\label{lem:full weight}
    Let $\omega \in \CS{2}{J}$ and suppose that all full weight summands of $\omega$ are summands of $\scomp \omega$.
    That is, $\wt(\omega) = \wt(\scomp\omega) > \wt(\ncomp\omega)$.
    Then some full weight basic summand of $d\omega$ contains a light factor.
    \begin{proof}
        Note that all basic summands of $\scomp \omega$ belong to $\CS{1}{J_1} \wedge \CS{1}{J_2} \wedge \CS{1}{J_3}$ for some $J_1 \le J_2 \le J_3$.
        Let $\Omega = \basis{x} \wedge \basis{y} \wedge \basis{z}$ be a full weight summand of $\scomp \omega$ for which $(J_1,J_2,J_3)$ is maximal under $\len{3}$.
        By \Cref{prop:degree 1 weights} we see in particular that $\wt(\omega) = J_1 + J_2 + J_3$.
        That is, if $\Omega' \in \CS{1}{J_1'} \wedge \CS{1}{J_2'} \wedge \CS{1}{J_3'}$ with $J_1' \le J_2' \le J_3'$ is another full weight summand of $\scomp\omega$, then $(J_1',J_2',J_3') \len{3} (J_1,J_2,J_3)$.
        Furthermore, up to a change of basis we may assume that if $\basis{x}' \wedge \basis{y} \wedge \basis{z}$ is a summand of $d\omega$ for which $(d\basis{x}') \wedge \basis{y} \wedge \basis{z}$ has $(\abcomp \basis{x}) \wedge \basis{y} \wedge \basis{z}$ as a summand, then $\basis{x}' = \basis{x}$.

        Suppose to the contrary $J_1 > 1$, and let $\xi \defeq (\abcomp \basis{x}) \wedge \basis{y} \wedge \basis{z}$.
        Suppose then $\Omega' = \basis{x}' \wedge \basis{y}' \wedge \basis{z}'$ is a summand of $d\omega$ for which $d\Omega'$ has a summand which cancels with $\xi$.
        By \Cref{prop:degree 1 weights} we also see that $\wt(\xi) = \wt(\Omega) = 1 + (J_1-1) + J_2 + J_3 = \wt(\omega)$, and thus $\wt(\Omega') \ge \wt(\xi)$.
        In particular, $\Omega'$ is full weight.
        By maximality of $(J_1,J_2,J_3)$ under full weight summands of $d\omega$, we see that $(J_1',J_2',J_3') \len{3} (J_1,J_2,J_3)$.

        See first that $J_3' \ge J_3$. 
        Otherwise $d\Omega'$ contains no summands in
        $\CS{1}{J_1} \wedge \CS{1}{J_2} \wedge \CS{1}{J_3}$ by consideration step of the filtration.
        By maximality, $J_3' = J_3$, and in fact $\basis{z}' = \basis{z}$.
        Similarly, once $J_3' = J_3$ we see that $J_2' = J_2$ in the same way and $\basis{y}' = \basis{y}$, and hence that $J_1' = J_1$ as well.
        In particular, we must then have $\Omega' = r\basis{x}' \wedge \basis{y} \wedge \basis{z}$ for some $r \in \RR$, where $(d\basis{x}') \wedge \basis{y} \wedge \basis{z}$ has $(\abcomp \basis{x}) \wedge \basis{y} \wedge \basis{z}$ as a summand.
        Hence, $\basis{x}' = \basis{x}$ by our initial assumptions.
        In particular, $\Omega' \in \sspan{\Omega}$.
        Cancellation is thus not external so that \Cref{lem:trivial or external} implies $\xi = 0$.
        But $\xi \in \CS{1}{1} \wedge \CS{1}{J_1-1} \wedge \CS{1}{J_2} \wedge \CS{1}{J_3}$ is trivial if and only if $\abcomp\basis{x} = 0$ by consideration of step of the filtration.
        Since $\abcomp$ is injective, $\basis{x} = 0$, contradicting the fact that $\basis{x} \in \CS{1}{J_1}$ with $J_1 > 1$.
    \end{proof}
\end{lemma}

\begin{lemma}\label{lem:degree 2 weights}
    Let $\omega \in \CS{2}{J}$.
    Then $\wt(\omega) \le J + 2c$.
    \begin{proof}
        Let $\omega \in \CS{2}{J}$.
        If $J = 1$ then $d\omega = \scomp \omega$.
        By \Cref{lem:full weight}, there is a full weight summand with a light factor. 
        Hence, \Cref{prop:degree 1 weights} implies that $\wt(\omega) \le 1 + 2c$.

        Inductively, assume the result for all $1 \le j < J$.
        If $\omega \in \CS{2}{J}$ has weight maximized only by terms of $\scomp\omega$ then \Cref{lem:full weight} again implies that $\wt(\omega) \le 1 + 2c < J + 2c$.
        On the other hand, if $\wt(\omega)$ is maximized by a summand of $\ncomp\omega$, then such a summand is of the form $\basis{t} \wedge \basis{y} \in \CS{1}{i_1} \wedge \CS{2}{i_2}$ for which $i_1 + i_2 \le J$ by \Cref{lem:cautious step degree 1}.
        Inductively then $\wt(\Omega) = \wt(\basis{t} \wedge \basis{y}) \le i_1 + (i_2 + 2c) \le J + 2c$.
    \end{proof}
\end{lemma}

\begin{corollary}\label{cor:degree 2 weights}
    If $\omega \in \pindual{2}$ then $\wt(\omega) \le 3c$.
\end{corollary}

Unfortunately, it is not clear that this process can be extended to degrees higher than $2$.
The following conjecture seems likely though:
\begin{conjecture}\label{conj:simple weights}
    Suppose the weight of $\omega \in \CS{n}{J}$ is maximized only by summands of $\scomp \omega$ (e.g., $J = 1$).
    Then some weight-maximizing summand of $d\omega$ has a light factor.
\end{conjecture}

If this holds, then \autoref{conj:better nilpotent bounds} holds as well.
We expect that this sort of analysis of the simple component can be done more carefully to lead to sharper upper bounds in each degree than the ones that follow.
We provide such an analysis in the following section for the situation when $c = 1$ (i.e. when the space is simple). 
For now we have the following general bounds:

\begin{proposition}\label{prop:nilpotent bounds}
    Let $\omega \in \pindual{N}$ for some $N \ge 2$ where $\pidual$ is $c$-step nilpotent up through degree $N$.
    If $\wt(\omega)$ is maximized by summands of $\scomp \omega$ then $\wt(\omega) \le N(4c-1) - 3(2c - 1)$.
    If $\omega \in \CS{N}{J}$, then $\wt(\omega) \le N(4c-1) - 3(2c-1) + (J - 1)$.
    In particular, $\wt(\omega) \le N(4c-1) - 5c + 2$.
    \begin{proof}
        By \Cref{cor:degree 2 weights} this is true for $N = 2$.
        Inductively, assume the result for all $2 \le n < N$.
        Suppose first $\omega \in \pindual{N}$ has weight maximized by a summand $\Omega = \basis{x}_1 \wedge \cdots \wedge \basis{x}_{k} \in \CS{n_1}{i_1} \wedge \cdots \wedge \CS{n_k}{i_k}$ of $\scomp\omega$.
        Let $r = \#\set{j \sothat n_j = 1}$, and write
        \[
            \wt(\Omega) 
            = \sum_{j} \wt(\basis{x}_{j}) 
            = \left(\sum_{\set{j \sothat n_j = 1}} \wt(\basis{x}_{j})\right) + \left(\sum_{\set{j \sothat n_j \neq 1}} \wt(\basis{x}_{j})\right).
        \]
        Therefore, since $\sum{n_j} = N+1$ one can check that
        \begin{align*}
            \wt(\Omega) 
                & \le \sum_{\set{j \sothat n_j = 1}}cn_j + \sum_{\set{j \sothat n_j \neq 1}} \left(n_j(4c-1) - 5c + 2\right) \\
                & = N(4c-1) + (4c - 1) - r(3c-1) - (k-r)(5c - 2).
        \end{align*}
        First, suppose $r = 0$.
        Since $k \ge 2$ and $5c - 2 > 0$ we have
        \begin{align*}
            \wt(\Omega) 
                & \le N(4c-1) + (4c - 1) - 2(5c - 2) \\
                & = N(4c-1) - 3(2c - 1).
        \end{align*}
        Next, suppose $r = 1$.
        Since $\Omega$ is a summand of $\scomp\omega$ we must have $k \ge 3$.
        Hence, since $c \ge 1$ we have
        \begin{align*}
            \wt(\Omega) 
                & \le N(4c+1) + (4c-1) - r(3c-1) - (k-r)(5c-2) \\
                & \le N(4c+1) - 8c + 3 \\
                & < N(4c+1) - 3(2c - 1).
        \end{align*}
        Now suppose $2 \le r < k$. 
        Then
        \begin{align*}
            \wt(\Omega) 
                & \le N(4c+1) + (4c-1) - r(3c-1) - (k-r)(5c-2) \\
                & \le N(4c+1) + (4c-1) - (6c-2) - (5c-2) \\
                & < N(4c+1) - 3(2c - 1).
        \end{align*}
        Finally if $r = k$ we must in fact have $r \ge 4$ since $N \ge 3$ so that
        \begin{align*}
            \wt(\Omega) 
                & \le N(4c+1) + (4c-1) - r(3c-1) \\
                & \le N(4c+1) - 8c - 3 \\
                & < N(4c+1) - 3(2c - 1).
        \end{align*}
        This completes the proof that if $\wt(\omega)$ is maximized only by terms of $\scomp\omega$ then $\wt(\omega) \le N(4c+1) - 3(2c - 1)$.
        In particular, this applies to $\omega \in \CS{N}{1}$.

        Inductively on $1 \le j < J$ (the result for $j = 1$ being part of what was just established) assume that if $\basis{x} \in \CS{N}{j}$ then $\wt(\basis{x}) \le N(4c+1) - 3(2c - 1) + (j-1)$.
        Let $\omega \in \CS{N}{J}$, and assume $\wt(\omega)$ is maximized by a component $\Omega = \basis{t} \wedge \basis{x} \in \CS{1}{a} \wedge \CS{N}{b}$ of $\ncomp \omega$.
        By \Cref{prop:degree 1 weights} we have $a + b \le J$ so that 
        \begin{align*}
            \wt(\omega) 
                & \le \wt(\basis{t}) + \wt(\basis{x}) \\
                & \le a + N(4c+1) - 3(2c - 1) + (b-1) \tag{induction} \\
                & \le N(4c + 1) - 3(2c-1) + J - 1 \tag{$a + b \le J$}
        \end{align*}
        as required.
        When $J = c$, the final result is as claimed.
    \end{proof}
\end{proposition}

We thus have the following, which consequently yields \Cref{cor:nilpotent isolip upper bounds}.

\begin{proof}[Proof of \Cref{prop:nilpotent lipschitz bounds}]
    By \Cref{prop:quantitative relative obstruction}, $\|\Phi\res{\pindual{n}}\|_{op}$ has an upper bound estimated in terms of $\|\Phi\res{d\pindual{n}}\|_{op} + O(L^n)$, and in particular we can see that for $\omega \in \pindual{n}$ we have $\|\Phi(v)\|_{op} \le CL^{\wt(\omega)}$ for some constant $C$ depending only on $n$, on $Y$, on the norms of $\pindual{n}$, and the algebraic structure of the differential.
    But $\omega \in \pindual{n}$ implies that $\wt(\omega) \le n(4c - 1) - 5c + 2$, which implies the result.
\end{proof}

\subsection{Bounds on Simple Spaces}\label{ssec:simple bounds}
We continue using our prior notation from \Cref{ssec:filtrations}.
As an immediate corollary of \Cref{prop:nilpotent bounds} we have the following naive observation:

\begin{corollary}
    If $\minimal{Y}$ is simple, then $\wt(\omega) \le 3(n-1)$ for all indecomposables $\omega \in \pindual{n}$.
\end{corollary}

However, we will see in the rest of this section that this upper bound can be improved.
We expect these methods to generalize to higher step actions as well, but the proofs and statements even in the current context are cumbersome and inelegant. 
Hence, we expect other methods exist for refining these bounds.

Now, we refine the above arguments in the event that $c = 1$ (or equivalently, in the event that $\ncomp = 0$ and $d = \scomp$).
In particular, this means that we are assuming our space is simple.
As in the above section, to prove \Cref{prop:simple lipschitz bounds} it is enough to show that the weight of any indecomposable in $\pindual{n}$ is bounded by $2n-1$. 
Therefore, this section is dedicated to the proof of \Cref{lem:simple bounds}.

\begin{definition}
    An indecomposable $\omega \in \pindual{n}$ is \term{overweight} if $\wt(\omega) \ge 2n-1$. 
\end{definition}

Note that in the setting of a simple space $\pindual{n} = \pindualupto{n}{1}$ for all $n$.
Many of the arguments made in this next lemma are closely related to those of the claims of \Cref{lem:full weight}, but with some modifications to accommodate for the fact that we are working with $c = 1$.

\begin{lemma}\label{lem:simple bounds}
    Suppose $\wedge V$ is simple, and let $\omega \in \pindual{n}$ be indecomposable. 
    If $\omega$ is overweight, then $\wt(\omega) = 2n - 1$, and the weight maximizing summands of $d\omega$ are of the form $\basis{p} \wedge \basis{q}$, or $\basis{a} \wedge \basis{b} \wedge \basis{c}$.
    In particular, in the quadratic case, at most one of $\basis{p}$ or $\basis{q}$ is overweight.
    In the cubic case, each factor is overweight.
    \begin{proof}
        The result is trivial for indecomposables of $\pindual{1}$.
        If $\omega \in \pindual{2}$, then either $d\omega = 0$ or every summand of $d\omega$ belongs to $\pindual{1} \wedge \pindual{1} \wedge \pindual{1}$. 
        As every element of $\pindual{1}$ is (by definition) overweight, the result holds practically trivially for $n = 2$ as well.

        Let $n > 2$ and assume inductively the result is true for all elements of degree less than $n$.
        Let $\omega \in \pindual{n}$ be overweight, so that $\wt(\omega) \ge 2n - 1$.
        In a similar approach to many of the previous propositions, the primary tool for analysis is the fact that $d^2\omega = 0$.

        Suppose to the contrary that $d\omega$ contains a nonzero summand $\basis{x} \wedge \basis{y}$ for which both $\basis{x}$ and $\basis{y}$ are overweight.
        Among all quadratic terms with overweight factor $\basis{x}$, take $\basis{x}$ to be minimal degree. 
        That is, if $\basis{x}' \wedge \basis{y}'$ is another summand with $\basis{x}'$ overweight, then $|\basis{x}'| \ge |\basis{x}|$.
        Since $|\basis{y}| + |\basis{x}| = n+1$, this implies that $|\basis{y}'| \le |\basis{y}|$ as well, so that the degree of $\basis{y}$ is maximal among such overweight factors.
        Moreover, we may assume without loss of generality there are no other summands $\basis{x}' \wedge \basis{y}$ of $d\omega$ for which $d\basis{x}'$ contains $\xi$ as a summand.
        Otherwise, replace $\basis{x}$ with $\basis{x} + \basis{x}'$.

        Since $\basis{x}$ is overweight, by induction $d\basis{x}$ contains an overweight summand of the form $\xi_{2} = \basis{p} \wedge \basis{q}$ with $\basis{p}$ overweight and $\basis{q}$ not overweight, or of the form $\xi_{3} = \basis{a} \wedge \basis{b} \wedge \basis{c}$ where each factor is an overweight basic indecomposable.
        We assume that the coefficient of this summand is $1$ by rescaling $\omega$ if necessary, and we will write $\xi \in \set{\xi_2,\xi_3}$ when dealing with both cases simultaneously.

        We seek a contradiction to the fact that $d^2\omega = 0$.
        Consider the summand $\xi \wedge \basis{y}$ of $(d\basis{x}) \wedge \basis{y}$.
        By \Cref{lem:internal} either $\xi \wedge \basis{y} = 0$ or $\xi \wedge \basis{y}$ is canceled externally.
        If $\xi \wedge \basis{y} = 0$ then $\basis{y}$ must be a factor of $\xi$. 
        But each factor of $\xi$ has degree strictly less than $|\basis{x}| \le |\basis{y}|$, a contradiction.
        Therefore, the only option is external cancellation.

        To this end, suppose $\basis{x}' \wedge \basis{y}'$ is a quadratic summand of $d\omega$ for which $d(\basis{x}' \wedge \basis{y}')$ contains a cancelling summand of $\xi \wedge \basis{y}$.
        Suppose without loss of generality that $|\basis{x}'| \le |\basis{y}'|$.
        If $(d\basis{x}') \wedge \basis{y}'$ cancels with $\xi \wedge \basis{y}$, then by consideration of degrees we must have $\basis{y}' = \basis{y}$.
        On the other hand, suppose $\basis{x}' \wedge d\basis{y}'$ cancels with $\xi \wedge \basis{y}$.
        If $\basis{x}'$ is overweight, then either $\basis{x}'$ agrees with $\basis{y}$ (in which case, we were actually considering internal cancellation), or else $\basis{x}'$ is an overweight factor of $\xi$ which contradicts the minimality assumptions on $|\basis{x}|$.
        On the other hand, if $\basis{x}'$ is not overweight, then $\xi = \xi_{2}$ must have $\basis{x}'$ as a factor.
        That is, $\xi_2 = \basis{p} \wedge \basis{q} = \basis{p} \wedge \basis{x}'$.
        But then $\inner{\basis{x}' \wedge d\basis{y}', \xi_2 \wedge \basis{y}} = \inner{\basis{x}' \wedge d\basis{y}', \basis{p} \wedge \basis{x}' \wedge \basis{y}} \neq 0$ implies that $\basis{p} \wedge \basis{y}$ is a nonzero summand of $d\basis{y}'$.
        Since $\basis{y}$ and $\basis{p}$ are both overweight, this contradicts the inductive assumption.

        Finally, suppose there is a cubic summand $\basis{x}' \wedge \basis{y}' \wedge \basis{z}'$ in $d\omega$ whose differential contains (without loss of generality) a cancelling summand $\basis{x}' \wedge \basis{y}' \wedge \basis{p}' \wedge \basis{q}'$ which was produced by $\basis{x}' \wedge \basis{y}' \wedge d\basis{z}'$.
        Inductively, since this term came from a quadratic differential on $\basis{z}'$, we may assume that $\basis{q}'$ is not overweight.
        If $\xi = \xi_2$, then $\inner{\basis{x}' \wedge \basis{y}' \wedge \basis{p}' \wedge \basis{q}',\xi_2 \wedge \basis{y}} = 0$ by consideration of word lengths, and so it is enough to consider cancellation when $\xi = \xi_3$.
        But in this case, note that each factor of $\basis{a} \wedge \basis{b} \wedge \basis{c} \wedge \basis{y}$ is overweight by assumption.
        This implies that each factor of $\basis{x}' \wedge \basis{y}' \wedge \basis{p}' \wedge \basis{q}'$ is overweight as well, contradicting the weight of $\basis{q}'$.     
        
        Therefore, if the weight of $\omega$ is maximized by a summand of the form $\basis{x} \wedge \basis{y}$, at most one factor is overweight, and so without loss of generality $\wt(\basis{x}) \le 2|\basis{x}| - 2$ and $\wt(\basis{y}) \le 2|\basis{y}| - 1$. In particular, 
        \begin{equation*}
            \wt(\omega) = \wt(\basis{x}) + \wt(\basis{y}) \le 2(|\basis{x}| + |\basis{y}|) - 3 = 2(|\omega| + 1) - 3 = 2|\omega| - 1.
        \end{equation*}
        On the other hand, if the weight of $\omega$ is maximized by a summand of the form $\basis{a} \wedge \basis{b} \wedge \basis{c}$, each factor may be overweight, which inductively implies that 
        \begin{equation*}
            \wt(\omega) \le (2|\basis{a}| - 1) + (2|\basis{b}| - 1) + (2|\basis{c}| - 1) = 2(|\omega| + 1) - 3 = 2|\omega| - 1,
        \end{equation*}
        and $\omega$ is overweight when $\wt(\omega)$ is maximized in one of these ways.
    \end{proof}
\end{lemma}

We record this improvement as a separate result: 
\begin{corollary}
    If $\minimal{Y} = \wedge \pidual$ is simple, then $\wt(\omega) \le 2n - 1$ for all $\omega \in \pindual{n}$. 
\end{corollary}

In particular the bounds in \Cref{prop:simple lipschitz bounds} hold, and thus so does \Cref{cor:simple isolip upper bounds}.

\subsection{Bounds on Coformal Spaces}\label{ssec:coformal bounds}
We conclude this section by reducing the upper bounds for coformal spaces.
By \cite[Prop. 3.3]{neisendorfer1978formal}, we may take as our definition that a space $Y$ is \textbf{coformal} if $(\minimal{Y},d)$ is such that $d = d_1 : V \to V \wedge V$ is quadratic.

\begin{lemma}\label{lem:coformal weights}
    If $d$ is quadratic and $\omega \in \CS{n}{J}$, then 
    \[
        d\omega 
        \in \Dirsum_{\substack{n_1 \le n_2 \\ n_1 + n_2 = n+1}}\left(\Dirsum_{i_1 + i_2 \le J + c} \CS{n_1}{i_1} \wedge \CS{n_2}{i_2}\right).
    \]
    \begin{proof}
        Note that \textit{a priori}, we know only that $i_1,i_2 \le J$ and $i_1,i_2 \le c$ in the above expression.
        For instance, since $d$ is quadratic, we can write 
        \[
            d\omega = \sum_{\substack{n_1 \le n_2 \\ n_1 + n_2 = n+1}} \left(\sum_{i_1,i_2 \le c} \Omega_{(i_1,i_2,n_1,n_2)}\right) \in \Dirsum_{\substack{n_1 \le n_2 \\ n_1 + n_2 = n+1}}\left(\Dirsum_{i_1, i_2 \le c} \CS{n_1}{i_1} \wedge \CS{n_2}{i_2}\right).
        \]
        The result is obvious for $n \le 2$.
        Inductively assume the result is true for all $k < n$ for some $n > 2$.
        By \Cref{lem:cautious step higher degrees}, if $\omega \in \CS{n}{J}$, then $\ncomp \omega \in \Dirsum_{i + j \le J}\CS{1}{i} \wedge \CS{n}{j}$ so that it is enough to consider terms of $\scomp\omega$.
        To this end, let $(N_1,N_2)$ be such that $\Omega_{(i_1,i_2,N_1,N_2)} \neq 0$ for some $(i_1,i_2)$.
        In particular, assume $1 < N_1 \le N_2 < n$.
        Moreover, define 
        \[
            \ell \defeq \max\set{i_1 + i_2 \sothat \Omega_{(i_1,i_2,N_1,N_2)} \neq 0}
        \]
        and fix 
        \[
            (I_1,I_2) \defeq \max_{\len{2}} \set{(i_1,i_2) \sothat i_1 + i_2 = \ell,\; \Omega_{(i_1,i_2,N_1,N_2)} \neq 0}.
        \]
        Let $\Omega = \basis{x} \wedge \basis{y} \in \CS{N_1}{I_1} \wedge \CS{N_2}{I_2}$ be a nonzero basic summand of $d\omega$, and define $\xi \defeq (\abcomp \basis{x}) \wedge \basis{y} \in \CS{1}{1} \wedge \CS{N_1}{I_1-1} \wedge \CS{N_2}{I_2}$. 
        
        By \Cref{lem:trivial or external}, either $\xi = 0$ or $\xi$ is cancelled externally.
        However, injectivity of $\abcomp$ and the fact that $(I_1-1,N_1) \lessn{2} (I_2,N_2)$ implies nearly immediately that $\xi \neq 0$.
        Hence, there is some external summand $\Omega' = \basis{x}' \wedge \basis{y}' \in \CS{N_1'}{I_1'}\wedge\CS{N_2'}{I_2'}$ of $d \omega$ for which $d\Omega'$ has a summand $\xi'$ that cancels with $\xi$.
        Here as before we order degrees so that $N_1' \le N_2'$.

        Since $\xi'$ is a summand of $d\Omega'$ and since $\xi'$ cancels with $\xi$ we must have
        \[
            \xi' \in \left[\left(\im{d\res{\CS{N_1'}{I_1'}}} \wedge \CS{N_2'}{I_2'}\right) + \left(\CS{N_1'}{I_1'} \wedge \im{d\res{\CS{N_2'}{I_2'}}}\right)\right] \cap \left[\CS{1}{1} \wedge \CS{N_1}{I_1-1}\wedge\CS{N_2}{I_2}\right].
        \]
        In particular there are $a_1,b_1,a_2,b_2$ for which
        \[
            \xi' \in \left(\CS{n_{1,1}}{a_1} \wedge \CS{n_{1,2}}{b_1} \wedge \CS{N_2'}{I_2'}\right) + \left(\CS{N_1'}{I_1'} \wedge \CS{n_{2,1}}{a_2} \wedge \CS{n_{2,2}}{b_2}\right) 
        \]
        where $n_{1,1} + n_{1,2} = N_1' + 1$ and $n_{2,1} + n_{2,2} = N_2' + 1$.
        Suppose first that 
        \begin{equation*}
            \xi \in \left(\CS{n_{1,1}}{a_1} \wedge \CS{n_{1,2}}{b_1} \wedge \CS{N_2'}{I_2'}\right) \cap \left(\CS{1}{1} \wedge \CS{N_1}{I_1-1} \wedge \CS{N_2}{I_2}\right). 
            \tag{A1}
            \label{assumption 1}
        \end{equation*}

        \begin{claim}
            Assuming \eqref{assumption 1}, then $\Omega' \in (I_1,I_2,N_1,N_2)$. 
            \begin{proofofclaim}
                By considering degrees, and the assumption that $N_1' \le N_2'$, we see that we must have $N_2' = N_2$. 
                Moreover, without loss of generality we have $n_{1,1} = 1$ and $N_1' = n_{1,2} = N_1$, since $1 < N_1$.
                In particular, this implies that $a_1 = 1$.
                Furthermore, since
                \[
                    \xi' \in \im{d\res{\CS{N_1}{I_1'}}} \wedge \CS{N_2}{I_2'},
                \] 
                this implies that in fact $\xi'$ came from the nilpotent component of the differential.
                By definition we must then have $b_1 \le I_1'-1$.

                Suppose first that $N_1 < N_2$.
                Then $I_2' = I_2$ and $b_1 = I_1 - 1$.
                But then $I_1 = 1 + b_1 \le I_1'$.
                Thus since $(N_1',N_2') = (N_1,N_2)$, the maximality of $(I_1,I_2)$ implies that in fact $(I_1',I_2') = (I_1,I_2)$.

                If $N_1 = N_2$, there are two possibilities.
                If $b_1 = I_1 - 1$, then the above proof shows $(I_1',I_2',N_1',N_2') = (I_1,I_2,N_1,N_2)$.
                On the other hand, if $b_1 = I_2$ then $I_1' \ge b_1 + 1 > I_2$, which contradicts the maximality of $I_2$.
            \end{proofofclaim}
        \end{claim}

        By injectivity of $\abcomp$ this implies that $\Omega' \in \sspan{\Omega}$, which means that cancellation was internal, a contradiction.

        Suppose next that 
        \begin{equation*}
            \xi' \in 
                \left(\CS{N_1'}{I_1'} \wedge \CS{n_{2,1}}{a_2} \wedge \CS{n_{2,2}}{b_2}\right) 
                \cap 
                \left(\CS{1}{1} \wedge \CS{N_1}{I_1-1} \wedge \CS{N_2}{I_2}\right).
                \tag{A2}
                \label{assumption 2}
        \end{equation*}
        Note that by induction we have $a_2 + b_2 \le I_2' + c$.
        In the previous case, we could essentially ignore the situation $N_1' = 1$. 
        Now however we have either $N_1' = 1$ or (without loss of generality) $n_{2,1} = 1$.
        We consider these in the reverse order:
        
        \begin{claim}
            Assuming \eqref{assumption 2}, if $n_{2,1} = 1$ then $(I_2',I_1',N_1',N_2') = (I_1,I_2,N_1,N_2)$.
            \begin{proofofclaim}
                In this situation, we see by the same arguments as in the previous claim, that $N_1' = N_1$ and $N_2' = n_{2,2} = N_2$, and that furthermore $a_2 = 1$ so that $b_2 \le I_2' - 1$.

                In the event that $N_1 < N_2$ we have therefore that $b_2 = I_2$ which means $I_2' \ge 1 + b_2 > I_2$, contradicting the maximality of $I_2$.

                If $N_1 = N_2$ there are two possibilities. If $b_2 = I_2$ then as above we contradict the maximailty of $I_2$.
                On the other hand, if $b_2 = I_1 - 1$ then $I_2' \ge b_2 + 1 = I_1$.
                Hence (noting the switch of position of $I_2'$ and $I_1'$ in the lefthand term) we must have $(I_2',I_1',N_1',N_2') = (I_1,I_2,N_1,N_2)$.
            \end{proofofclaim}
        \end{claim}
        As before, in the above setting injectivity of $\abcomp$ implies then that $\Omega' \in \sspan{\Omega}$.

        \begin{claim}
            Assuming \eqref{assumption 2}, if $N_1' = 1$, then $I_1 + I_2 \le J + c$.
            \begin{proofofclaim}
                In this case, because $1 < N_1$, we know that $N_1' \neq N_1$.
                Therefore, without loss of generality, we have $n_{2,1} = N_1$ and $n_{2,2} = N_2$,
                which implies that $a_2 = I_1 - 1$ and $b_2 = I_2$.
                Moreover we must have $I_1' = 1$, and thus
                $\Omega'$ is a term in $\ncomp\omega$.
                In particular, $I_2' = J - 1$.
                By induction, therefore, we see that $(I_1 - 1) + I_2 =     a_2 + b_2 \le (J-1) + c$.
            \end{proofofclaim}
        \end{claim}
        Hence $\ell = I_1 + I_2 \le J + c$ as required.
    \end{proof}
\end{lemma}

The above result indicates that we might reasonably hope to refine the bounds on the weights of indecomposables in the event that $d$ is quadratic.
Indeed, we can.

\begin{lemma}
    Let $\omega \in \CS{n}{J}$ for some $c$-step nilpotent DGA $(\wedge \pidual, d)$ with quadratic differential.
    Then for $n \ge 2$ we have
    \[ 
        \wt(\omega) \le (c+1)(n-1) + J - c.
    \]
    In particular any $\omega \in \pindual{n}$ has $\wt(\omega) \le (c+1)(n-1)$.
    \begin{proof}
        For $n = 2$ and $J = 1$ the result is trivial because $d\res{\CS{2}{1}} = 0$.
        Note that for $\basis{t} \in \CS{1}{J}$ we have $\wt(\basis{t}) = J$.
        Since $d$ is quadratic we have $d\res{\CS{2}{J}} = \ncomp\res{\CS{2}{J}}$, so that \Cref{lem:cautious step higher degrees} implies inductively that if $\omega \in \CS{2}{J}$ then $\wt(\omega) \le J + 1$ which is the desired form when $n = 2$.
        In particular, $\wt(\omega) \le c + 1$ for all $\omega \in \pindual{2}$.
        
        Inductively assume the result holds for all $(j,k) \lessn{2} (J,n)$ for some $(J,n)$.
        Suppose that $\basis{x} \wedge \basis{y} \in \CS{n_1}{i_1} \wedge \CS{n_2}{i_2}$ is a summand of $d\omega$.
        By \Cref{lem:coformal weights} we have $i_1 + i_2 \le J+c$.
        By induction, we have 
        \begin{align*}
            \wt(\basis{x} \wedge \basis{y}) 
                & \le \Big((c+1)(n_1-1) + i_1 - c \Big)+ \Big((c+1)(n_2-1) + i_2 - c\Big) \\
                & \le (c + 1)(n - 1) + J - c.
        \end{align*}
        Since $\basis{x} \wedge \basis{y}$ was an arbitrary summand of $d\omega$ we have $\wt(\omega) \le (c+1)(n-1) + J - c$ as desired.
    \end{proof}
\end{lemma}

\begin{corollary}\label{cor:coformal weights}
    If $Y$ is coformal then $\wt(\omega) \le (c+1)(n-1)$ for all $\omega \in \pindual{n}$.
\end{corollary}

In particular the bounds in \Cref{prop:coformal lipschitz bounds} hold, and thus so does \Cref{cor:coformal isolip upper bounds}.

\section{High Volume Nullhomotopies in Coformal Spaces}\label{sec:coformal-example}
We conclude by demonstrating the (near) sharpness of the bounds in \Cref{cor:coformal isolip upper bounds}.
This section may be viewed as motivation for further study of the general bounds of \Cref{thm:nilpotent upper bounds}.

We take as inspiration a family of maps into simply connected examples which must necessarily have high-volume nullhomotopies, as constructed in \cite[\textsection 7]{Chambers_MW_2018_quantitative_nullhomotopy}.
There, one constructs a space $X_n$ by attaching $(n + 1)$-cells to $\SS^2 \wedgesum \SS^2$ in order to kill $\pi_n(\SS^2 \wedgesum \SS^2) \tensor \QQ$. By taking iterated Whitehead brackets of $L$-Lipschitz maps one can find a nullhomotopic $L$-Lipschitz map $S^n \to X_n$ for which every nullhomotopy has degree $\Omega(L^{2n-2})$ over some $(n + 1)$-cell.
We seek a generalization of this example to nilpotent spaces, where the volume fundamentally relies on the nilpotence class of the space in question.
\Cref{cor:coformal isolip upper bounds} shows that the family we construct in \Cref{thm:coformal example} is (nearly) optimal among coformal spaces.

By \cite[Proposition 3.3]{neisendorfer1978formal}, we may view a coformal space as one whose (rational) higher order Whitehead brackets contain zero.
(See \cite{andrews1978sullivan} for more on the connection between higher order Whitehead brackets and minimal models.)
In this subsection we build $L$-Lipschitz maps into a $c$-step space, and maps into that space through (standard) Whitehead brackets.
In \Cref{thm:coformal example} we show that the resulting map has sufficiently high degree to require nullhomotopies of volume $\Omega(L^{(c+1)(n-1)})$.

Define $\wedgeofs \defeq \bigvee_{i=1}^{c} \SS^{2}_{i}$ to be a boquet of $c$-many $2$-spheres. Note that $H_2(\wedgeofs) \aisom \ZZ^{c}$.
Define $\map{\splitt}{\wedgeofs}{\wedgeofs}$ to be any self-map inducing the $(c \times c)$-matrix 
\begin{equation}\label{eq:c-by-c action}
    \begin{bmatrix}
        1       & 0     & 0      & \cdots    & 0 & 0 \\
        1       & 1     & 0      & \cdots    & 0 & 0 \\
        0       & 1     & 1      & \cdots    & 0 & 0 \\
        \vdots  &\ddots & \ddots & \ddots & \vdots & \vdots \\
        0       &0      &0        & \cdots & 1  & 0 \\
        0       &0      &0      & \cdots    & 1    & 1
        \end{bmatrix}
    \in \operatorname{Aut}(\ZZ^{c}) \aisom \operatorname{Aut}(H_2(\wedgeofs)).
\end{equation}
For each $q \in \set{A,B}$ let $\torusofsplit{q}$ be a copy of the mapping torus of $\splitt$:
\begin{equation*}
    \torusofsplit{q} \defeq \frac{\l(\wedgeofs \times I\r) \cup_{\splitt} \l(\wedgeofs \times \set{1}\r)}{(x,0) \sim (x,1)}.
\end{equation*}
Given $* \in \wedgeofs$ common to each copy of $\SS^{2}$ define $\gamma = \set{*} \times \SS^{1} \sse \torusofsplit{q}$ to be the circle traced out by the basepoint.
Denote $\wedgeofs[2c] = (\wedgeofs[c])_{A} \wedgesum (\wedgeofs[c])_{B} = (\bigvee_{i=1}^{c} \SS^{2}_{A,i}) \wedgesum (\bigvee_{i=1}^{c} \SS^{2}_{B,i})$.
Define $\gluedt$ to be the mapping torus of ${\splitt_{A}} \wedgesum {\splitt_{B}} : \wedgeofs[2c] \to \wedgeofs[2c]$ under the obvious choices of domain and codomain for each copy of $\splitt$.
Note that $\gluedt = \torusofsplit{A} \cup_{\gamma} \torusofsplit{B}$. 
It is straightforward to see that $\pi_1(\gluedt) \aisom \ZZ$ is generated by the natural inclusion $\map{t}{\SS^{1}}{\gamma} \into \gluedt$.

\begin{lemma}\label{prop:action of t}
    Let $\alpha_j \in \pi_{2}(\torusofsplit{A})$ be the homotopy class of the map $\SS^{2} \to \SS^{2}_j \into \torusofsplit{A}$. 
    Then $t \bullet \alpha_j = \alpha_{j} + \alpha_{j+1}$, and $[t,\alpha_{j}] = \alpha_{j+1}$ for $j < c$ while $t \bullet \alpha_c = \alpha_{c}$ and $[t,\alpha_{c}] = 0$.
    Identical results hold for $\beta_j$ generators of $\pi_2(\torusofsplit{B})$.
    \begin{proof}
        By definition the generator $t \in \pi_1(\gluedt)$ acts on $\pi_2(\torusofsplit{A}) \tensor \QQ \sse \pi_2(\gluedt) \tensor \QQ$ under the basis $\set{\alpha_1,\alpha_{2},\dots,\alpha_{c}}$ by the $(c \times c)$ matrix in \eqref{eq:c-by-c action} above.
    \end{proof}
\end{lemma}

Inductively use Whitehead brackets to define
\begin{equation*}
    \highdeg{k}{j} \defeq
    \begin{cases}
        \beta_1                        & k = 2, j = 1, \\
        [\alpha_1,\highdeg{k-1}{c}]    & k > 2, j = 1, \\
        [t,\highdeg{k}{j-1}]           & j > 1.
    \end{cases}
\end{equation*}
Conflating maps with the homotopy elements they represent, observe that $\highdeg{k}{j} \in \pi_{k}(\gluedt)$.

\begin{lemma}\label{lem:jacobi}
    For all $k > 2$ and $j < c$ we have $[t,[\alpha_j,\highdeg{k-1}{c}]] = [\alpha_{j+1},\highdeg{k-1}{c}]$ and for all $k \ge 2$ we have $[t,\highdeg{k}{c}] = 0$.
    \begin{proof}
        The fact $[t,\highdeg{2}{c}] = 0$ is true by \Cref{prop:action of t} (applied, of course, to $\beta_{j}$ instead of $\alpha_{j}$).
        Let $k > 2$ and assume for induction that the result holds for all lower indices.
        By the Jacobi identity we have
        \begin{align*}
            [t,[\alpha_j,\highdeg{k-1}{c}]] 
                & = -[\alpha_{j},[\highdeg{k-1}{c},t]] - [\highdeg{k-1}{c},[t,\alpha_j]] \\
                & = -[\alpha_{j},0]
                    - [\highdeg{k-1}{c},[t,\alpha_j]] \tag{induction} \\
                & = [[t,\alpha_j],\highdeg{k-1}{c}]. \tag{antisymmetry} 
        \end{align*}
        In the event $j < c$, \Cref{prop:action of t} implies then that $[t,[\alpha_j,\highdeg{k-1}{c}]]  = [\alpha_{j+1},\highdeg{k-1}{c}]$ while if $j = c$, the same result indicates that $[t,[\alpha_j,\highdeg{k-1}{c}]] = [0,\highdeg{k-1}{c}] = 0$ as required.
    \end{proof}
\end{lemma}

\begin{corollary}\label{cor:jacobi}
    $\highdeg{k}{j} = [\alpha_{j},\highdeg{k-1}{c}]$ for all $k$ and $j$.
    In particular $\highdeg{k}{c} = [\alpha_{c},\highdeg{k-1}{c}]$.
    \begin{proof}
        Inductively by the previous lemma, $\highdeg{k}{j} = [t,\highdeg{k}{j-1}] = [t,[\alpha_{j-1},\highdeg{k-1}{c}]] = [\alpha_{j},\highdeg{k-1}{c}]$, where the base case $j = 1$ is true by definition.
    \end{proof}
\end{corollary}

\begin{lemma}\label{lem:not nullhomotopic}
    For all $k$, one has $\highdeg{k}{c} \neq 0$.
    \newcommand{\mcyl}{\ms{M}}
    \begin{proof}
        Let $\mcyl = \mcyl_{\splitt_{A} \wedgesum \splitt_{B}}$ denote the mapping cylinder of $\splitt_{A} \wedgesum \splitt_{B}$.
        The universal cover $\tilde{\gluedt}$ of $\gluedt$ may be viewed as $\ZZ$ copies of $\mcyl$ identified in the obvious ways.
        Let $\mcyl_{p}$ denote the copies of $\mcyl$ in $\tilde{\gluedt}$ over the interval $[-p,p]$ so that $\gluedt = \bigcup_{p \in \ZZ} \mcyl_{p}$.
        
        Let $\tilde{\alpha_{c}} : \SS^{2}_{A,c} \into \wedgeofs[2c] \times \set{0} \into \tilde{\gluedt}$ lift $\alpha_{c} : \SS^{2}_{A,c} \into \gluedt$ in the obvious way.
        Similarly denote $\tilde{\beta_{c}}$ and $\tilde{\zeta_{k,c}}$.
        It is not hard to see that $\mcyl_{p}$ deformation retracts to $\wedgeofs[2c] \times \set{p}$ through standard deformation retractions of each copy of $\mcyl$ so that $\pi_{n}(\mcyl_{p}) \isom \pi_{n}(\wedgeofs[2c])$ for all $p$.
        Moreover, because $t \bullet \alpha_{c} = \alpha_{c}$ (and similarly for $\beta_{c}$) we see that this retract sends $(\SS^{2}_{A,c} \wedgesum \SS^{2}_{B,c}) \times \set{0}$ identically onto $(\SS^{2}_{A,c} \wedgesum \SS^{2}_{B,c}) \times \set{p}$.
        Then the retract identifies the class of $\tilde{\alpha_{c}}$ with that of $\Id_{A,c} : \SS^{2}_{A,c} \into \wedgeofs[2c]$.
        The same is true for $\tilde{\beta_{c}}$ under the obvious changes.
        That is to say, the isomorphism $\pi_{2}(\mcyl_{p}) \aisom \pi_{2}(\wedgeofs[2c])$ induced by the retraction will identify $\tilde{\alpha_{c}}$ and $\tilde{\beta_{c}}$ with $\Id_{A,c}$ and $\Id_{B,c}$ respectively.
        By induction, since $\zeta_{k,c} = [\alpha_{c},\zeta_{k-1,c}]$ and since $t \bullet \zeta_{k,c} = \zeta_{k,c}$, we see that the retraction of $\mcyl_{p}$ identifies the class of $\tilde{\zeta_{k,c}}$ in $\pi_{k}(\mcyl_{p})$ with the nonzero class $[\Id_{A,c},[\Id_{A,c},\cdots[\Id_{A,c},\Id_{B,c}]\cdots]] \in \pi_{k}(\wedgeofs[2c])$.
        Hence for any $p$, the map $\tilde{\zeta_{k,c}}$ is not nullhomotopic in $\mcyl_{p}$.

        If to the contrary $\zeta_{k,c}$ were nullhomotopic by some $H : \SS^{k} \times I \to \gluedt$, there would be a lift $\tilde{H} : \SS^{k} \times I \to \tilde{\gluedt}$ nullhomotopoing $\zeta_{k,c}$ inside $\tilde{\gluedt}$.
        By compactness there is some $p$ for which this nullhomotopy is contained entirely within $\mcyl_{p}$, a contradiction.
    \end{proof}
\end{lemma}

\begin{lemma}\label{lem:whitehead lip maps}
    Let $Y$ be a space for which $\map{f}{\SS^{n_1}}{Y}$ and $\map{g}{\SS^{n_2}}{Y}$ are $O(L)$-Lipschitz maps.
    Then $[f,g] \in \pi_{n_1 + n_2 -1}(Y)$ has an $O(L)$-Lipschitz representative.
    \begin{proof}
        The map $[f,g]$ is formed by the composition
        \begin{equation*}
            \SS^{n_1 + n_2 - 1} \to \SS^{n_1} \wedgesum \SS^{n_2} \rto{f \wedgesum g} Y.
        \end{equation*}
        The first map is $O(1)$-Lipschitz, and the second is $O(L)$-Lipschitz, and hence the composition is $O(L)$-Lipschitz.
    \end{proof}
\end{lemma}

\begin{proposition}
    For all $k$ there is an $O(L)$-Lipschitz representative $\highdegrep{k}{c}{L}$ of $L^{(c+1)(k-1)}\highdeg{k}{c}$.
    \begin{proof}
        Let $\alpha_{1}^{(L)}$ be the map $\SS^{2} \rto{\times L} \SS^{2} \rto{\alpha_1} \gluedt$ be the $O(L)$-Lipschitz map formed by precomposing $\alpha_1$ with an $L$-Lipschitz self-map of $\SS^{2}$ of degree $L^2$. 
        That is, $\alpha_{1}^{(L)}$ is an $L$-Lipschitz representative of $L^2\alpha_1$.
        Similarly, define $t^{(L)}$ and $\beta_{1}^{(L)}$ to be $L$-Lipschitz representatives of $Lt$ and $L^2\beta_1$ respectively.
        Inductively, we attain $L$-Lipschitz representatives $\alpha_{j}^{(L)}$ and $\beta_{j}^{(L)}$ of $L^{j+1}\alpha_{j}$ and $L^{j+1}\beta_{j}$ respectively by 
        \begin{equation*}
            \alpha_{j}^{(L)} \defeq [t^{(L)},\alpha_{j-1}^{(L)}] = [Lt,L^{j}\alpha_{j-1}] = L^{j+1}\alpha_{j}
        \end{equation*}
        and similarly for $\beta_{j}^{(L)}$.
        In particular there is an $O(L)$-Lipschitz representative $\highdegrep{2}{c}{L} \defeq \beta_{c}^{(L)}$ of $L^{c+1}\highdeg{2}{c} = L^{c+1}\beta_{2}$.
        This establishes the base case $k = 2$.

        Assume the result inductively for $k - 1$ for some $k > 2$.
        By \Cref{lem:whitehead lip maps} there is an $O(L)$-Lipschitz representative $\highdegrep{k}{c}{L}$ of
        \begin{align*}
            [\alpha_{c}^{(L)},\highdegrep{k-1}{c}{L}]
                & = [L^{c+1}\alpha_{c},L^{(c+1)(k-2)}\highdeg{k-1}{c}] \tag{induction} \\
                & = L^{(c+1) + (c+1)(k-2)}[\alpha_{c},\highdeg{k-1}{c}] \tag{homogeneity} \\
                & = L^{(c+1)(k-1)}\highdeg{k}{c}, \tag{\Cref{cor:jacobi}}
        \end{align*}
    which completes the proof.
    \end{proof}
\end{proposition}

\begin{theorem}\label{thm:coformal example}
    There is a space $Y$ which is $c$-step nilpotent and coformal up through degree $n$, together with a nullhomotopic, $L$-Lipschitz map $f : \SS^{n} \to Y$ such that any nullhomotopy of $f$ has volume $\Omega(L^{(c+1)(n-1)})$.
    \begin{proof}
        The argument here is identical to that of \cite[\textsection 7]{Chambers_MW_2018_quantitative_nullhomotopy}.
        Let $f = \highdegrep{n}{c}{L}$, and let $Y$ be the space formed by attaching $(n+1)$-cells to $\gluedt$, whose attaching maps form a basis for $\pi_{n}(\gluedt) \tensor \QQ$.
        Ensure that $Y$ is coformal through degree $n$ by attaching higher cells to kill any nontrivial higher order Whitehead products.
        In particular, fix a standard nullhomotopy $N_{n,c}$ of $\highdeg{n}{c}$ in $Y$.
        Such a nullhomotopy has degree $C \neq 0$ (in the sense of relative homology) over some $(n+1)$ cell $q$ of $Y$ because (by \Cref{lem:not nullhomotopic}), $\highdeg{n}{c}$ is not nullhomotopic in $\gluedt$.
        This yields a closed $(n+1)$-form $\omega$ on $Y$ such that $\int_{\SS^{n} \times I}N_{n,c}^*\omega = C \neq 0$.
        Let $f' = \highdeg{n}{c} \circ \phi_{n,L}$ where $\SS^{n} \to \SS^{n}$ is a self-map of degree $L^{(c+1)(n-1)}$.
        Note that $f' \homotopic f$ as both are representatives of $L^{(c+1)(n-1)}\highdeg{n}{c}$ by definition.
        Hence a nullhomotopy $N_{n,c,L}$ of $f$ may be defined by first homotoping $f$ to $f'$, followed by the standard nullhomotopy $N_{n,c}$.
        This nullhomotopy has volume bounded below by $\int_{\SS^{n} \times I}N_{n,c,L}^*\omega = L^{(c+1)(n-1)}C = \Omega(L^{(c+1)(n-1)})$.

        Now suppose that $N'$ is some other nullhomotopy of $f$. Let $\tilde{N}$ be the map $\SS^{n+1} \to Y$ formed by performing $N_{n,c,L}$ along the northern hemisphere and $N'$ along the southern hemisphere, connected along the $\SS^{n}$ equator along the map $f$ itself.
        The (rational) Heurewicz map sends $\pi_{n+1}(Y)$ to $0$ and so the total degree of $\tilde{N}$ must be $0$, and hence $\tilde{N}$ has zero degree on any $(n+1)$-cell.
        In particular, this means that $\int_{\SS^{n} \times I}(N')^*\omega = \int_{\SS^{n} \times I}N_{n,c,L}^*\omega = \Omega(L^{(c+1)(n-1)})$.
    \end{proof}
\end{theorem}

%-----------------------
% Glossary and References
%-----------------------

\nocite{*}
\bibliographystyle{alpha}
\bibliography{sources.bib}

\end{document}